\documentclass[11pt]{amsart}                    
\usepackage{amssymb}  

\input xy
\xyoption{all}
 
\theoremstyle{plain}
\newtheorem{theorem}{Theorem}[section]
\newtheorem{corollary}[theorem]{Corollary}
\newtheorem{lemma}[theorem]{Lemma}
\newtheorem{subclaim}[theorem]{Subclaim}
\newtheorem{conjecture}[theorem]{Conjecture}
\newtheorem{proposition}[theorem]{Proposition}
\newtheorem{fact}[theorem]{Fact}
\newtheorem{claim}[theorem]{Claim}
\newtheorem{question}[theorem]{Question}
\theoremstyle{definition}
\newtheorem{definition}[theorem]{Definition}
\newtheorem{remark}[theorem]{Remark}

\newtheorem{assumption}[theorem]{Assumption}
\newtheorem{notation}[theorem]{Notation}
\theoremstyle{remark}

\newcommand{\St}{\operatorname{\mathfrak{S}t}}

\newcommand{\cf}{\operatorname{cf}}
\newcommand{\otp}{\operatorname{otp}}

\newcommand{\dom}{\operatorname{dom}}
\newcommand{\Aut}{\operatorname{Aut}}
\def\rg{\operatorname{rg}}
\newcommand{\ftp}{\operatorname{tp}}
\newcommand{\tp}{\operatorname{ga-tp}}
\newcommand{\gaS}{\operatorname{ga-S}}
\newcommand{\Hanf}{\operatorname{Hanf}}
\newcommand{\LS}{\operatorname{LS}}

\newcommand{\conc}{\Hat{\ }}

\newcommand{\sq}[2]{\sideset{^{#1}}{}{\operatorname{#2}}}

\newcommand{\Mod}{\operatorname{Mod}}

\renewcommand{\phi}{\varphi}

\newcommand{\Union}{\bigcup}
\newcommand{\initial}\lessdot
 
\newcommand{\K}{\operatorname{\mathcal{K}}}

\newcommand{\C}{\mathfrak C}

\def\?{?\vadjust
{\vbox to 0pt{\vskip-7pt\hbox to 1.1\hsize{\hfill\huge ?!}}}}

\begin{document}
 
\title[Abstract Elementary Classes with No Maximal Models]{Categoricity in
Abstract Elementary Classes with No Maximal Models}

\author
{Monica VanDieren}
\email[Monica VanDieren]{mvd@umich.edu}
\address{Department of Mathematics\\
University of Michigan\\
Ann Arbor MI 48109}

\date{September 7, 2005}

\begin{abstract}
The results in this paper are in a context of abstract
elementary classes identified by Shelah and Villaveces in which the
amalgamation property is not assumed.  
 The long-term goal is to
solve Shelah's Categoricity Conjecture in this context.  Here we tackle a
problem of Shelah and Villaveces by proving that in their context, the
uniqueness of limit models follows from categoricity under the assumption
that the subclass of amalgamation bases is closed under unions of bounded,
$\prec_{\K}$-increasing chains.

\end{abstract}
\maketitle
\bigskip

\section*{Introduction} \label{s:introduction}
The origins of much of pure model theory can be traced back to
\L o\'{s}' Conjecture \cite{Lo}.
This conjecture was resolved by 
M. Morley in his Ph.D. thesis in 1962 \cite{Mo}.
Morley then questioned the status of the conjecture for
uncountable theories.
Building on work of W. Marsh, F. Rowbottom and J.P. Ressayre,
S. Shelah proved the statement for uncountable
theories in 1970 \cite{Sh 31}.
Out
of Morley and Shelah's proofs the program of
\emph{stability theory} or
\emph{classification theory} evolved.

While first-order logic has far reaching applications to other fields
of mathematics, there are several interesting frameworks which cannot be
captured by first-order logic.  
A classification theory for non-elementary classes
will open the door potentially to
a multitude of applications of model theory to classical mathematics and
provide insight into first-order model theory.

 Shelah posed a
generalization of 
\L o\'{s}' Conjecture to $L_{\omega_1,\omega}$ as a test question
to measure progress in non-first-order model theory.  
%
%
%
%
Focus on non-elementary classes began to shift in the late seventies
when Shelah, influenced by B. J\'{o}nsson's work in universal algebra (see
\cite{Jo1},
\cite{Jo2}), identified the notion of 
\emph{abstract elementary class (AEC)}  to capture many
non-first-order  logics \cite{Sh 88} including $L_{\omega_1,\omega}$ and
$L_{\omega_1,\omega}(\mathbf{Q})$.   
An abstract elementary class is a  class of
structures of the same similarity type endowed with a morphism satisfying
natural properties such as closure under directed limits.  

\begin{definition}\label{aec defn}\index{abstract elementary class}\index{AEC}

$\K$ is an \emph{abstract elementary class (AEC)} iff
$\K$ is a class of models for some vocabulary which is denoted by $L(\K)$,
and the class is equipped with a partial order, 
$\preceq_{\K}$ satisfying the following:
\begin{enumerate}
\item\label{closure under iso axiom} Closure under isomorphisms.
\begin{enumerate}
 
\item 
For every $M\in \mathcal{K}$
 and every $L(\K)$-structure $N$ if $M\cong N$ then $N\in \mathcal{K}$.
\item Let $N_1,N_2\in\K$ and  $M_1,M_2\in \K$ such that
there exist $f_l:N_l\cong M_l$ (for $l=1,2$) satisfying
$f_1\subseteq f_2$ then  $N_1\prec_{\K}N_2$ 
implies that  $M_1\prec_{\K}M_2$. 

\end{enumerate}
 
\item $\preceq_{\K}$ refines the submodel relation.

\item \label{funny axiom}
If $M_0, M_1\preceq_{\K} N$ and $M_0$ is a submodel of $M_1$,
then $M_0\preceq_{\K} M_1$.
\item\label{DLS} (Downward L\"{o}wenheim-Skolem Axiom) There is a
\emph{L\"{o}wenheim-Skolem number of $\K$}, denoted $LS(\K)$ which is the minimal
$\kappa$ such that for every $N\in\K$ and every 
$A\subset N$, there exists $M$ with $A\subseteq M\prec_{\K}N$ of cardinality
$\kappa+|A|$.
\item\label{union axiom} If $\langle M_i\mid i<\delta\rangle$ is a
$\prec_{\K}$-increasing and chain of models in $\K$ 
\begin{enumerate}
\item $\Union_{i<\delta}M_i\in\K$, 
\item for every $j<\delta$, $M_j\prec_{\K}\Union_{i<\delta}M_i$  and 
\item if $M_i\prec_{\K}N$ for every $i<\delta$, then 
$\Union_{i<\delta}M_i\prec_{\K}N$.
\end{enumerate}
\end{enumerate}
\end{definition}

\begin{definition}

For $M,N\in \K$ a monomorphism $f:M\rightarrow N$ is
called a \emph{$\prec_{\K}$-embedding} or a \emph{$\prec_{\K}$-mapping}
iff 
$f[M]\preceq_{\K} N$.

\end{definition}

\begin{notation}
We write $\K_{\mu}:=\{M\in\K\mid \|M\|=\mu\}$.  
\end{notation}

\begin{remark}
The Hanf number of $\K$ will be formally defined in Definition \ref{hanf
defn}.  It is bounded by
$\beth_{(2^{2^{LS(\K)}})^+}$.
\end{remark}

Shelah extended his categoricity conjecture for
$L_{\omega_1,\omega}$-theories in the following form in \cite{Sh 702}, see
also \cite{Shc}:

\begin{conjecture}[Shelah's Categoricity Conjecture]\label{cat
conj}\index{Shelah's Categoricity Conjecture!for abstract elementary classes} 
Let
$\K$ be an abstract elementary class.  If
$\K$ is categorical in some $\lambda\geq \Hanf(\K)$, then for every 
$\mu\geq\Hanf(\K)$, $\K$ is categorical in $\mu$.
\end{conjecture}

\begin{definition}

We say \emph{$\K$ is categorical in $\lambda$}\index{categorical} whenever there
exists exactly one model in $\K$ of cardinality $\lambda$ up to isomorphism.

\end{definition}

Despite the existence of 
over 1000 published pages of partial results towards this conjecture,
it remains open.
Since the mid-eighties, model theorists have
approached Shelah's conjecture from two different directions (see
\cite{Gr1} for a short history).  Shelah,  M. Makkai
and O. Kolman attacked the conjecture with set theoretic assumptions
\cite{MaSh}, \cite{KoSh}, \cite{Sh 472}.  On the
other hand, Shelah also looked at the conjecture under 
model theoretic assumptions in \cite{Sh 394}, \cite{Sh 576} and \cite{Sh
600}.  The approach of Shelah and A. Villaveces in \cite{ShVi 635}
involved a balance between set theoretic and model theoretic
assumptions.  This paper further investigates the
context of \cite{ShVi 635} which we delineate here:
\begin{assumption}\label{assm intro}\index{context of Shelah and Villaveces}
\begin{enumerate}
\item $\K$ is an AEC with no maximal models
with respect to the relation $\prec_{\K}$,
\item $\K$ is categorical in some fixed $\lambda\geq\Hanf(\K)$,
\item GCH holds and
\item a form of the weak diamond holds, namely
$\Phi_{\mu^+}(S^{\mu^+}_{\cf(\mu)})$ holds for every
$\mu$ with $\mu<\lambda$ (see Definition \ref{weak diamond defn}).
\end{enumerate}
\end{assumption}

The purpose of \cite{ShVi 635} was to begin investigating the conjecture
that the amalgamation property follows from categoricity in a large
enough cardinality.  
All of the other attempts to prove Conjecture \ref{cat conj} have made use
of the assumption of the amalgamation property which is a sufficient
condition to define a reasonable notion of (Galois)-type (see Section
\ref{s:amalg}).

\begin{definition}
Let $\K$ be an abstract elementary class and $\mu$ a cardinal
$\geq\LS(\K)$.
\begin{enumerate}
\item 
We say that $M\in\K_\mu$ is an \emph{amalgamation base} 
if 
for every $N_1,N_2\in\K_{\mu}$
and $g_i:M\rightarrow N_i$ for $(i=1,2)$, there are
$\prec_{\K}$-embeddings $f_i$, $(i=1,2)$ and a model $N$ such that
the following diagram commutes:
 
\[
\xymatrix{\ar @{} [dr] N_1
\ar[r]^{f_1}  &N \\
M \ar[u]^{g_1} \ar[r]_{g_2} 
& N_2 \ar[u]_{f_2} 
}
\]

\item
An abstract elementary class $\K$ satisfies the \emph{amalgamation
property} iff 
every $M\in\K$ 
is an amalgamation base.
\item
We write $\K^{am}$\index{$\K^{am}$}\index{amalgamation base!$\K^{am}$} for
the class of amalgamation bases which are in
$\K$.  We also use $\K^{am}_\mu$ to denote the class of amalgamation
bases of cardinality $\mu$.

\end{enumerate}

\end{definition}

\begin{remark}\begin{enumerate}
\item 
The definition of amalgamation base varies across the literature.
Our definition of amalgamation base is weaker than an alternative formulation
which does not put any restriction on the cardinality of $N_1$ and $N_2$. 
Under the assumption
of the amalgamation property, these definitions are known to be
equivalent.  However, in this context where the amalgamation property is
not assumed, we cannot guarantee the existence of the stronger form of
amalgamation bases.
\item
We get an equivalent definition of amalgamation base, if we additionally
require that $g_i\restriction M=id_M$ for $i=1,2$, in the definition above.
See \cite{Gr2} for details.
\end{enumerate}
\end{remark}

%
%
%
%

It is conjectured that categoricity in a large enough cardinality implies
the amalgamation property.  However, there are examples of abstract
elementary classes which are categorical in $\omega$ successive
cardinals, but fail to have the amalgamation property in larger
cardinalities \cite{HaSh}, \cite{ShVi2}.  Shelah 
constructs an abstract elementary
class  whose models are
bipartite random graphs.  Models of cardinality $\aleph_1$ in this class
witness the failure of amalgamation.  Intriguingly, under the assumption
of Martin's Axiom, this class of bipartite graphs is categorical in
$\aleph_0$ and $\aleph_1$.  On the other hand, if one assumes a
version of the weak diamond, Shelah proves that categoricity in
$\aleph_0$ and $\aleph_1$ implies amalgamation in $\aleph_1$ (\cite{Sh
88} or see
\cite{Gr1} for an exposition).  There are other natural examples of
abstract elementary classes which do not satisfy the amalgamation
property but are unstable such as the class of locally finite groups
\cite{Ha}.

Limited progress has been made to prove that amalgamation follows from
categoricity.  Kolman and Shelah manage to prove this for  AECs that
can be axiomatized by a $L_{\kappa,\omega}$ sentence with
$\kappa$ a measurable cardinal \cite{KoSh}.  They first introduce limit models
as a substitute for saturated models, and then prove the uniqueness of
limit models (see Definition \ref{defn limit}).

%

 To better understand the relationship between the
amalgamation property, categoricity and the uniqueness of limit models,
consider the questions of uniqueness and existence of limit models in
classes which satisfy the amalgamation property, but not are not
necessarily categorical:
\begin{remark}
Even under the amalgamation property, the uniqueness and existence of
limit models do not come for free.  The existence requires 
stability (see \cite{Sh 600} or \cite{GrVa}).
The question of uniqueness of limit models is tied into (super)stability
as well.  Even in first-order logic, the uniqueness of limit models fails
for un-superstable theories (see \cite{GrVa} or \cite{Sh 394} for
examples).  The uniqueness of limit models has  been proven
in AECs under the assumption of categoricity (\cite{KoSh}, \cite{Sh 394},
and here, Theorem
\ref{uniqueness of limit
models}).  Recently Grossberg, VanDieren and Villaveces identified
sufficient conditions (which are consequences of superstability) for the
uniqueness of limit models in classes with the amalgamation property
\cite{GrVaVi}.

\end{remark}

The motivation for this paper is to elaborate on recent work of Shelah
and Villaveces in which they strive to prove under weaker assumptions
than Kolman and Shelah that the amalgamation property follows from
categoricity above the Hanf number.  The first step in proving
amalgamation is to show the uniqueness of limit models.

The uniqueness of limit models under Assumption \ref{assm intro}
generalizes Theorem 6.5 of \cite{Sh 394} where Shelah assumes the full
amalgamation property.  The amalgamation property is used in
\cite{Sh 394} in several forms including the fact that saturated models
are model homogeneous and that all reducts of Ehrenfeucht-Mostowski
models are amalgamation bases.  Shelah then uses the uniqueness of limit
models to prove that the union of a chain of $\mu$-saturated models is
$\mu$-saturated, provided that the chain is of length $<\mu^+$.  This is
one of the main steps in proving a downward categoricity transfer theorem
for classes with the amalgamation property.

In the Fall of 1999, we identified several problems with Shelah and
Villaveces' proof of the uniqueness of limit models from \cite{ShVi 635}. 
After two years of correspondence, Shelah and Villaveces conceded that
they were not able to resolve these problems.   
While these issues are undertaken in this paper, to date the proof
of the uniqueness of limit models has resisted a complete solution under
Assumption \ref{assm intro}.  After presenting a partial solution (Theorem
\ref{uniqueness of limit models}) of the uniqueness of limit
models and discussing this with  Shelah at a Mid-Atlantic Mathematical
Logic Seminar in the Fall of 2001,  we were not able to remove
the extra hypothesis.   The extra hypothesis was weakened in \cite{Va}. 
This paper provides a complete proof of an intermediate uniqueness result
patching a gap that was found in
\cite{Va} in the Fall of 2002.  The partial solution  to the uniqueness of
limit models described here is in the context identified in
\cite{ShVi 635} (Assumption \ref{assm intro}) under the hypothesis:

\begin{description}
\item[Hypothesis 1]  Every continuous tower inside $\C$ has an amalgamable
extension inside $\C$ (see Sections
\ref{s:limit models} and \ref{s:tower} for the definitions).
\end{description}

\begin{remark}
The model $\C$ in Hypothesis 1 is not the usual monster model.  It is a
weak substitute for the monster model and is introduced in Section
\ref{s:limit models}.  Monster models, as we know them in first-order
logic, are model homogeneous.  In the absence of the amalgamation
property, model homogeneous models may not exist.
\end{remark}

In the context of \cite{ShVi 635}, Hypothesis 1 is a consequence of the
more natural Hypothesis 2 (see Section \ref{s:<b extension property}).

\begin{description}
\item[Hypothesis 2]
For $\mu<\lambda$, the class of amalgamation bases of cardinality $\mu$
(denoted by $\K^{am}_\mu$) is closed under unions of
$\prec_{\K}$-increasing chains of length $<\mu^+$.
\end{description}


It seems reasonable to consider a weakening of Grossberg's Intermediate
Categoricity Conjecture which captures Hypothesis 2:
\begin{conjecture}
Let $\K$ be an AEC.  If there exists a $\lambda\geq\Hanf(\K)$ such that
$\K$ is categorical in $\lambda$, then $\K^{am}_\mu$ is closed under
unions of length $<\mu^+$ for all $\mu$ with $LS(\K)\leq\mu<\lambda$.
\end{conjecture}

 Although  Theorem
1.11 of Chapter 4 in
\cite{Sh 300} addresses a similar problem to Hypothesis 2, this
statement may be too ambitious to prove.  An alternative hypothesis which
also implies Hypothesis 1 is

\begin{description}
\item[Hypothesis 3]
The union of a $\prec_{\K}$-increasing chain of length $<\mu^+$ of limit
models of cardinality $\mu$ is a limit model.
\end{description}

Hypothesis $3$ may be more approachable as it is a relative of the first-order consequence of superstability that the union of a
$\prec$-increasing chain of
$\kappa(T)$-many saturated models is saturated.

Hypothesis 1 has relatives in the literature as well.  Indeed in \cite{Sh
88} where the amalgamation property is not assumed, Shelah identifies the
link between the existence of maximal elements of $\K^3_{\aleph_0}$
(a specialization of towers of length 1) and $2^{\aleph_1}$ non-isomorphic
models in
$\aleph_1$.

We thank R. Grossberg for his kind and generous guidance on
Chapter II of the Ph.D.
thesis \cite{Va} on which this paper is based. We are indebted to J.
Baldwin for his unending dedication to improving the clarity and
correctness of this paper and for chairing the thesis examination
committee for \cite{Va}.   A. Villaveces' comments and discussions on this
paper and on \cite{Va} were invaluable.   Additionally, we express our
gratitude to the referee whose patience and comments contributed
to a much needed revision of the original work.

This paper is divided into 3 parts outlined below.

\noindent {\bf Part I.}  The first part summarizes the
necessary definitions and background material.  It also includes some new
results on $\mu$-splitting.
\begin{enumerate}
\item[\S \ref{s:amalg}] Galois-types
\item[\S \ref{s:limit models}] Limit Models

\item[\S \ref{s:ab}] Limit Models are Amalgamation Bases

\item[\S \ref{s:mu-split}] $\mu$-Splitting

\item[\S \ref{s:tower}] Towers

\end{enumerate}
\noindent{\bf Part II.}  Here we provide a complete proof
of the uniqueness of limit models under Hypothesis 1 and Assumption
\ref{assm intro}.
\begin{enumerate}
\item[\S \ref{s:full}] Relatively Full Towers
%

\item[\S \ref{s:reduced and full}]  Continuous
$<^c_{\mu,\alpha}$-Extensions

\item[\S \ref{s:dense towers}]  Refined Orderings on Towers

\item[\S \ref{s:unique limits}]  Uniqueness of Limit Models

\end{enumerate}
\noindent{\bf Part III.}
In this part of the paper we include a partial result in the direction of
Hypothesis 1 and discuss reduced towers.
\begin{enumerate}
\item[\S \ref{s:<b extension property}]$<^c_{\mu,\alpha}$-Extension
Property for Nice Towers

\item[\S \ref{s:reduced new section}] Reduced Towers  

\end{enumerate}

\bigskip

\part{Preliminaries}\label{prelim part}

Throughout this paper, unless otherwise stated, we will make Assumption
\ref{assm intro} and $\mu$ will be a cardinal satisfying
$LS(\K)\leq\mu<\lambda$ where $\lambda$ is the categoricity cardinal.

We introduce the necessary
definitions and background from \cite{ShVi 635}.  The reader familiar
with \cite{ShVi 635} may skim through Section \ref{s:limit models} where
the monster model is introduced and then proceed to Section
\ref{s:mu-split} which includes some new results on $\mu$-splitting.

\renewcommand{\thetheorem}{\Roman{part}.\arabic{section}.\arabic{theorem}}

\bigskip
\section{Galois-types}\label{s:amalg}
 
In this section we discuss problems that arise when working without the
amalgamation property in AECs.  The first obstacle is to identify a
reasonable notion of type.  Because of the category-theoretic definition
of abstract elementary classes, the first-order notion of 
formulas and types cannot
be applied.  To overcome this barrier, Shelah has suggested identifying
types, not with formulas, but with the orbit of an element under
the group of automorphisms fixing a given structure.  In order to
carry out this definition of type, the following 
binary relation $E$
must be an equivalence relation on triples $(a,M,N)$.  
In order to avoid confusing this new notion of ``type'' with the
conventional
one (i.e. set of formulas) we will follow 
\cite{Gr1} and \cite{Gr2} and
 introduce it below under
the name of \emph{Galois-type}.  

\begin {definition}\index{$E$, binary relation}
For triples $(\bar a_l, M_l, N_l)$ where $\bar a_l\in N_l$
and
$M_l\preceq_{\K}N_l\in\K$ for
$l=1,2$, we define a binary relation $E$ as follows:
$(\bar a_1, M_1, N_1)E(\bar a_2, M_2, N_2)$ iff 
$M:=M_1=M_2$ and there exists $N\in\K$ and $\prec_{\K}$-mappings
$f_1,f_2$ such that
$f_l:N_l\rightarrow N$ and $f_l\restriction M=id_M$ for $l=1,2$ and
$f_1(\bar a_1)=f_2(\bar a_2)$:
 
\[
\xymatrix{\ar @{} [dr] N_1
\ar[r]^{f_1}  &N \\
M \ar[u]^{id} \ar[r]_{id} 
& N_2 \ar[u]_{f_2} 
}
\]

\end{definition}

To prove that $E$ is an equivalence relation (more specifically, that $E$
is transitive), we need to restrict ourselves to 
amalgamation bases.  


\begin{remark}
$E$ is an equivalence relation on the set of triples of the form
$(\bar a, M, N)$ where $M\preceq_{\K}N$, $\bar a\in N$ and 
$M,N\in\K^{am}_\mu$ for fixed $\mu\geq LS(\K)$.  
To see that $E$ is transitive, consider $(a_1,M,N_1)E(a_2,M,N_2)$ and
$(a_2,M,N_2)E(a_3,M,N_3)$ where $M,N_1,N_2,N_3\in\K^{am}_\mu$.  Let
$N_{1,2}$ and
$f_1$,
$f_2$ be such that 
$f_1:N_1\rightarrow N_{1,2}$; $f_2:N_2\rightarrow N_{1,2}$ and
$f_1\restriction M=f_2\restriction M=id_M$ with
$f_1(a_1)=f_2(a_2)$. Similarly define $g_2$, $g_3$ and $N_{2,3}$ with
$g_2(a_2)=g_3(a_3)$.  By the Downward L\"{o}wenheim-Skolem Axiom, we may
assume that
$N_{1,2}$ and $N_{2,3}$ have cardinality $\mu$. Consider the
following diagram of this situation.

\[
\xymatrix{\ar @{} [dr] N_1
\ar[r]^{f_1}  &N_{1,2} \\
M \ar[u]^{id} \ar[r]_{id}\ar[d]_{id} 
& N_2 \ar[u]_{f_2} \ar[d]^{g_2}\\
N_3\ar[r]_{g_3}&N_{2,3}
}
\]

Since $N_2$ was chosen to be an amalgamation base, we can amalgamate
$N_{1,2}$ and $N_{2,3}$ over $N_2$ with mappings $h_1$ and $h_3$ 
and an amalgam $N^*$ giving us the following diagram:

\[
\xymatrix{\ar @{} [dr] N_1
\ar[r]^{f_1}  &N_{1,2}\ar[dr]^{h_1} &\\
M \ar[u]^{id} \ar[r]_{id} \ar[d]_{id}
& N_2 \ar[u]_{f_2} \ar[d]^{g_2}&N^*\\
N_3\ar[r]_{g_3}&N_{2,3}\ar[ru]_{h_3}&
}
\]
Notice that $h_1(f_1(a_1))=h_3(g_3(a_3))$.  Thus $h_1\circ f_1$ and
$h_3\circ g_3$ witness that $(a_1,M,N_1)E(a_3,M,N_3)$.

\end{remark}


\begin{remark}[Invariance]\label{f preserves ab}
If $M$ is an
amalgamation base and
$f$ is an
$\prec_{\K}$-embedding, then 
$f(M)$ is an amalgamation base.
\end{remark}

In AECs with the amalgamation property, we are often limited to speak of
types only over models.  Here we are further restricted to deal with
types only over models which are amalgamation bases.
 
\begin{definition}Let $\mu\geq LS(\K)$ be given.
\begin{enumerate}
\item\index{Galois-type}\index{type, Galois}\index{$\tp(a/M,N)$} For
$M,N\in\K^{am}_\mu$ with $M\preceq_{\K}N$ and
$\bar a\in
\sq{\omega>}|N|$, the \emph{Galois-type of $\bar a$ in $N$ over $M$},
written
$\tp(\bar a/M,N)$, is defined to be $(\bar a,M,N)/E$.

\item\index{$\gaS(M)$}\index{Galois-type!$\gaS(M)$} For $M\in\K^{am}_\mu$,
$$\gaS^1(M):=\{\tp(a/M,N)\mid M\preceq_{\K} N\in\K^{am}_\mu,  a\in N\}.$$

\item \index{Galois-type!realized}We say \emph{$p\in \gaS(M)$ is realized
in
$M'$} whenever
$M\prec_{\K}M'$ and there exist $\bar a\in M'$ and $N\in\K^{am}_{\mu}$
such that
$p=(\bar a,M,N)/E$.

\item \index{Galois-type!restriction}For $M'\in\K^{am}_\mu$ with
$M\prec_{\K}M'$ and
$q=\tp(\bar a/M',N)\in
\gaS(M')$, we define \emph{the restriction of $q$ to $M$} as
$q\restriction M:=\tp(\bar a/M,N)$.
\item \index{Galois-type!extension}For $M'\in\K^{am}_\mu$ with
$M\prec_{\K}M'$, we say that
\emph{$q\in
\gaS(M')$ extends $p\in \gaS(M)$} iff $q\restriction M=p$.
\item $p\in\gaS(M)$ is said to be \emph{non-algebraic} if no $a\in M$
realizes
$p$.

\end{enumerate}
\end{definition}
 
%

\begin{notation}\index{$\tp(a/M)$}
We will often abbreviate a Galois-type, $\tp(a/M,N)$ as $\tp(a/M)$, when
the role of $N$ is not crucial or is clear.  This occurs mostly when we
are working inside of a fixed structure $\C$, which we define in
Section \ref{s:limit models}.
\end{notation}

\begin{fact}[see \cite{Gr2}]
When $\K=\Mod(T)$ for $T$ a complete first-order theory, the above
definition of $\tp(a/M,N)$ coincides with the classical first-order
definition where $c$ and $a$ have the same type over $M$ iff for every
first-order formula $\varphi(x,\bar b)$ with parameters $\bar b$ from
$M$, 
$$N\models\varphi(c,\bar b)\text{ iff }N\models\varphi(a,\bar b).$$
\end{fact}

We will now define Galois-stability in an analogous way:

\begin{definition}\index{stable}
We say that $\K$ is 
\emph{Galois-stable in $\mu$} if for 
every $M\in\K^{am}_\mu$,
$|\gaS^1(M)|=\mu$.
\end{definition}

\begin{fact}[Fact 2.1.3 of \cite{ShVi 635}]\label{cat implies stab}
If $\K$ is categorical in $\lambda$, then for every $\mu<\lambda$,
we have that $\K$ is Galois-stable in $\mu$.
\end{fact}

By combining results from \cite{ShVi 635}, \cite{GrVa} and \cite{BaKuVa}
it is possible to improve this to conclude Galois-stability in some 
cardinals $\geq\lambda$, but it remains open whether or not in AECs
categoricity implies Galois-stability in all cardinalities above $LS(\K)$.

\begin{definition}
Let $\mu>LS(\K)$, $M$ is said to be \emph{$\mu$-saturated} if for every
$N\prec_{\K}M$ with
$N\in\K^{am}_{<\mu}$ and every Galois-type $p$ over $N$, we have that $p$
is realized in
$M$.
\end{definition}

The following fact is proved by showing the equivalence of model
homogeneous models and saturated models in classes which satisfy the
amalgamation property \cite{Sh 576}.

\begin{fact}

Suppose that $\K$ satisfies the amalgamation property.  
If
$M_1$ and $M_2\in\K_\mu$ are
$\mu$-saturated and there exists $N\prec_{\K}M_1,M_2$ with
$N\in\K_{<\mu}$, then
$M_1\cong M_2$.

\end{fact}

Since we will be working in a context where the amalgamation property is
not assumed, we do not have the uniqueness of saturated models at hand. 
In fact even the existence of saturated models is questionable.  The
purpose of this paper is to identify a suitable substitute for saturation
that is unique up to isomorphism in every cardinality.  The candidate is
limit model discussed in the following section.  
Later we will give an alternative characterization of limit models as the
union of a relatively full tower (see Section \ref{s:full}).  This
characterization plays the role of $\mathbf F^a_{\kappa}$-saturated
models from first-order model theory (see Chapter IV of \cite{Shc}).

\section{Limit Models}\label{s:limit models}

In this section we define limit models and discuss their uniqueness and
existence.
A local substitute for the monster model is also introduced.

We begin with universal extensions which are central in the
definition of limit models.  A universal extension captures some
properties of saturated models without referring explicitly to types.
The notion of universality over countable models was first analyzed by
Shelah in Theorem 1.4(3) of \cite{Sh 87a}.

\begin{definition}
\begin{enumerate}
\item \index{universal over!$\kappa$-universal over}
Let $\kappa$ be a cardinal $\geq LS(\K)$.
We say that $N$ is \emph{$\kappa$-universal over $M$} iff 
for every $M'\in\K_{\kappa}$ with $M\prec_{\K}M'$ there exists
a $\prec_{\K}$-embedding $g:M'\rightarrow N$ such that
$g\restriction M=id_{M}$:

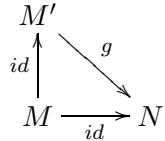
\[
\xymatrix{\ar @{} [dr] M'
\ar[dr]^{g}  & \\
M \ar[u]^{id} \ar[r]_{id} 
& N  
}
\]

\item \index{universal over}

We say $N$ is \emph{universal over $M$} or $N$ is \emph{a universal
extension of $M$} iff 
$N$ is $\|M\|$-universal over $M$.

\end{enumerate}
\end{definition}

\begin{notation}\index{$
M\stackrel{id}{\Longrightarrow}N$}
In diagrams, we will indicate that $N$ is universal over $M$, by
writing $
\def\labelstyle{\scriptstyle}
{\xymatrix 
{M\ar@2{->}[r]^{id}&N}}$.
\end{notation}

\begin{remark}
Notice that 
the definition of \emph{$N$ universal over $M$} requires all extensions
of $M$ of cardinality $\|M\|$ to be embeddable into $N$.  First-order
variants of this definition in the literature often involve
$\|M\|<\|N\|$.  We will  be considering the case when $\|M\|=\|N\|$.
\end{remark}

\begin{remark}
Suppose that $T$ is a first-order complete theory that is stable in some
regular $\mu$. Then every model $M$ of $T$ of cardinality $\mu$ has
an elementary extension $N$ of cardinality $\mu$ which is universal over
$M$. To see this, define an elementary-increasing and continuous chain of
models of $T$ of cardinality $\mu$,
$\langle N_i\mid i<\mu\rangle$ such that $N_{i+1}$ realizes all types
over $N_i$.  Let $N=\Union_{i<\mu}N_i$.  By a back-and-forth
construction, one can show that $N$ is universal over $M$.
\end{remark}

The existence of universal extensions in AECs follows from categoricity
in $\lambda$ and GCH or categoricity and uses the presentation of the
model of cardinality
$\lambda$ as a reduct of an EM-model.

\begin{fact}[Theorem 1.3.1 from \cite{ShVi 635}]\label{exist
univ}\index{universal over!existence of universal extensions} Let $\mu$
be such that $LS(\K)\leq\mu<\lambda$.  Then every element of $\K^{am}_\mu$
has a universal extension in
$\K^{am}_\mu$.
\end{fact}

Another existence result that does not use GCH or categoricity can be
proved under the assumption of Galois-stability  and the
amalgamation property (\cite{Sh 600} or see
\cite{GrVa} for a proof).

Notice that the following observation asserts that it is unreasonable 
to prove a stronger existence statement than Fact \ref{exist univ},
without having proved the amalgamation property.

\begin{proposition}\label{univ ext implies ab}\index{amalgamation base}
If $M\in\K_\mu$ has a universal extension, then $M$ is an amalgamation
base.
\end{proposition}

We can now define the principal concept of this paper:

\begin{definition}\label{defn limit}
For $M',M\in\K_\mu$ and $\sigma$ a limit ordinal with
$\sigma<\mu^+$,
we say that $M'$ is a \emph{$(\mu,\sigma)$-limit over $M$}\index{limit
model!$(\mu,\sigma)$-limit model over $M$} iff there exists a
$\prec_{\K}$-increasing and continuous sequence of models $\langle M_i\in
\K_{\mu}\mid i<\sigma
\rangle$
such that 
\begin{enumerate}
\item $M= M_0$,
\item $M'=\Union_{i<\sigma}M_i$ 
\item\label{ab cond in defn} for $i<\sigma$, $M_i$ is an amalgamation base
and
\item\label{univ cond in defn} $M_{i+1}$ is universal over $M_i$.
\end{enumerate}
\end{definition}

\begin{remark}
\begin{enumerate}
\item Notice that in Definition \ref{defn limit}, for $i<\sigma$ and
$i$ a limit ordinal, $M_i$ is a $(\mu,i)$-limit model.
\item Notice that Condition (\ref{univ cond in defn}) implies Condition
(\ref{ab cond in defn}) of Definition 
\ref{defn limit}.  In our constructions, since the question of whether a
particular model is an amalgamation base becomes crucial, we choose to
list this as a separate condition.
\end{enumerate}
\end{remark}

\begin{definition}\label{limit defn}
We say that $M'$ is a \emph{$(\mu,\sigma)$-limit}\index{limit
model!$(\mu,\sigma)$-limit} iff there is some
$M\in\K$ such that $M'$ is a $(\mu,\sigma)$-limit over $M$.
\end{definition}


 While limit models were used is \cite{KoSh} and \cite{Sh 394}, their use
extends to other contexts.  There is evidence that the
uniqueness of limit models  provides a basis for the development of a
notion of non-forking and a stability theory for abstract elementary
classes. Limit
models are used in \cite{GrVa} to develop the notion of non-splitting in
tame, Galois-stable AECs. The uniqueness of limit models implies the
existence of superlimits in
\cite{Sh 576}. 
Additionally, in \cite{Sh 600} the uniqueness of limit models appears as
an axiom for good frames and the limit models are closely related to
brimmed-models.
 In all of these applications, limit models provide a substitute for
Galois-saturated models.

By repeated applications of Fact \ref{exist univ}, the 
existence of $(\mu,\omega)$-limit models can be proved:

\begin{fact}[Theorem 1.3.1 from \cite{ShVi 635}]\label{exist
mu,theta}\index{limit model!existence} Let $\mu$ be a cardinal such
that
$\mu<\lambda$.
For every $M\in\K_\mu^{am}$, there is a $(\mu,\omega)$-limit over $M$.
\end{fact}

In order to extend this argument further to yield the existence
of $(\mu,\sigma)$-limits for arbitrary limit ordinals $\sigma<\mu^+$,
we need to be able to verify that limit models are in fact amalgamation bases.
We will examine this in Section \ref{s:ab}.

While the existence of limit models can be derived
from the categoricity and weak diamond assumptions, 
the uniqueness of limit models is more difficult.  Here we recall two
easy uniqueness facts which state that limit models of the same length
are isomorphic.  They are proved using the natural back-and-forth
construction of an isomorphism.

\begin{fact}[Fact 1.3.6 from \cite{ShVi 635}]
\label{unique limits}\index{uniqueness of limit models!of the same length}
Let $\mu\geq LS(\K)$ and $\sigma<\mu^+$.  
If $M_1$ and $M_2$ are $(\mu,\sigma)$-limits over $M$, then
there exists an isomorphism $g:M_1\rightarrow M_2$ such that
$g\restriction M = id_M$.  Moreover if $M_1$ is a $(\mu,\sigma)$-limit
over $M_0$; $N_1$ is a $(\mu,\sigma)$-limit over $N_0$ and
$g:M_0\cong N_0$, then there exists a $\prec_{\K}$-mapping, $\hat g$,
extending $g$ such that $\hat g:M_1\cong N_1$.
\end{fact}

\[
\xymatrix{\ar @{} [dr] M_1
\ar@.[r]^{\hat g}  &N_1 \\
M_0 \ar@2{->}[u]^{id} \ar[r]_{g} 
& N_0 \ar@2{->}[u]^{id}
}
\]

\begin{fact}[Fact 1.3.7 from \cite{ShVi 635}]
\label{sigma and cf(sigma) limits}\index{uniqueness of limit models!of the
same cofinality} Let $\mu$ be a cardinal and $\sigma$ a limit ordinal with
$\sigma<\mu^+\leq\lambda$.  If $M$ is a $(\mu,\sigma)$-limit
model, then $M$ is a $(\mu,cf(\sigma))$-limit model.
\end{fact}

A more challenging uniqueness question is to prove that two limit models
of different lengths ($\sigma_1\neq\sigma_2$) are isomorphic:
\begin{conjecture}
Suppose that $\K$ is categorical in some $\lambda\geq\Hanf(\K)$ and $\mu$
is a cardinal with $LS(\K)\leq\mu<\lambda$. Let $\sigma_1$ and $\sigma_2$
be limit ordinals
$<\mu^+$. Suppose $M_1$ and $M_2$ are $(\mu,\sigma_1)$- and
$(\mu,\sigma_2)$-limits over $M$, respectively.  Then $M_1$ is isomorphic
to $M_2$ over $M$.
\end{conjecture}

The main result of this paper, Theorem \ref{uniqueness of limit
models},  is a solution to this conjecture under Assumption \ref{assm
intro} and Hypothesis 1.

We will need one more notion of limit model, which will later serve as a
substitute for a monster model.  This is a natural extension of the
limit models already defined:

\begin{definition}\index{limit
model!$(\mu,\mu^+)$}\index{$(\mu,\mu^+)$-limit model} Let
$\mu$ be a cardinal
$<\lambda$, we say that $\check M$ is a 
\emph{$(\mu,\mu^+)$-limit over $M$} iff there exists a
$\prec_{\K}$-increasing and continuous chain of models $\langle
M_i\in\K^{am}_\mu\mid i<\mu^+\rangle$ such that
 $M_0=M$,
 $\Union_{i<\mu^+}M_i=\check M$, and 
for $i<\mu^+$, $M_{i+1}$ is universal over $M_i$.

\end{definition}

\begin{remark}
While it is known that in our context $(\mu,\theta)$-limit models are
amalgamation bases when $\theta<\mu^+$, it is open whether or not
$(\mu,\mu^+)$-limits are amalgamation bases.  
To avoid confusion between these two concepts of limit models,
we will denote $(\mu,\mu^+)$-limit models with a $\check{}$
above the model's name (i.e. $\check M$).  Later we will avoid this
confusion by fixing a $(\mu,\mu^+)$-limit model and denoting it by $\C$,
since it will substitute the usual notion of a monster model.
\end{remark}

The existence of $(\mu,\mu^+)$-limit models follows from the fact
that $(\mu,\theta)$-limit models are amalgamation bases when 
$\theta<\mu^+$, see Corollary \ref{exist limit}.\index{limit
model!existence}  The
uniqueness of $(\mu,\mu^+)$-limit models (Proposition \ref{mu,mu^+-limits
are wmh}) can be shown using an easy back-and-forth
construction as in the proof of Fact \ref{unique limits}.

  The following theorem indicates that $(\mu,\mu^+)$-limits provide
some level of homogeneity.  First we recall an exercise regarding
amalgamation.

\begin{remark}\label{renaming amalgam}
Suppose that $M_0$, $M_1$ and $M_2$ can be amalgamated, then by renaming
elements, we can choose the amalgam to be a $\prec_{\K}$-extension of
$M_2$.
\end{remark}

\begin{theorem}\label{mu,mu+ limit is univ}
If $\check M$ is a $(\mu,\mu^+)$-limit, then for every $N\prec_{\K}\check
M$ with $N\in\K^{am}_\mu$, we have that $\check M$ is universal over
$N$.  Moreover, $\check M$ is a $(\mu,\mu^+)$-limit
over $N$.
\end{theorem}
\begin{proof}
Suppose that $\check M$ is a $(\mu,\mu^+)$-limit model and
$N\prec_{\K}\check M$ is such that $N\in\K^{am}_\mu$.  Let $N'$ be an
extension of $N$ of cardinality $\mu$.  Let $\langle M_i\mid
i<\mu^+\rangle$ witness that $\check M$ is a $(\mu,\mu^+)$-limit model. 
Since $N$ has cardinality $\mu$, there exists $i<\mu^+$, such that
$N\prec_{\K}M_i$.  Since $N$ is an amalgamation base, we can amalgamate
$M_i$ and $N'$ over $N$ with amalgam $M'\in\K_\mu$.  By Remark
\ref{renaming amalgam}, we may assume that
$M_i\prec_{\K}M'$.

\[
\xymatrix{\ar @{} [dr] N'
\ar[r]^{h}  &M' \\
N \ar[u]^{id} \ar[r]_{id} 
& M_i \ar[u]_{id} 
}
\]

Since $M_{i+1}$ is universal over $M_i$, there is $g:M'\rightarrow
M_{i+1}$ such that $g\restriction M_i=id_{M_i}$.  Then $g\circ h$ give us
the desired mapping from $N'$ into $\check M$ over $N$.
\[
\xymatrix{\ar @{} [dr] N'
\ar[r]^{h}  &M'\ar[dr]^{g}& \\
N \ar[u]^{id} \ar[r]_{id} 
& M_i \ar[u]_{id} \ar@2{->}[r]_{id} &M_{i+1}
}
\]

\end{proof}

\begin{remark}
If $N$ is not an
amalgamation base, then there are no universal models over $N$.
\end{remark}

It is immediate that $\C$ realizes many types:
\begin{corollary}\label{C is saturated}

For every $M\in\K^{am}_\mu$ with $M\prec_{\K}\C$, we have that $\C$ is
saturated over $M$.  
\end{corollary}

\begin{corollary}\label{mu,mu^+-limits are
wmh}\index{$(\mu,\mu^+)$-limit model!homogeneity} Suppose $\check M_1$ and
$\check M_2$ are
$(\mu,\mu^+)$-limits over $M_1,M_2\in\K^{am}_\mu$, respectively.
If there exists an isomorphism $h:M_1\cong M_2$, then
$h$ can be extended to an isomorphism $g:\check M_1\cong \check M_2$.
\end{corollary}

Since $(\mu,\mu^+)$-limit models are unique and are universal over all
amalgamation bases of cardinality $\mu$, they are in some sense
homogeneous.  We will see that 
if $\check M$ is a $(\mu,\mu^+)$-limit
model and $\tp(a/M,\check M)=\tp(b/M,\check M)$, then there is an automorphism
$f$ of
$\check M$ fixing $M$ such that $f(a)=b$ (Corollary \ref{type aut}).
In some ways, $(\mu,\mu^+)$-limit models behave like monster models in
first-order logic if we restrict ourselves to amalgamation bases and
models of cardinality $\mu$.  This justifies the following notation.

\begin{notation}\label{C note} We  fix a cardinal $\mu$ with
$LS(\K)\leq\mu<\lambda$ and a
$(\mu,\mu^+)$-limit model and denote it by $\C$.  For $M\prec_{\K}\C$ we
abbreviate $$\{f\mid f\text{ is an automorphism of }\C\text{ with
}f\restriction M=id_M\}$$ by $\Aut_M(\C)$. 
\end{notation}

While it is customary to work entirely inside of a fixed monster model
$\C$ in first-order logic, we will sometime need to consider structures
outside of $\C$ since we do not have the full power of model homogeneity
in this context.

We now recall a result from \cite{ShVi 635} which will be used in our
proof of Corollary \ref{type aut}.  Although
Shelah and Villaveces work without the amalgamation property
as an assumption, using weak diamond they prove a weak amalgamation
property, which they refer to as \emph{density of amalgamation bases}.

\begin{fact}[Theorem 1.2.4 from \cite{ShVi 635}]\label{density of
ab}\index{amalgamation base!density of} Every $M\in\K_{<\lambda}$
has a proper $\K$-extension of the same cardinality which is an
amalgamation base. 

\end{fact}

We can now improve Fact \ref{exist univ} slightly.  
This improvement is used throughout 
this paper.

\begin{lemma}\label{univ with a first}
For every $\mu$ with 
$LS(\K)\leq\mu<\lambda$, if $M\in\K^{am}_{\mu}$, $N\in\K$ and 
$\bar a\in\sq{\mu^+>}|N|$ are such that $M\prec_{\K}N$, then
there exists $M^{\bar a}\in\K^{am}_\mu$ such that $M^{\bar a}$ is universal
over $M$ and $M\Union\bar a\subseteq M^{\bar a}$.
\end{lemma}
\begin{proof}
By Axiom \ref{DLS} of AEC, we can find $M'\prec_{\K}N$
of cardinality $\mu$ containing $M\Union\bar a$.  Applying 
Fact \ref{density of ab}, there exists an amalgamation base of cardinality
$\mu$, say $M''$, extending $M'$.  By Fact 
\ref{exist univ}
we can find a universal extension of $M''$ of cardinality $\mu$, say 
$M^{\bar a}$.

Notice that $M^{\bar a}$ is also universal over $M$.  Why?  Suppose 
$M^*$ is an
extension of $M$ of cardinality $\mu$.  Since $M$ is an amalgamation base we
can amalgamate $M''$ and $M^*$ over $M$.  WLOG we may assume that the
amalgam, $M^{**}$, is an extension of $M''$ of cardinality $\mu$ and a
$\prec_{\K}$-mapping
$f^{*}:M^{*}\rightarrow M^{**}$ with $f^*\restriction M=id_M$.

\[
\xymatrix{\ar @{} [dr] M^*
\ar[r]^{f^{**}}  &M^{**}\ar[dr]^{g}& \\
M \ar[u]^{id} \ar[r]_{id} 
& M^{''} \ar[u]_{id} \ar@2{->}[r]_{id}&M^{\bar a} 
}
\]

Now, since $M^{\bar a}$ is universal over $M''$, there exists a 
$\prec_{\K}$-mapping $g$ such that $g:M^{**}\rightarrow M^{\bar a}$ with
$g\restriction M''=id_{M''}$.  Notice that $g\circ f^{*}$ 
gives us the desired mapping of $M^*$ into $M^{\bar a}$.

\end{proof}

Notice that Lemma \ref{univ with a first} is a step closer to proving
that $\K^{am}$ satisfies Axiom \ref{DLS} of the definition of AEC as it
gives a weak downward L\"{o}wenheim-Skolem property.  It is an open
question whether or not
$\K^{am}$ is an AEC\footnote{The main difficulty is Axiom \ref{union
axiom}.}. 

An alternative  version  of Lemma \ref{univ with a
first} gives us
\begin{lemma}\label{univ ext containing a}
Given amalgamation bases of cardinality $\mu$, $M_1$ and $M_2$.  If
$M_1,M_2\prec_{\K}\C$, then there exists an amalgamation base
$M'\prec_{\K}\C$ of cardinality $\mu$ that is universal over both $M_1$
and $M_2$.

\end{lemma}
\begin{proof}
Let $\langle M'_i\mid i<\mu^+\rangle $ witness that $\C$ is a
$(\mu,\mu^+)$-limit model.  Then there exists $i<\mu^+$ such that
$M_1,M_2\prec_{\K}M'_i$.  Notice that by choice of the sequence $\langle
M'_j\mid j<\mu^+\rangle$, we have that $M'_{i+1}$ is universal over
$M'_i$.  Thus, using the assumption that $M_1$ and $M_2$ are
amalgamation bases, $M'_{i+1}$ is universal over
$M_1$ and
$M_2$.
\end{proof}
 
The following is a corollary of Theorem \ref{mu,mu+ limit is
univ} and justifies our choice of notation for
$(\mu,\mu^+)$-limit models.

\begin{corollary}\label{type aut}
 If
$\tp(a/M,\C)=\tp(b/M,\C)$, then there is an automorphism
$f$ of
$\C$ fixing $M$ such that $f(a)=b$.
\end{corollary}
\begin{proof}
Suppose that $\tp(a/M,\C)=\tp(b/M,\C)$.  
By Theorem \ref{mu,mu+ limit is univ}, $\C$ is a $(\mu,\mu^+)$-limit over
$M$.  Let $\langle M_i\in\K^{am}_\mu\mid i<\mu^+\rangle$ witness this. 
There exists an $i<\mu^+$ such that $a,b\in M_i$.  Denote $M_i$ by
both $M_a$ and
$M_b$.  By definition of types, there is a model $N$
of cardinality
$\mu$ and $\prec_{\K}$-mappings $g,h$ such that $g(a)=h(b)$ and the
following diagram commutes:
\[
\xymatrix{\ar @{} [dr] M_a
\ar[r]^{g}  &N \\
M \ar[u]^{id} \ar[r]_{id} 
& M_b \ar[u]_{h} 
}
\]

Notice that $\C$ is universal over $M_b$.  So
there is a
$\prec_{\K}$-mapping,  $f':N\rightarrow \C$ such that the following
diagram commutes:

\[
\xymatrix{\ar @{} [dr] M_a
\ar[r]^{g}  &N\ar[dr]^{f'}& \\
M \ar[u]^{id} \ar[r]_{id} 
& M_b \ar[u]_{h} \ar@2{->}[r]_{id}&\C
}
\]

Consider $f'\circ g$.  Notice that it is a partial automorphism with
domain $M_a$.  By Corollary \ref{mu,mu^+-limits are
wmh} applied to $f'\circ g(M_a)$ and $M_a$, the mapping $f'\circ g$ can be
extended to an automorphism of
$\C$, call such an extension $f$. Then,
$f\restriction M=id_M$ and $f(a)=f'\circ g(a)=f'(h(b))= b$, as
required.
\end{proof}

\bigskip
 
\section{Limit Models are Amalgamation Bases}\label{s:ab}

While Fact \ref{density of ab} 
asserts the existence of 
amalgamation bases, it is useful to identify what
other features are sufficient for   a model to be an amalgamation base.  
Makkai and Shelah were able to
prove that all existentially closed models are amalgamation bases for
$L_{\kappa,\omega}$ theories with $\kappa$ above a strongly compact
cardinal (Corollary 1.6 of
\cite{MaSh}).  
 Kolman and Shelah
identified a concept called
\emph{niceness} which implied amalgamation in categorical
$L_{\kappa,\omega}$ theories with $\kappa$ above a measurable cardinal. 
(Note: Their notion of niceness is not related to the notion of nice
towers appearing in Section \ref{s:tower}).
 They then showed that every model of
cardinality
$<\lambda$ was nice (see
\cite{KoSh}).  These results relied heavily on set theoretic
assumptions.

In a more general context, Shelah and Villaveces state that every
limit model is an amalgamation base (Fact 1.3.10 of \cite{ShVi 635}),
using
$\Diamond_{\mu^+}(S^{\mu^+}_{\cf(\mu)})$.  
For completeness, we provide a proof that
every
$(\mu,\theta)$-limit model with $\theta<\mu^+$ is an amalgamation base
under a  weaker version of diamond
$(\Phi_{\mu^+}(S^{\mu^+}_{\cf(\mu)}))$.  This is the content of Theorem
\ref{limits are ab}.

Let us first recall the set theoretic and model theoretic machinery
necessary for the proof.

\begin{definition}\index{$S^{\mu^+}_\theta}
Let $\theta$ be a regular ordinal $<\mu^+$.  We denote
$$S^{\mu^+}_\theta:=\{\alpha<\mu^+\mid\cf(\alpha)=\theta\}.$$
\end{definition}

The $\Phi$-principle defined next is known as \emph{Devlin and Shelah's
weak diamond} \cite{DS}.
\begin{definition}\label{weak diamond defn}

For $\mu$ a cardinal and $S\subseteq\mu^+$ a stationary set,
the weak diamond, denoted by $\Phi_{\mu^+}(S)$ \index{$\Phi_{\mu^+}(S)$},
is said to hold iff for all $F :
{}^{\mu^+ >} 2 \rightarrow 2$ there exists
$g: \mu^+ \rightarrow 2$ such that for every 
$f: \mu^+ \rightarrow 2$
the set 
$$\{ \delta \in S \mid F (
f \restriction\delta ) = g ( \delta ) \}\text{ is
stationary.}$$
\end{definition}

We will be using a consequence of $\Phi_{\mu^+}(S)$,
called $\Theta_{\mu^+}(S)$ (see \cite{Gr2}).

\begin{definition}
For $\mu$ a cardinal $S\subseteq\mu^+$ a stationary set,
$\Theta_{\mu^+}(S)$\index{$\Theta_{\mu^+}(S)$} is said to hold if and
only if for all families of functions
$$\{f_\eta\;:\;\eta\in\sq{\mu^+}{2}\text{ where }
f_\eta:\mu^+\rightarrow
{\mu^+}\}$$ 
and for every club $C\subseteq\mu^+$,  there exist
$\eta\neq\nu\in\sq{\mu^+}{2}$ and there exists a 
$\delta\in C\cap S$ such
that
\begin{enumerate}
\item $\eta\restriction\delta=\nu\restriction\delta$, 
\item $f_\eta\restriction\delta=f_\nu\restriction\delta$ and 
\item $\eta(\delta)\neq\nu(\delta)$.
\end{enumerate}
\end{definition}

The relative strength of these principles is provided below.  See
\cite{Gr2} for details.

\begin{fact}\index{$\Diamond_{\mu^+}(S)$}\label{diamond impl}
For $S$ a stationary subset of $\mu^+$, 
$\Diamond_{\mu^+}(S)\Longrightarrow \Phi_{\mu^+}(S)\Longrightarrow
\Theta_{\mu^+}(S)$.
\end{fact}

For most regular $\theta<\mu^+$, Fact \ref{diamond impl} and the following
imply that
$\Phi_{\mu^+}(S^{\mu^+}_\theta)$ follows from GCH:
\begin{fact}[\cite{Gy} for $\mu$ regular and \cite{Sh 108} for $\mu$
singular] For every $\mu>\aleph_1$,
 $GCH
\Longrightarrow$ $\Diamond_{\mu^+}(S)$ where $S=S^{\mu^+}_\theta$
for any regular $\theta\neq\cf(\mu)$.
\end{fact}

Thus, from $\mathrm{GCH}$ and $\Phi_{\mu^+}(S^{\mu^+}_{\cf(\mu)})$ we have
that $\Phi_{\mu^+}(S^{\mu^+}_\theta)$ holds
for every regular $\theta<\mu^+$.

In addition to the weak diamond, we will be using Ehrenfeucht-Mostowski 
models.  Let us recall some facts here.

The following gives a characterization of AECs as PC-classes.  Fact 
\ref{pres thm} is often referred to as Shelah's Presentation Theorem.

\begin{definition}\index{PC-class}
A class $\K$ of structures is called a \emph{PC-class} if there exists
a language $L_1$, a first-order theory $T_1$ in the language $L_1$ and
a collection of types without parameters, $\Gamma$, 
such that
$L_1$ is an expansion of $L(\K)$ and
$$\K=PC(T_1,\Gamma,L):=\{M\restriction L: M\models T_1\text{ and }
M\text{ omits all types from }\Gamma\}.$$
When $|T_1|+|L_1|+|\Gamma|+\aleph_0=\chi$, we say that $\K$ is $PC_\chi$. 
$PC$-classes are sometimes referred to as \emph{projective
classes}\index{projective class} or
\emph{pseudo-elementary classes}\index{pseudo-elementary class}.
\end{definition}

\begin{fact}[Lemma
1.8 of \cite{Sh 88} or see \cite{Gr2}]\label{pres thm} If
$(\K,\prec_{\K})$ is an AEC, then there exists $\chi\leq 2^{LS(\K)}$
such  that $\K$ is
$PC_\chi$.
\end{fact}

The representation of AECs as PC-classes allows us to construct\\
Ehrenfeucht-Mostowski models if there are arbitrarily large models in our
class.

\begin{definition}\label{hanf defn}\index{Hanf number}
Given an AEC $\K$, we define the \emph{Hanf number of $\K$}, abbreviated
$\Hanf(\K)$, as the minimal $\kappa$ such that for every
$PC_{2^{LS(\K)}}$-class, $\K'$, if there exists a model $M\in\K'$ of
cardinality $\kappa$, then there are arbitrarily large models in $\K'$.
\end{definition}

\begin{fact}[Claim 0.6 of \cite{Sh 394} or see \cite{Gr2}]
\label{EM models exist} Assume that $\K$ is an AEC
that contains a model of cardinality $\geq \beth_{(2^{2^{LS(\K)}})^+}$. 
Then, there is a $\Phi$, proper for linear orders\footnote{Also known as a
blueprint, see Definition 2.5 of Chapter VII, \S5 of \cite{Shc} for a
formal definition.}, such that for all linear orders
$I\subseteq J$ we have that
\begin{enumerate}
\item 
$EM(I,\Phi)\restriction L(\K)\prec_{\K}
EM(J,\Phi)\restriction L(\K)$ and
\item $\|EM(I,\Phi)\restriction L(\K)\|=|I| +LS(\K)$. 
\end{enumerate}
\end{fact}

It is theorem of
C.C. Chang based on a theorem of Morley that
$\Hanf(\K)\leq\beth_{(2^{2^{LS(\K)}})^+}$ (see Section 4 of Chapter VII
of \cite{Shc}).  Morley's proof \cite{Mo2} gives a better upper bound
in certain situations: for a class $\K$ that is $PC_{\aleph_0}$, the Hanf
number of $\K$ is $\leq \beth_{\omega_1}$.

In our context,
since $\K$ has no maximal models, $\K$ has a model of cardinality
$\Hanf(\K)$.  Then by Fact \ref{EM models exist}, we can construct
Ehrenfeucht-Mostowski models.\index{Hanf
number}\index{Ehrenfeucht-Mostowski models}\index{EM-models}

We describe an index set which appears often in papers about the
categoricity conjecture.  This index set appears in several places
including 
\cite{KoSh},
\cite{Sh 394} and \cite{ShVi 635}.  
 
\begin{notation}Let $\alpha<\lambda$ be given. 

  For $X\subseteq\alpha$, we
define $$I_X:=\big\{\eta\in\sq\omega X : \{n<\omega\mid \eta(n)\neq 0\}
\text{ is finite}\}\big\}.$$\index{$I_X$}

%
\end{notation}

The following fact is proved in several 
papers e.g. \cite{ShVi 635}.

\begin{fact}\label{EM is relatively univ}
If 
$M\prec_{\K}EM(I_\lambda,\Phi)\restriction L(\K)$ is a model of
cardinality $\mu^+$ with $\mu^+<\lambda$, then there exists
a $\prec_{\K}$-mapping 
$f:M\rightarrow EM(I_{\mu^+},\Phi)\restriction L(\K)$.
\end{fact}

A variant of this universality property is (implicit in Lemma 3.7 of 
\cite{KoSh} or see \cite{Ba}):

\begin{fact}
Suppose $\kappa$ is a regular cardinal.
If 
$M\prec_{\K}EM(I_{\kappa},\Phi)\restriction L(\K)$ is a model of
cardinality $<\kappa$ and 
$N\prec_{\K}EM(I_\lambda,\Phi)\restriction L(\K)$ is an extension 
of $M$ of cardinality $\|M\|$, then there exists a 
$\prec_{\K}$-embedding $f:N\rightarrow 
EM(I_{\kappa},\Phi)\restriction L(\K)$ such that
$f\restriction M=id_M$.
\end{fact}

We now prove  that limit models are amalgamation bases.

\begin{theorem}\label{limits are ab}\index{limit model!are amalgamation
bases}\index{amalgamation bases!limit models} Under Assumption \ref{assm
intro}, if
$M$ is a
$(\mu,\theta)$-limit for some
$\theta$ with
$\theta<\mu^+\leq\lambda$, then
$M$ is an amalgamation base.
\end{theorem}

\begin{proof}
Given $\mu$, suppose that $\theta$ is the minimal infinite ordinal 
$<\mu^+$ such that there exists a model $M$ which is a $(\mu,\theta)$-limit
and not an amalgamation base.  Notice that by Fact 
\ref{sigma and cf(sigma) limits}, we may assume that
$\cf(\theta)=\theta$.  We assume that all models have as their universe a
subset of $\mu^+$.
  
For this proof we will make use of the following notation.  We will
consider binary sequences ordered by initial segment and denote this
ordering by $\lessdot$.  For $\eta\in\sq\alpha 2$ we use $l(\eta)$ as
an abbreviation for the length of $\eta$, in this case $l(\eta)=\alpha$.

With the intention of eventually applying
$\Theta_{\mu^+}(S^{\mu^+}_\theta)$, we will define a tree
of structures $\langle M_\eta\in\K_\mu\mid\eta\in\sq{\mu^+>}2\rangle$ such
that when $l(\eta)$ has cofinality $\theta$, $M_\eta$ will be a
$(\mu,\theta)$-limit model and $M_{\eta\conc 0},M_{\eta\conc 1}$ will
witness that $M_\eta$ is not an amalgamation base.  After this tree of
structures is defined we will embed each chain of models into a universal
model of cardinality $\mu^+$.  We will apply
$\Theta_{\mu^+}(S^{\mu^+}_\theta)$ to these embeddings.
$\Theta_{\mu^+}(S^{\mu^+}_\theta)$ will provide an amalgam for
$M_{\eta\conc 0}$ and $M_{\eta\conc 1}$ over $M_\eta$ for some
sequence $\eta$ whose length has cofinality $\theta$, giving us a
contradiction.

In order to construct such a tree of models, we will need several
conditions to hold throughout the inductive construction: 
\begin{enumerate}
\item $M\preceq_{\K}M_{\langle\rangle}$
\item for $\eta\lessdot\nu\in\sq{\mu^+>}2$, $M_\eta\prec_{\K}M_\nu$
\item\label{cont condition} for $l(\eta)$ a limit ordinal with
$\cf(l(\eta))\leq\theta$,
$M_\eta=\Union_{\alpha<l(\eta)}M_{\eta\restriction\alpha}$
\item for $\eta\in\sq{\alpha}2$ with
$\alpha\in S^{\mu^+}_\theta$, 
\begin{enumerate}
\item $M_\eta$ is a $(\mu,\theta)$-limit model
\item\label{non ab} $M_{\eta\conc  \langle 0\rangle}, M_{\eta\conc
\langle 1\rangle}$ cannot be amalgamated over
$M_\eta$ 
\item $M_{\eta\conc \langle 0\rangle}$ and $M_{\eta\conc \langle
1\rangle}$ are amalgamation bases of cardinality $\mu$
\end{enumerate}

\item for $\eta\in\sq{\alpha}2$ with
$\alpha\notin S^{\mu^+}_\theta$,
\begin{enumerate}
\item $M_\eta$ is an amalgamation base

\item\label{univ cond on non theta}$M_{\eta\conc  \langle 0\rangle},
M_{\eta\conc \langle 1\rangle}$ are universal over $M_\eta$ and
\item $M_{\eta\conc  \langle 0\rangle}$ and $M_{\eta\conc \langle
1\rangle}$ are amalgamation bases of cardinality $\mu$ 
(it may be that $M_{\eta\conc  \langle 0\rangle}=M_{\eta\conc \langle
1\rangle}$ in this case).
\end{enumerate}
\end{enumerate}

This construction is possible:

\emph{$\eta=\langle\rangle$:}
By Fact \ref{density of ab}, we can find $M'\in\K^{am}_\mu$ such that
$M\prec_{\K}M'$.
Define $M_{\langle\rangle}:=M'$.

\emph{$l(\eta)$ is a limit ordinal:}
When $\cf(l(\eta))>\theta$, let 
$M'_\eta:=\Union_{\alpha<l(\eta)}M_{\eta\restriction\alpha}$.
$M'_\eta$ is not necessarily an amalgamation base, but for the purposes
of this construction, continuity at such limits is not important.  Thus
by Fact \ref{density of
ab}
we can find an extension of $M'_\eta$, say $M_\eta$, of cardinality $\mu$ 
such that $M_\eta$ is an amalgamation base.

For $\eta$ with $\cf(l(\eta))\leq\theta$, we require continuity.
Define $M_\eta:=\Union_{\alpha<l(\eta)}M_{\eta\restriction\alpha}$.
We need to verify that if $l(\eta)\notin S^{\mu^+}_\theta$, then
$M_\eta$ is an amalgamation base.  In fact, we will show that
such a $M_\eta$ will be a $(\mu,\cf(l(\eta)))$-limit model.
Let $\langle\alpha_i\mid i<\cf(l(\eta))\rangle$ be an increasing and
continuous sequence of ordinals converging to 
$l(\eta)$ such that $\cf(\alpha_i)<\theta$ for every
$i<\cf(l(\eta))$.  Condition (\ref{univ cond on non theta})
guarantees that for $i<\cf(l(\eta))$, 
$M_{\eta\restriction\alpha_{i+1}}$ is universal over
$M_{\eta\restriction\alpha}$.  Additionally, condition
(\ref{cont condition}) ensures us that
$\langle M_{\eta\restriction\alpha_i}\mid i<\cf(l(\eta))\rangle$ is
continuous.  This sequence of models witnesses that $M_\eta$ is a
$(\mu,\cf(l(\eta)))$-limit model. By our minimal choice of $\theta$ and
our assumption that in this phase of the construction
$\cf(l(\eta))\lneq\theta$, we have that
$(\mu,\cf(l(\eta)))$-limit models are amalgamation bases.  Thus $M_\eta$
is an amalgamation base.

\emph{$\eta\conc \langle i\rangle$ where $l(\eta)\in
S^{\mu^+}_\theta$:}
We first notice
that $M_\eta:=\Union_{\alpha<l(\eta)}M_{\eta\restriction\alpha}$ is a
$(\mu,\theta)$-limit model.  Why? Since
$l(\eta)\in S^{\mu^+}_\theta$ and $\theta$ is regular, we can
find an increasing and continuous sequence of ordinals,
$\langle\alpha_i\mid i<\theta\rangle$ converging to $l(\eta)$
such that for each 
$i<\theta$ we have that $\cf(\alpha_i)<\theta$.  Condition
(\ref{univ cond on non theta}) of the construction guarantees that
for each $i<\theta$, $M_{\eta\restriction\alpha_{i+1}}$ is universal
over $M_{\eta\restriction\alpha_i}$.  Thus
$\langle M_{\eta\restriction\alpha_i}\mid i<\theta\rangle$ witnesses
that $M_\eta$ is a $(\mu,\theta)$-limit model.

Since $M_\eta$ is a $(\mu,\theta)$-limit, we can fix an
isomorphism $f:M\cong M_\eta$.  By Remark \ref{f preserves ab}, $M_\eta$
is not an amalgamation base.  Thus there exist $M_{\eta\conc 0}$ and
$M_{\eta\conc 1}$ extensions of $M_\eta$ which cannot be amalgamated
over $M_\eta$.  WLOG, by the Density of Amalgamation Bases, we can choose
$M_{\eta\conc \langle 0\rangle}$ and
$M_{\eta\conc \langle 1\rangle}$ to be elements of $\K^{am}_\mu$.

\emph{$\eta\conc \langle i\rangle$ where $l(\eta)\notin
S^{\mu^+}_\theta$:}
Since $M_\eta$ is an amalgamation base, we can choose
$M_{\eta\conc \langle 0\rangle}$ and $M_{\eta\conc\langle 1\rangle}$ to be
extensions of
$M_\eta$  such that $M_{\eta\conc \langle l\rangle}\in\K^{am}_\mu$ and 
$M_{\eta\conc \langle l\rangle}$ is universal over
$M_{\eta}$, for $l=0,1$.

This completes the construction.
Let $C$ be a club containing $\{\alpha<\mu^+\mid M_\alpha\text{ has
universe }\alpha\}$.

For every $\eta\in\sq{\mu^+}2$, define 
$M_\eta:=\Union_{\alpha<\mu^+}M_{\eta\restriction\alpha}$.
Notice that by condition (\ref{univ cond on non theta}) of the
construction, each
$M_\eta$ has cardinality $\mu^+$. By categoricity in $\lambda$
and Fact \ref{EM is relatively univ}, we can fix a
$\prec_{\K}$-mapping
$g_\eta:M_\eta\rightarrow EM(I_{\mu^+},\Phi)\restriction L(\K)$
for each $\eta\in\sq{\mu^+}2$.  Now apply 
$\Theta_{\mu^+}(S^{\mu^+}_\theta)$ to find
$\eta,\nu\in\sq{\mu^+}2$ and $\alpha\in S^{\mu^+}_\theta\cap C$ such that
\begin{itemize}
\item $\rho:=\eta\restriction\alpha=\nu\restriction\alpha$,
\item $\eta(\alpha)=0$, $\nu(\alpha)=1$ and
\item $g_\eta\restriction M_\rho=g_\nu\restriction M_\rho$.
\end{itemize}


Let $N:= EM(I_{\mu^+},\Phi)\restriction L(\K)$.
Then the following diagram commutes:
\[
\xymatrix{\ar @{} [dr] M_{\rho\conc \langle 1\rangle}
\ar[r]^{g_\nu\restriction M_{\rho\conc \langle 1\rangle}} 
&N \\ M_\rho \ar[u]^{id}
\ar[r]_{id}  & M_{\rho\conc \langle 0\rangle} \ar[u]_{g_\eta\restriction
M_{\rho\conc
\langle 0\rangle}}  }
\]

Notice that $g_\eta\restriction M_{\rho\conc \langle 0\rangle}$ and
$g_\nu\restriction M_{\rho\conc \langle 1\rangle}$  witness that
$M_{\rho\conc \langle 0\rangle}$ and $M_{\rho\conc \langle 1\rangle}$ can
be amalgamated over
$M_\rho$.  Since $l(\rho)=\alpha\in S^{\mu^+}_\theta$, $M_{\rho\conc
\langle  0\rangle}$ and $M_{\rho\conc \langle 1\rangle}$ were chosen so
that they cannot be amalgamated over
$M_\rho$.  Thus, we contradict condition (\ref{non ab}) of the
construction.
\end{proof}

Now that we have verified that limit models are amalgamation bases, we
can use the existence of universal extensions to construct
$(\mu,\theta)$-limit models for arbitrary $\theta<\mu^+$.

\begin{corollary}[Existence of limit models] \label{exist
limit}\index{limit model!existence}
For every cardinal $\mu$ and limit ordinal $\theta$
with
$\theta\leq\mu^+\leq\lambda$, if $M$ is an amalgamation base of
cardinality 
$\mu$, then there exists  a
$(\mu,\theta)$-limit over
$M$.
\end{corollary}
\begin{proof}
By repeated applications of Fact \ref{exist univ} (existence of
universal extensions) and  Theorem \ref{limits are ab}.
\end{proof}

In addition to the fact that limit models are amalgamation bases, we will
use an even stronger amalgamation property of limit models. 
It
is a result of Shelah and Villaveces.  The
 argument provided is a simplification of the original and was suggested
by J. Baldwin. 

\begin{fact}[Weak Disjoint Amalgamation 
\cite{ShVi 635}]\label{wda} Given $\lambda>\mu\geq LS(\K)$ and
$\alpha,\theta_0<\mu^+$  with $\theta_0$ regular.
If $M_0$ is a $(\mu,\theta_0)$-limit and 
$M_1,M_2\in\K_\mu$ are
$\prec_{\K}$-extensions of
$M_0$, 
then for every $\bar b\in\sq\alpha(M_1\backslash M_0)$, there
exist $M_3$, a model,  
and $h$, a
$\prec_{\K}$-embedding, such that
\begin{enumerate}
\item $h:M_2\rightarrow M_3$;
\item $h\restriction M_0=id_{M_0}$ and
\item $h(M_2)\cap\bar b=\emptyset$ (equivalently $h(M_2)\cap
M_1=M_0$). 
\end{enumerate}
\end{fact}
 

\begin{proof}
Let $M_0$, $M_1$ and $M_2$ be given as in the statement of the
claim. First notice that we may assume that $M_0$, $M_1$ and $M_2$  are
such that there is a $\delta<\mu^+$ with
$M_0=M_1\cap(EM(I_{\delta},\Phi)\restriction L(\K))$ and
$M_1,M_2\prec_{\K}EM(I_{\mu^+},\Phi)\restriction L(\K)$.  Why?  Define
$\langle N_i\in\K_\mu\mid i<\mu^+\rangle$ a $\prec_{\K}$-increasing and
continuous chain of amalgamation bases such that
\begin{enumerate}
\item $N_0=M_0$ and
\item $N_{i+1}$ is universal over $N_i$.
\end{enumerate}
Let $N_{\mu^+}=\Union_{i<\mu^+}N_i$.  By categoricity and Fact \ref{EM is
relatively univ}, there exists a
$\prec_{\K}$-mapping $f$ such that $f:N_{\mu^+}\rightarrow
EM(I_{\mu^+},\Phi)\restriction L(\K)$.  Consider the club
$C=\{\delta<\mu^+\mid f(N_{\mu^+})\cap(EM(I_{\delta},\Phi)\restriction
L(\K))=f(N_\delta)\}$.  Let $\delta\in C\cap S^{\mu^+}_{\cf(\theta_0)}$. 
Notice that $f(N_\delta)$ is a $(\mu,\cf(\theta_0))$-limit model.  Since
$M_0$ is also a $(\mu,\cf(\theta_0))$-limit model, there exists
$g:M_0\cong f(N_\delta)$.  Since $f(N_{\delta+1})$ is universal over
$f(N_\delta)$, we can extend $g$ to $g'$ such that
$g':M_1\rightarrow f(N_{\delta+1})$ with
$g'(M_1)\cap EM(I_{\delta},\Phi)\restriction L(\K)=g'(M_0)$.  Thus we may
take $M_0$, $M_1$ and $M_2$ with
$M_0=M_1\cap EM(I_{\delta},\Phi)\restriction L(\K)$.

Let $\delta$ be such that $M_1\cap
(EM(I_{\delta},\Phi)\restriction
L(\K))=M_0$ and let
$\delta^*<\mu^+$ be such that
$M_1,M_2\prec_{\K}EM(I_{\delta^*})\restriction L(\K)$. Let $h$ be the $\K$
mapping from $EM(I_{\delta^*})\restriction L(\K)$ into
$EM(I_{\mu^+},\Phi)\restriction L(\K)$ induced by
$$\alpha\mapsto \delta^*+\alpha$$
for all $\alpha<\delta^*$.

We will show that if $b\in M_1\backslash M_0$ then $b\notin h(M_2)$.
Suppose for the sake of contradiction that $b\in M_1\backslash M_0$ and
$b\in h(M_2)$.  Let $\tau$ be a Skolem term and let $\bar \alpha$, $\bar
\beta$ be finite sequences such that $\bar \alpha\in I_{\delta}$ and
$\bar \beta\in I_{\delta^*}\backslash I_{\delta}$, satisfying $b=\tau(\bar
\alpha,\bar \beta)$.

Since $b\in h(M_2)$, there exists a Skolem term $\sigma$ and finite
sequences $\bar \alpha'\in I_{\delta}$ and $\bar \beta'\in
I_{\mu^+}\backslash I_{\delta^*}$ satisfying
$b=\sigma(\bar \alpha',\bar \beta')$.

Since $\bar \beta'$ and $\bar\beta$ are disjoint, we can find $\bar
\gamma'$ and $\bar \gamma\in I_{\delta}$ such that
the type of $\bar\beta'\conc\bar\beta$ is the same as the type of
$\bar \gamma'\conc\bar \gamma$ over $\bar\alpha'\conc\bar\alpha$ with
respect to the lexicographical order of $I_{\mu^+}$.
Notice then that the type of $\bar \beta'$ and $\bar\gamma'$ over
$\bar\gamma\conc\bar\alpha'\conc\bar\alpha$ are the same with respect to
the lexicographical ordering.

Recall $$EM(I_{\mu^+},\Phi)\restriction L(\K)\models
b=\tau(\bar\alpha,\bar\beta)=\sigma(\bar\alpha',\bar\beta').$$ 
Thus
$$EM(I_{\mu^+},\Phi)\restriction L(\K)\models
\tau(\bar\alpha,\bar\gamma)=\sigma(\bar\alpha',\bar\gamma').$$
Since $\bar\gamma'$ and $\bar\beta'$ look the same over
$\bar\gamma\conc\bar\alpha'\conc\bar\alpha$, we also have
$$EM(I_{\mu^+},\Phi)\restriction L(\K)\models
\tau(\bar\alpha,\bar\gamma)=\sigma(\bar\alpha',\bar\beta').$$
Combining the implications gives us a representation of $b$ with
parameters from $I_{\delta}$.  Thus $b\in EM(I_{\delta},\Phi)\restriction
L(\K)$.  Since $M_0=M_1\cap (EM(I_{\delta},\Phi)\restriction
L(\K))$, we get that $b\in M_0$ which contradicts our choice of $b$.

\end{proof}

Let us state an easy corollary of Fact \ref{wda} that will
simplify future constructions:

\begin{corollary}\label{wda improv1}
Suppose $\mu$, $M_0$, $M_1$, $M_2$ and $\bar b$ are as in the
statement of Fact \ref{wda}.  If $M_1\prec_{\K}\C$,
then there exists a $\prec_{\K}$-mapping $h$
such that
\begin{enumerate}
\item $h:M_2\rightarrow \C$,
\item $h\restriction M_0=id_{M_0}$ and
\item $h(M_2)\cap \bar b=M_0$ (equivalently $h(M_2)\cap
M_1=\emptyset$).
\end{enumerate}
\end{corollary}

\begin{proof}
By Fact \ref{wda}, there exists a $\prec_{\K}$-mapping
$g$ and a model $M_3$ of cardinality $\mu$ such that
\begin{itemize}
\item $g:M_2\rightarrow M_3$
\item $g\restriction M_0=id_{M_0}$
\item $g(M_2)\cap\bar b=M_0$ and
\item $M_1\prec_{\K}M_3$.
\end{itemize}
Since $\C$ is universal over $M_1$, we can fix
a $\prec_{\K}$-mapping $f$ such that
 $f:M_3\rightarrow \C$ and
 $f\restriction M_1=id_{M_1}$.
Notice that $h:=g\circ f$ is the desired mapping from $M_2$
into $\C$.

\end{proof}

\bigskip
\section{$\mu$-splitting}\label{s:mu-split}
Appearing in \cite{Sh 394} is $\mu$-splitting, which is
a generalization of the first-order notion of splitting to AECs.
Most results concerning $\mu$-splitting are proved under the assumption of
categoricity.  Just recently Grossberg and VanDieren have made progress
without categoricity by considering
$\mu$-splitting in Galois-stable, tame AECs (see \cite{GrVa}). 

In this section we
will develop non-$\mu$-splitting as our dependence relation and  prove
the extension and uniqueness properties for non-$\mu$-splitting types.

Before defining $\mu$-splitting we need to describe what is meant by the
image of a Galois-type:

\begin{definition}
Let $M$ be an amalgamation base and $p\in \gaS(M)$.  If $h$ is a
$\prec_{\K}$-mapping with domain $M$ we can define $h(p)$\index{$h(p)$} as
follows.  Since  $\C$ is saturated over $M$ (Corollary \ref{C is
saturated}),  we can fix
$a\in
\C$ realizing $p$.  By Proposition
\ref{mu,mu^+-limits are wmh}, we can extend $h$ to $\check h$ an
automorphism of $\C$.  Denote by 
$$h(p):=\tp(\check
h(a)/h(M)).$$  The verification that this definition does not
depend on our choices of $\check h$ and $a$ is left to the
reader.
\end{definition}

\begin{definition}\label{mu-split defn}
Let $\mu$ be a cardinal with $\mu<\lambda$.  For $M\in\K^{am}$ and
$p\in \gaS(M)$, we say that \emph{$p$ $\mu$-splits over
$N$}\index{$\mu$-splits}\index{Galois-type!$\mu$-splits}
iff
$N\prec_{\K}M$ and there exist amalgamation bases $N_1,N_2\in\K_\mu$ and a
$\prec_{\K}$-mapping $h:N_1\cong N_2$ such that
\begin{enumerate}

\item $N\prec_{\K}N_1,N_2\prec_{\K}M$,
\item $h(p\restriction N_1)\neq p\restriction N_2$
 and
\item $h\restriction N= id_N$.
\end{enumerate}
\end{definition}

\begin{remark}
If $T$ is a first-order theory stable in $\mu$ and $M$ is saturated, then
for all $N\prec M$ of cardinality $\mu$, the first-order type,
$\ftp(a/M)$, does not split (in the first-order sense) over $N$ iff
$\tp(a/M)$ does not $\mu$-split over $N$.
\end{remark}

Let us state some easy facts concerning $\mu$-splitting.
\begin{remark}

Let $N\prec_{\K}M\prec_{\K}M'$ be amalgamation bases of cardinality
$\mu$ such that $\tp(a/M')$ does not $\mu$-split over $N$.
\begin{enumerate}
\item (Monotonicity)\index{$\mu$-splits!monotonicity}
Then $\tp(a/M)$ does
not $\mu$-split over $N$.
\item (Invariance)\index{$\mu$-splits!invariance}
If $h$ is a $\prec_{\K}$-mapping with domain $M'$, $h(\tp(a/M'))$ does not
$\mu$-split over $h(N)$.
\end{enumerate}
\end{remark}

The following appears in \cite{Sh 394} under the assumption of the
amalgamation property.  The same conclusion holds in this context.
\begin{fact}
[Claim 3.3.1 of \cite{Sh 394}]
If $\K$ is $\mu$-Galois stable and $\K$ satisfies the amalgamation
property, then for every
$M\in\K_{\geq\mu}$ and every $p\in\gaS(M)$, there exists a  $N\prec_{\K}
M$ of cardinality $\mu$ such that $N\in\K$ and $p$ does not $\mu$-split
over
$N$.
\end{fact}

Shelah
and Villaveces draw connections between categoricity and
superstability properties
using $\mu$-splitting.
Let us recall some first-order consequences of superstability.

\begin{fact}\label{fo split thm}
Let $T$ be a countable first-order theory.  Suppose $T$ is superstable.
\begin{enumerate}
\item\label{fo no split chain}
If $\langle M_i\mid i\leq\sigma\rangle$ is a $\prec$-increasing and
continuous chain of models and $\sigma$ is a limit ordinal, then
for every $p\in S(M_\sigma)$, there exists $i<\sigma$ such that $p$
does not fork over $M_i$.
\item\label{fo fin char of split}
Let $T$ be a countable first-order theory.  Suppose $T$ is superstable.
Let $\langle M_i\mid i\leq\sigma\rangle$ be a $\prec$-increasing and
continuous chain of models with $\sigma$ a limit ordinal.  If $p\in
S(M_\sigma)$ is such that for every $i<\sigma$, $p\restriction M_i$ does
not fork over
$M_0$, then
$p$ does not fork over $M_0$.
\end{enumerate}
\end{fact}

These results are consequences of
$\kappa(T)=\aleph_0$\footnote{$\kappa(T)$ is the locality cardinal of
non-forking see Definition 3.1 in Chapter III \S3 of \cite{Shc} } and the
finite character of forking (see Chapter III \S 3 of
\cite{Shc}). It is interesting that Shelah and Villaveces manage to prove
analogs of these theorems without having the finite character of
$\mu$-splitting or the compactness theorem.

Fact
\ref{non-split thm} is an analog of Fact \ref{fo split thm}(\ref{fo no
split chain}), restated: under the assumption of categoricity there are no
long splitting chains.  The proof of this fact relies on a combinatorial
blackbox principle (see Chapter III of
\cite{Shg}.)

\begin{fact}[Theorem 2.2.1 from \cite{ShVi 635}]\label{non-split thm}
Under Assumption \ref{assm intro}, suppose that 
\begin{enumerate}
\item $\langle M_i\mid i\leq\sigma\rangle$ is 
$\prec_{\K}$-increasing and continuous,
\item for all $i\leq\sigma$, $M_i\in\K^{am}_\mu$,
\item for all $i<\sigma$, $M_{i+1}$ is universal over $M_i$ and
\item $p\in \gaS(M_\sigma)$.
\end{enumerate}
Then there exists an $i<\sigma$ such that
$p$ does not $\mu$-split over $M_i$.

\end{fact}

Implicit in Shelah and Villaveces' proof of Fact \ref{non-split thm} is a
statement similar to Fact
\ref{fo split thm}(\ref{fo fin char of split}). The proof of Fact
\ref{non-split thm} is by contradiction.  If Fact
\ref{non-split thm} fails to be true, then there is a counter-example
that has one of three properties (cases (a), (b), and (c) of their
proof).  Each case is separately refuted.  Case (a) yields:

\begin{fact}\label{non-split goes up}
Under Assumption \ref{assm intro}, suppose that 
\begin{enumerate}
\item $\langle M_i\mid i\leq\sigma\rangle$ is 
$\prec_{\K}$-increasing and continuous,
\item for all $i\leq\sigma$, $M_i\in\K^{am}_\mu$,
\item for all $i<\sigma$, $M_{i+1}$ is universal over $M_i$,
\item $p\in \gaS(M_\sigma)$ and
\item $p\restriction M_i$ does not $\mu$-split over $M_0$ for all
$i<\sigma$.
\end{enumerate}
Then 
$p$ does not $\mu$-split over $M_0$.

\end{fact}

The proofs of Fact \ref{non-split thm} and Fact \ref{non-split goes
up} use the full power of the categoricity assumption.  In
particular, Shelah and Villaveces use the fact that every model can be
embedded into a reduct of an Ehrenfeucht-Mostowski model.  It is open
as to whether or not the categoricity assumption can be removed:
\begin{question}
Can statements similar to Facts \ref{non-split thm} and \ref{non-split
goes up} be proved under the assumption of  any of the
working definitions of Galois superstability?
\end{question}

We now derive the extension and uniqueness properties for non-splitting
types (Theorem \ref{ext property for non-splitting} and Theorem
\ref{unique ext}). 
These results do not rely on any assumptions on the class.  We will use
these properties to find extensions of towers, but they are also useful
for developing a stability theory for tame abstract elementary classes in
\cite{GrVa}.

\begin{theorem}[Extension of non-splitting
types]\index{non-$\mu$-splitting!extension property}
\label{ext property for non-splitting}
Suppose that $M\in\K_\mu$ is universal over $N$ and
$\tp(a/M,\C)$ does not $\mu$-split over $N$, when $\C$ is a
$(\mu,\mu^+)$-limit containing
$a\Union M$.  

Let
$M'\in\K^{am}_\mu$ be an extension of $M$ with $M'\prec_{\K}\C$.
Then there exists 
a $\prec_{\K}$-mapping $g\in\Aut_M(\C)$ such that
$\tp(a/g(M'))$ does not $\mu$-split over $N$.  Equivalently,
$g^{-1}\in\Aut_{M}(\C)$ is such that
$\tp(g^{-1}(a)/M')$ does not $\mu$-split over $N$.
\end{theorem}

\begin{proof}
Since $M$ is universal over $N$, there exists a $\prec_{\K}$-mapping
$h':M'\rightarrow M$ with $h'\restriction N=id_N$.  
By Proposition \ref{mu,mu^+-limits are wmh}, we can extend $h'$ to an
automorphism $h$ of
$\C$.
Notice that by monotonicity, $\tp(a/h(M'))$ does not $\mu$-split over $N$.
By invariance, 
$$(*)\quad\tp(h^{-1}(a)/M')\text{ does not }\mu\text{-split over }N.$$

\begin{subclaim}
$\tp(h^{-1}(a)/M)=\tp(a/M)$.
\end{subclaim}
\begin{proof}
We will use the notion of $\mu$-splitting to prove this subclaim.  So let
us rename the models in such a way that our application of the definition
of $\mu$-splitting will become transparent.
Let $N_1:=h^{-1}(M)$ and $N_2:=M$.  
Let $p:=\tp(h^{-1}(a)/h^{-1}(M))$.
Consider the mapping $h:N_1\cong N_2$.
By invariance, $p$ does not $\mu$-split over $N$.
Thus, 
$h(p\restriction N_1)=p\restriction N_2$.
Let us calculate this
$$h(p\restriction N_1) = \tp(h(h^{-1}(a))/h(h^{-1}(M))) = \tp(a/M).$$
While,
$$p\restriction N_2 = \tp(h^{-1}(a)/M).$$
Thus $\tp(h^{-1}(a)/M)=\tp(a/M)$ is as required.
\end{proof}

From the subclaim, we can find a $\prec_{\K}$-mapping
$g\in\Aut_{M}(\C)$ such that $g\circ h^{-1}(a)=a$.
Notice that by applying $g$ to $(*)$ we get
$$(**)\quad\tp(a/g(M'),\C)\text{ does not }\mu\text{-split over
}N.$$
Applying $g^{-1}$ to $(**)$ gives us the \emph{equivalently} clause:
$$\tp(g^{-1}(a)/M',\C)\text{ does not }\mu\text{-split over
}N.$$ 
Since $g\restriction M=id_M$, we have that
$$\tp(g(a)/M)=\tp(g^{-1}(a)/M)=\tp(a/M).$$
\end{proof}

Not only do non-splitting extensions exist, but they are unique:
\begin{theorem}[Uniqueness of non-splitting extensions]\label{unique
ext}\index{non-$\mu$-splitting!uniqueness} Let $N,M,M'\in\K^{am}_\mu$ be
such that
$M'$ is universal over
$M$ and
$M$ is universal over $N$.  If $p\in \gaS(M)$ does not $\mu$-split over
$N$, then there is a unique $p'\in\gaS(M')$ such that $p'$ extends $p$
and $p'$ does not $\mu$-split over $N$.
\end{theorem}
\begin{proof}
By Theorem \ref{ext property for non-splitting}, there exists $p'\in
\gaS(M')$ extending $p$ such that $p'$ does not $\mu$-split over $N$. 
Suppose for the sake of contradiction that there exists $q'\neq
p'\in\gaS(M')$ extending $p$ such that $q'$ does not $\mu$-split over
$N$. Let $a, b$ be such that $p'=\tp(a/M')$ and $q'=\tp(b/M')$.
Since $M$ is universal over $N$, there exists a $\prec_{\K}$-mapping
$f:M'\rightarrow M$ with $f\restriction N=id_N$.
Since $p'$ and $q'$ do not $\mu$-split over $N$ we have
$$(*)_a\quad\tp(a/f(M'))=\tp(f(a)/f(M'))\text{ and}$$
$$(*)_b\quad\tp(b/f(M'))=\tp(f(b)/f(M')).$$
On the other hand, since $p'\neq q'$, we have that
$$(*)\quad\tp(f(a)/f(M'))\neq\tp(f(b)/f(M')).$$
Combining $(*)_a$, $(*)_b$ and $(*)$, we get
$$\tp(a/f(M'))\neq\tp(b/f(M')).$$
Since $f(M')\prec_{\K}M$, this inequality witnesses that
$$\tp(a/M)\neq\tp(b/M),$$
contradicting our choice of $p'$ and $q'$ both extending $p$.
\end{proof}

\begin{remark}\label{non-alg non-split}
Notice that the following follows from the existence and uniqueness of
non-splitting extensions:
Let $N,M,M'\in\K^{am}_\mu$ with $M$ universal over $N$ and
$M\prec_{\K}M'$.  If
$p\in\gaS(M)$ does not
$\mu$-split over $N$ and is non-algebraic, then any $q\in\gaS(M')$ which
extends $p$ and does not $\mu$-split over $N$ is also non-algebraic.
\end{remark}

The following is a corollary of the existence and uniqueness for
non-splitting types.  It allows us to extend an increasing chain of
non-splitting types.  Recall that in an AEC, a type $p$ extending an
increasing chain of types $\langle p_i\mid i<\theta\rangle$ does not
always exist and may not even be unique when it does exist (see
\cite{BaKuVa}.)
\begin{corollary}\label{chain split new cor}
Suppose that $\langle M_i\in\K^{am}_\mu\mid i<\theta\rangle$ is an
$\prec_{\K}$-increasing chain of models and $\langle p_i\in \gaS(M_i)\mid
i<\theta\rangle$ is an increasing chain of types such that for every
$i<\theta$, $p_i$ does not $\mu$-split over $M_0$ and $M_1$ is universal
over $M_0$.  If $M=\Union_{i<\theta}M_i$ is an amalgamation
base, then there exists
$p\in \gaS(M)$ such that for each $i<\theta$ $p_i\subset p$.  Moreover,
$p$ does not $\mu$-split over $M_0$.
\end{corollary}

\begin{proof}
Suppose that $M$ is an amalgamation base.
By Theorem \ref{ext property for non-splitting}, there is $p\in\gaS(M)$
extending $p_1$ such that $p$ does not $\mu$-split over $M_0$.  By 
Theorem \ref{unique ext}, we have that for
every $i<\theta$, $p_i=p\restriction M_i$.
\end{proof}

\bigskip

\section{Towers}\label{s:tower}

While 
Theorem \ref{ext property for non-splitting} allows us to find extensions
of a non-splitting Galois type in any AEC, Sections 
 \ref{s:reduced and full} and \ref{s:<b extension property} are
dedicated to the difficult task of finding non-splitting extensions of
$\alpha$-many types simultaneously under categoricity.  The mechanics
used to do this include towers.

Shelah introduced chains of towers in \cite{Sh 48} and \cite{Sh 87b}
as a tool to build a model of cardinality
$\mu^{++}$ from models of cardinality $\mu$. 
Towers are also used in \cite{BaSh} to handle abstract classes
which satisfy Axioms \ref{closure under iso axiom}-\ref{DLS} of AECs, but
for which the union axiom, Axiom
\ref{union axiom}, is not assumed.  A particular example of such
classes is the class of Banach Spaces.

We follow the notation introduced in \cite{ShVi 635}.  In \cite{ShVi
635} several other towers were defined. The superscript $c$ in the
ordering $<^c_{\mu,\alpha}$ and the superscripts $+$ and $*$ in the class
$\sq{+}{\K}^*_{\mu,\alpha}$ serve as parameters in their paper
to distinguish various definitions.  In this paper, we
will refer to only the towers in 
$\sq{+}{\K}^*_{\mu,\alpha}$ ordered by
$<^c_{\mu,\alpha}$.

\begin{definition}
\index{$\sq{+}{\K}^*_{\mu,\alpha}$}
\[{}^+\K^*_{\mu,\alpha}:=\left\{\begin{array}{ll}
(\bar M,\bar a,\bar N)
& 
\left|\begin{array}{l}
\bar M=\langle M_i\in\K_\mu\mid i<\alpha\rangle\text{ is }
\prec_{\K}\text{-increasing};\\
M_i \text{ is a }(\mu,\theta_i)\text{-limit model for some }
\theta_i<\mu^+;\\
 a_i\in M_{i+1}\backslash M_i\text{ for }i+1<\alpha;\\
\bar N=\langle N_i\in\K_\mu\mid i+1<\alpha\rangle\;\\
N_i \text{ is a }(\mu,\sigma_i)\text{-limit model for some
}\sigma_i<\mu^+;\\ 
\text{for every }i+1<\alpha,\; N_i\prec_{\K}M_i;
\\
M_i\text{ is universal over }N_i\text{ and}\\
\tp(a_i/M_i,M_{i+1})\text{ does not }\mu\text{-split over }N_i. 
\end{array}\right\}
\end{array}\right.\]

\end{definition}

\begin{remark}
The sequence $\bar M$ is not necessarily continuous.
The sequence $\bar N$ may not be
$\prec_{\K}$-increasing or continuous. 
\end{remark}

\begin{notation}
We will use the term \emph{continuous tower} to refer to towers of the
form
$(\bar M,\bar a,\bar N)$ with $\bar M$ a continuous sequence.  If $(\bar
M,\bar a,\bar N)\in\sq{+}\K^*_{\mu,\alpha}$, we say that
$\Union_{i<\alpha}M_i$ is the \emph{top of the tower} and that $(\bar
M,\bar a,\bar N)$ has \emph{length} $\alpha$.
\end{notation}

\begin{notation}
For $\theta$ a limit ordinal $<\mu^+$, we write
$\sq{+}\K^\theta_{\mu,\alpha}$ for the collection of all towers $(\bar
M,\bar a,\bar N)\in\sq{+}\K^*_{\mu,\alpha}$ where each
$M_i$ is a $(\mu,\theta)$-limit model.
\end{notation}
Our goal is to simultaneously extend the $\alpha$ non-splitting
Galois-types,
$\{\tp(a_i/M_i,M_{i+1})\mid i+1<\alpha\}$.  The
following ordering on towers captures this.

\begin{definition}\label{<c defn}\index{$<^c_{\mu,\alpha}$}
For $(\bar M,\bar a,\bar N)$ and $(\bar M',\bar a',\bar
N')\in\sq{+}\K^*_{\mu,\alpha}$, we say\\
$(\bar M,\bar a,\bar N)\leq^c_{\mu,\alpha}(\bar M',\bar a',\bar N')$
iff \begin{enumerate}
\item for $i<\alpha$ either $M'_i=M_i$ or $M'_i$ is universal over $M_i$,
\item $\bar a=\bar a'$ and
\item\label{N's equal cond} $\bar N=\bar N'$.

\end{enumerate}
We say $(\bar M,\bar a,\bar N)<^c_{\mu,\alpha}(\bar M',\bar a',\bar N')$
iff $(\bar M,\bar a,\bar N)\leq^c_{\mu,\alpha}(\bar M',\bar a',\bar N')$
and $M'_i\neq M_i$ for every $i<\alpha$.
\end{definition}
 
\begin{remark}
Notice that in Definition \ref{<c defn}, 
if $(\bar M,\bar a,\bar N)<^c_{\mu,\alpha}(\bar M',\bar a,\bar N)$, then
for every $i<\alpha$,
$\tp(a_i/M'_i,M'_{i+1})$ does not $\mu$-split over $N_i$.
\end{remark}

\begin{notation}
We will often be looking at extensions of an initial segment of a tower. 
We introduce the following notation for this.  Suppose $(\bar M,\bar
a,\bar N)\in\sq{+}\K^*_{\mu,\alpha}$.  Let $\beta<\alpha$.  We write
$\bar M\restriction \beta$ for the sequence $\langle M_i\mid
i<\beta\rangle$.  Similarly, $\bar a\restriction\beta=\langle a_i\mid
i+1<\beta\rangle$ and $\bar N\restriction\beta=\langle N_i\mid
i+1<\beta\rangle$.
Then
$(\bar M,\bar a,\bar N)\restriction\beta$ will represent the tower
$(\bar M\restriction\beta,\bar
a\restriction\beta,\bar N\restriction\beta)\in\sq{+}\K^*_{\mu,\beta}$. 
If $(\bar M',\bar a',\bar N')$ is a $<^c_{\mu,\beta}$-extension of $(\bar
M,\bar a,\bar N)\restriction \beta$, we refer to it as a \emph{partial
extension} of $(\bar M,\bar a,\bar N)$.
\end{notation}

The requirement that $M'_i$ is universal over $M_i$ in the definition of
$<^c_{\mu,\alpha}$ allows us to conclude that the models in the
union of a $<^c_{\mu,\alpha}$-increasing chain of towers are limit models.
In particular, the union of a $<^c_{\mu,\alpha}$-increasing chain (of
length $<\mu^+$) of towers is a tower.

\begin{definition}
We say that $\K$ satisfies the $<^c_{\mu,\alpha}$-extension property iff
every tower in $\sq{+}\K^*_{\mu,\alpha}$ has a
$<^c_{\mu,\alpha}$-extension.
\end{definition}

The
$<^c_{\mu,\alpha}$-extension property serves as a weak
substitute for the extension property of non-forking in first-order
model theory, but is much stronger than the extension property for
non-splitting.  Notice that for towers with
$\alpha=1$, Theorem \ref{ext property for non-splitting} and the
existence of universal extensions (Fact \ref{exist
univ})
 give the $<^c_{\mu,1}$ extension property.  
Actually it is
possible to derive the $<^c_{\mu,n}$-extension property for all
$n\leq\omega$ with no more than the existence of universal extensions and
the extension property for non-splitting types.

It is open whether or not every $\K$ satisfying Assumption \ref{assm
intro} has the
$<^c_{\mu,\alpha}$-extension property for $\alpha>\omega$.  The
difficulties concern discontinuous towers.  Notice that if $(\bar M,\bar
a,\bar N)$ is not continuous, then for some limit ordinal $i<\alpha$, we
may have that $\Union_{j<i}M_j$ is not an amalgamation base.  Suppose
that we have constructed a partial extension of $(\bar M,\bar a,\bar N)$
up to $i$.  It may be the case that this extension and $M_i$ may not be
amalgamated over
$\Union_{j<i}M_j$.  This would rule out much hope for using the partial
extension as a base for a longer extension of the entire tower $(\bar
M,\bar a,\bar N)$.  With this in mind, it is natural to restrict
ourselves to continuous towers.  However, it is not easy to prove that
every continuous tower has a continuous extension.  In fact, we can only
prove this under an extra assumption, Hypothesis 1, (see Section
\ref{s:reduced and full}).

In addition to the continuous towers, we have identified two subclasses of
$\sq{+}\K^*_{\mu,\alpha}$, amalgamable and nice towers, for which a
 $<^c_{\mu,\alpha}$-extension property can be proven.

\begin{definition}\label{nice-notation}
We say that $(\bar M,\bar a,\bar N)\in\sq{+}{\K}^*_{\mu,\alpha}$ is
\emph{nice} iff whenever $i<\alpha$ is a limit ordinal,
$\Union_{j<i}M_j$ is an amalgamation base.
\end{definition}

\begin{remark}
 Since every $M_i$ is a $(\mu,\theta_i)$-limit for some
limit ordinal $\theta_i<\mu^+$, by Theorem \ref{limits are ab}, we have
that every $M_i$ is also an amalgamation base.  So \emph{nice}
only is a requirement for limit ordinals $i$ where
$\bar M$ is not continuous at $i$. 
Thus, if $(\bar M,\bar a,\bar N)$ is a continuous tower, then $(\bar
M,\bar a,\bar N)$ is nice.

\end{remark}

Notice that the definition of nice does not require that the top of
the tower ($\Union_{i<\alpha}M_i$) be an amalgamation base.
For these towers we introduce the terminology:
\begin{definition}\label{amalgamable defn}
We say that $(\bar M,\bar a,\bar N)\in\sq{+}{\K}^*_{\mu,\alpha}$ is
\emph{amalgamable} iff it is nice and $\Union_{i<\alpha}M_i$ is an
amalgamation base.
\end{definition}

We use the word amalgamable to refer to such towers, because any two
$<^c_{\mu,\alpha}$-extensions of an amalgamable tower $(\bar M,\bar
a,\bar N)$ can be amalgamated over $\Union_{i<\alpha}M_i$.

 Notice that the classes of amalgamable and nice towers both avoid the
problematic towers described above.
 Namely, if $(\bar M,\bar a,\bar N)$ is discontinuous at $i$, we
require that $\Union_{j<i}M_j$ is an amalgamation base.  We can show that
every nice tower has an amalgamable extension (Theorem \ref{partial
ext}).  In particular, every continuous tower has an amalgamable
extension.  However, this amalgamable extension may not be continuous. 
Furthermore, if we instead restrict ourselves to amalgamable towers, we
will run into the difficulty that the union of a
$<^c_{\mu,\alpha}$-increasing chain of amalgamable towers need not be
amalgamable (or even nice).  But, with a little help from Hypothesis 1,
we are able to carry through the strategy of restricting ourselves to
continuous towers. By carefully stacking the amalgamable extensions in
Section
\ref{s:reduced and full},
we construct continuous extensions of continuous
towers.

\begin{notation}
Ultimately, we will be constructing a $<^c_{\mu,\alpha}$-extension,
$(\bar M',\bar a',\bar N')$ of a tower $(\bar M,\bar a,\bar N)$, but we
will allow the extension to live on a larger index set, $(\bar M',\bar
a',\bar N')\in\sq{+}\K^*_{\mu,\alpha'}$ for some $\alpha'>\alpha$.  We
will also like to arrange the construction so that $\alpha$ is not
identified with an initial segment of $\alpha'$, but as some other
scattered, increasing subsequence of $\alpha'$.  Therefore, we will
prefer to consider the relative structure of these index sets in addition
to their order types.  We make the following convention that will be
justified in later constructions.
When $\alpha$ and $\delta$ are ordinals, $\alpha\times\delta$ with the
lexicographical ordering ($<_{lex}$), is well ordered.  
Recall that $\otp(\alpha\times\delta,<_{lex})=\delta\cdot\alpha$ where 
$\cdot$ is ordinal multiplication.
For easier notation in future constructions, we will identify
$\alpha\times\delta$ with the interval of ordinals
$[0,\delta\cdot\alpha)$ and $\sq{+}\K^*_{\mu,\alpha\times\delta}$ will
refer to the collection of towers $\sq{+}\K^*_{\mu,\delta\cdot\alpha}$.
The notation will be more convenient when we compare towers in 
$\sq{+}\K^*_{\mu,\alpha\times\delta}$ with those in
$\sq{+}\K^*_{\mu,\alpha'\times\delta'}$ for $\alpha'\geq\alpha$ and
$\delta'\geq\delta$.
\end{notation}

We will make use of the following proposition concerning
$<^c_{\mu,\alpha}$ throughout the paper
\begin{proposition}\label{monotonicity of universality}
If $(\bar M',\bar a,\bar N)$ is a $<^c_{\mu,\alpha}$-extension of $(\bar
M,\bar a,\bar N)$, then for every $i\leq j<\alpha$, we have that $M'_j$
is universal over $M_i$.
\end{proposition}
\begin{proof}
By definition of $<^c_{\mu,\alpha}$, we have that $M'_i$ is universal
over $M_i$.  Since $\bar M'$ is increasing, $M'_i\preceq_{\K}M'_j$.  So
$M'_j$ is universal over $M_i$ as well.

\end{proof}

\bigskip

\part{Uniqueness of Limit Models}\label{p:uniqueness}

We will use towers to prove the uniqueness of
limit models by producing a model which is simultaneously a
$(\mu,\theta_1)$-limit model and a $(\mu,\theta_2)$-limit model.  The
construction of such a model is sufficient to prove the uniqueness of
limit models by Fact
\ref{unique limits} and involves building an increasing and continuous
chain of towers.

The idea is to build a two-dimensional
array (with the cofinality of the height $=\theta_1$ and the
cofinality of the width $=\theta_2$) of models such that the bottom
corner of the array ($M^*$) is a
$(\mu,\theta_1)$-limit model witnessed by the last column and a
$(\mu,\theta_2)$-limit model witnessed by the last row of the array.
The actual construction involves increasing the length of the towers
as we go from one row to the next.


The construction of this array is done by identifying each row of the
array with a tower and then building a $<_{\mu,\alpha}^c$-increasing and
continuous chain of towers (where $\alpha$ will vary throughout our
construction).

$$\tiny
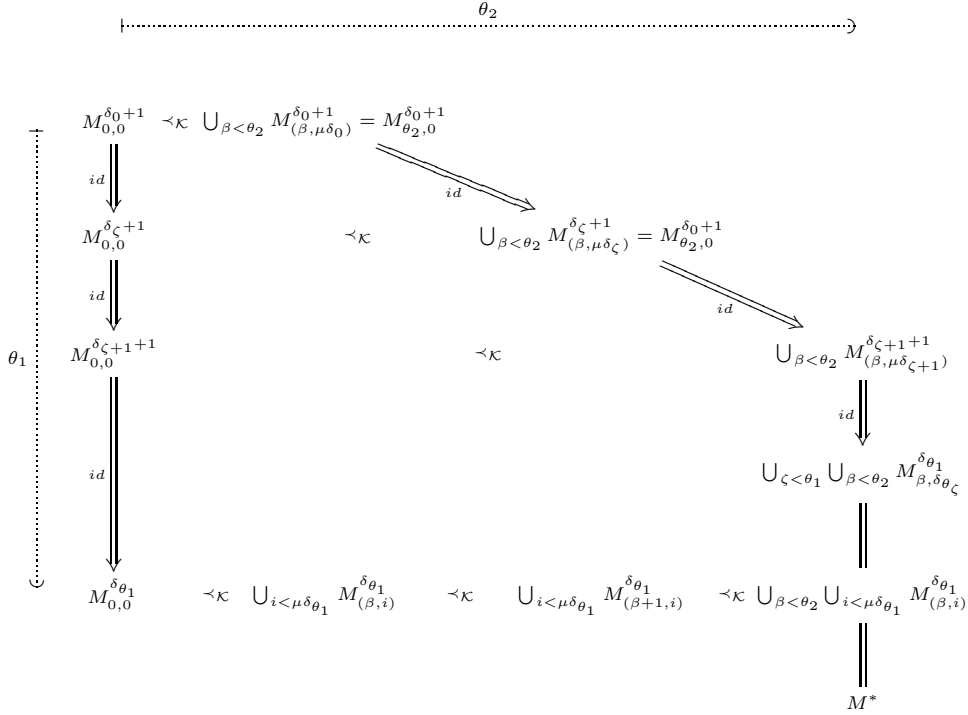
{
\xymatrix@C6pt{{}&{}\ar@{|.)}[rrr]^{\txt{$\theta_2$}}&&&{}\\
{}\ar@{|.)}[dddd]_{\txt{$\theta_1$}}
&M^{\delta_0+1}_{0,0}\ar@{}[r]|(.3){\prec_{\K}}\ar@2{->}[d]_{id}
&\;\;\Union_{\beta<\theta_2}M_{(\beta,\mu\delta_0)}^{\delta_0+1}=
M^{\delta_0+1}_{\theta_2,0}
\ar@2{->}[dr]_{id}&&\\
&M^{\delta_\zeta+1}_{0,0}\ar@{}[rr]|{\prec_{\K}}
\ar@2{->}[d]_{id}&
&\Union_{\beta<\theta_2}M_{(\beta,\mu\delta_\zeta)}^{\delta_\zeta+1}=
M^{\delta_0+1}_{\theta_2,0}\ar@2{->}[dr]_{id}
&\\
&
M^{\delta_{\zeta+1}+1}_{0,0}\ar@{}[rrr]|{\prec_{\K}}
\ar@2{->}[dd]_{id}&&
&\Union_{\beta<\theta_2}M_{(\beta,\mu\delta_{\zeta+1})}^{\delta_{\zeta+1}+1}
\ar@2{->}[d]_{id}\\
&&&&
\Union_{\zeta<\theta_1}\Union_{\beta<\theta_2}
M^{\delta_{\theta_1}}_{\beta,\delta_{\theta_\zeta}}\ar@{=}[d]\\
{}&M_{0,0}^{\delta_{\theta_1}} \ar@{}[r]|{\prec_{\K}}
&\;\;\Union_{i<\mu\delta_{\theta_1}}
M_{(\beta,i)}^{\delta_{\theta_1}}
\ar@{}[r]|{\prec_{\K}}
&
\Union_{i<\mu\delta_{\theta_1}}
M_{(\beta+1,i)}^{\delta_{\theta_1}}\ar@{}[r]|{\prec_{\K}}
& \Union_{\beta<\theta_2}\Union_{i<\mu\delta_{\theta_1}}
M_{(\beta,i)}^{\delta_{\theta_1}} \\
&&&&M^*\ar@{=}[u]
}
}
$$

In order to witness that $M^*$ is
a
$(\mu,\theta_1)$-limit model, we will need for our towers to be
\emph{increasing} in such a way that the models in the $\delta+1^{st}$
tower are universal over the models in the
$\delta^{th}$ tower.  This is possible if we can prove that every
continuous tower has a continuous
$<^c_{\mu,\alpha}$-extension.  This is the subject of Section
\ref{s:reduced and full} and related material appears in Section
\ref{s:<b extension property}.

While $M^*$ is built up by a chain of cofinality $\theta_2$, it may not
be a $(\mu,\theta_2)$-limit model.  
In order to conclude that
$M^*$ is a $(\mu,\theta_2)$-limit model, we show in Section \ref{s:full},
that the top of a continuous, relatively full tower of length $\theta_2$
is a
$(\mu,\theta_2)$-limit model.  
We will construct the relatively full tower by requiring that at every
stage of our construction of the array, we realize all the strong types
over the previous tower in a systematic way.  Section
\ref{s:dense towers} provides the technical machinery to carry this
through.  The actual construction of $M^*$ is carried out in Section
\ref{s:unique limits}.

\bigskip

\section{Relatively Full Towers}\label{s:full}

We begin this section by recalling the definition of \emph{strong types}
from \cite{ShVi 635}.
\begin{definition}[Definition 3.2.1 of \cite{ShVi 635}]\label{strong type
defn} For $M$ a $(\mu,\theta)$-limit model, \index{strong
types}\index{Galois-type!strong}\index{$\St(M)$}\index{$(p,N)$}
\begin{enumerate}
\item Let
$$\St(M):=\left\{\begin{array}{ll}
(p,N)
& 
\left|\begin{array}{l}
N\prec_{\K}M;\\
N\text{ is a }(\mu,\theta)-\text{limit model};\\
M\text{ is universal over }N;\\
p\in \gaS(M)\text{ is non-algebraic};\\
\text{and }p\text{ does not }\mu-\text{split over }N.
\end{array}\right\}
\end{array}\right .
$$ 
\item\label{sim defn}
For types $(p_l,N_l)\in\St(M)$ ($l=1,2$),
we say $(p_1,N_1)\sim (p_2,N_2)$ iff for every
$M'\in\K^{am}_\mu$ extending $M$ \index{$\sim$}
there is a $q\in \gaS(M')$ extending both $p_1$ and $p_2$ such that
$q$ does not $\mu$-split over $N_1$ and $q$ does not $\mu$-split over
$N_2$.
\end{enumerate}
\end{definition}

\begin{notation}
Suppose $M\prec_{\K}M'$ are amalgamation bases of cardinality $\mu$.
For $(p,N)\in \St(M')$, if $M$ is universal over $N$, we define the
restriction
$(p,N)\restriction M\in \St(M')$\index{$(p,N)\restriction M$} to be
$(p\restriction M,N)$.

We write $(p,N)\restriction M$ only when $p$ does not $\mu$-split
over $N$ and $M$ is universal over
$N$. 
\end{notation}

Notice that $\sim$ is an equivalence relation on $\St(M)$.  To see that
$\sim$ is a transitive relation on $\St(M)$, suppose that $(p_1,N_1)\sim
(p_2,N_2)$ and $(p_2,N_2)\sim(p_3,N_3)$.  Let $M'\in\K^{am}_\mu$ be an
extension of $M$ and fix $q_{ij}\in\gaS(M')$ extending both $p_i$ and
$p_j$ and $q_{ij}$ does not $\mu$-split over both $N_i$ and $N_j$ (for
$\langle i,j\rangle=\langle 1,2\rangle,
\langle 2,3\rangle$).  Since $p_2$ has a unique non-splitting extension
to $M'$ (Theorem \ref{unique
ext}), we know that
$q_{12}=q_{23}$.   Then $q_{12}$ witnesses that $(p_1,N_1)\sim (p_3,N_3)$
since it is an extension of both $p_1$ and $p_3$ and does not $\mu$-split
over both
$N_1$ and $N_3$.

The following lemma is used to provide a bound on the number of strong
types.

\begin{lemma}\label{enough to consider M'}Given $M\in\K^{am}_\mu$, and
$(p,N),(p',N')\in\St(M)$.  Let
$M'\in\K^{am}_\mu$ be a universal extension of $M$.  To show that
$(p,N)\sim (p',N')$ it suffices to find $q\in \gaS(M')$ such that
$q$ extends both $p$ and $p'$ and such that $q$ does not $\mu$-split over
$N$ and
$N'$.
\end{lemma}
\begin{proof}
Suppose $q\in \gaS(M')$ extends both $p$ and $p'$ and does not
$\mu$-split over $N$ and $N'$.  Let
$M^*\in\K^{am}_\mu$ be an extension of
$M$.  By universality of $M'$, there exists $f:M^*\rightarrow M'$ such
that
$f\restriction M=id_M$.  Consider $f^{-1}(q)$.  It extends $p$ and $p'$
and does not $\mu$-split over $N$ and $N'$ by invariance.  Thus
$(p,N)\sim (p',N')$.
\end{proof}

The following appears as a Fact 3.2.2(3) in \cite{ShVi 635}.  We provide
a proof here for completeness.

\begin{fact}\label{St small}
For $M\in\K^{am}_\mu$, $|\St(M)/\sim|\leq\mu$.
\end{fact}

\begin{proof}[Proof of Fact \ref{St small}]

Suppose for the sake of contradiction that \\
$|\St(M)/\sim|>\mu$.
  
Let $\{(p_i,N_i)\in\St(M)\mid i<\mu^+\}$ be pairwise non-equivalent.  By
Galois-stability (Fact \ref{cat implies stab}) and the pigeon-hole
principle, there exist
$p\in \gaS(M)$ and
$I\subset\mu^+$ of cardinality $\mu^+$ such that $i\in I$ implies
$p_i=p$.  Set
$p:=\tp(a/M)$ with $a\in\C$.

Fix $M'\in\K^{am}_\mu$ a universal extension of $M$ inside $\C$. 
We will show that there are $\geq\mu^+$ types over $M'$.  This will
provide us with a contradiction since $\K$ is Galois-stable in $\mu$ (Fact
\ref{cat implies stab}).

For each
$i\in I$, by the extension property of non-splitting (Theorem
\ref{ext property for non-splitting}), there exists $f_i\in\Aut_M(\C)$
such that 
\begin{itemize}
\item
$\tp(f_i(a)/M')\text{ does not }\mu\text{-split over }N_i$
and 
\item$\tp(f_i(a)/M')$ extends $\tp(a/M)$.
\end{itemize}
\begin{claim}\label{getting different types claim}
 For $i\neq j\in I$, we have that the types,
$\tp(f_i(a)/M')$ and $\tp(f_j(a)/M')$, are not equal. 
\end{claim}
\begin{proof}[Proof of Claim \ref{getting different types claim}]

Otherwise $\tp(f_i(a)/M')$ does not $\mu$-split over $N_i$ and does not
$\mu$-split over $N_j$.  By Lemma \ref{enough to consider M'}, this
implies that
$(p,N_i)\sim(p,N_j)$ contradicting our choice of non-$\sim$-equivalent
strong types.
\end{proof}
This completes the proof as $\{\tp(f_i(a)/M')\mid i\in I\}$ is a set of
$\mu^+$ distinct types over $M'$, contradicting $\mu$-Galois-stability.

\end{proof}

\medskip
We can now consider
towers which are saturated with respect to strong types (from
$\St(M)$). These towers are called relatively full.

\begin{definition}\label{full defn}\index{full relative to}
Let $\alpha,\delta$ and $\theta$ be limit ordinals $<\mu^+$. 
Suppose $\langle\bar M_{\beta,i}\mid(\beta,i)\in\alpha\times\delta\rangle$
is such that each $\bar M_{\beta,i}$ is a sequence of limit models, 
$\langle M^\gamma_{\beta,i}\mid \gamma<\theta\rangle$,
with
$M^{\gamma+1}_{\beta,i}$ universal over $M^\gamma_{\beta,i}$ for
all $(\beta,i)\in\alpha\times\delta$.

A tower $(\bar M,\bar a,\bar
N)\in\sq{+}\K^\theta_{\mu,\alpha\times\delta}$ is said to be \emph{full
relative to
$\langle\bar M^\gamma\mid\gamma<\theta\rangle$} iff for all
$(\beta,i)\in\alpha\times\delta$
\begin{enumerate}
\item $\bar M_{\beta,i}$ witnesses that $M_{\beta,i}$ is a
$(\mu,\theta)$-limit model and
\item for all
$(p,N^*)\in\St(M_{\beta,i})$ with $N^*=M^\gamma_{\beta,i}$ for some
$\gamma<\theta$, 
there is a $j<\delta\text{ such that }
(\tp(a_{\beta+1,j}/M_{\beta+1,j}),N_{\beta+1,j})\restriction
M_{\beta,i}\sim (p,N^*).$
\end{enumerate}
\end{definition}

\[\xymatrix@C=75pt@!R@!C{
M^0_{0,0}\ar@{}[r]|*+{\prec_{\K}\dots\prec_{\K}} \ar@2{->}[d]^{id}
& M^0_{\beta,i}\ar@{}[r]|*+{\prec_{\K}\dots}\ar@2{->}[d]^{id}&\\
M^\gamma_{0,0}
\ar@{}[r]|*+{\prec_{\K}\dots\prec_{\K}}\ar@2{->}[d]^{id}&
M^\gamma_{\beta,i}\ar@{}[r]|*+{\prec_{\K}\dots}\ar@2{->}[d]^{id}&\\
M^{\gamma+1}_{0,0}
\ar@{}[r]|*+{\prec_{\K}\dots\prec_{\K}}\ar@2{->}[d]^{id}&
M^{\gamma+1}_{\beta,i}\ar@{}[r]|*+{\prec_{\K}\dots}
\ar@2{->}[d]^{id}&
\\  M_{0,0} = \Union_{\gamma<\theta}M^\gamma_{0,0}
\ar@{}[r]|*+{\prec_{\K}\dots\prec_{\K}}&
M_{\beta,i}=\Union_{\gamma<\theta}M^\gamma_{\beta,i}\ar@{}[r]|*+{\prec_{\K}\dots}&
}\]

\begin{notation}\index{relatively full}
We say that $(\bar M,\bar a,\bar N)\in\K^\theta_{\mu,\alpha\times\delta}$
is
\emph{relatively full} iff there  exists $\langle\bar
M_{\beta,i}\mid(\beta,i)\in\alpha\times\delta\rangle$ as in Definition
\ref{full defn} such that $(\bar M,\bar a,\bar N)$ is full relative to
$\langle\bar
M_{\beta,i}\mid(\beta,i)\in\alpha\times\delta\rangle$.
\end{notation}

\begin{remark}
A strengthening  of Definition \ref{full defn} appears 
in \cite{ShVi 635} under the name full towers (see Definition 3.2.3 of
their paper).   Consider the statement:
$$(*)\quad\forall
M\in\K^{am}_\mu\text{ and }\forall(p,N),(p',N')\in\St(M),\;
(p,N)\sim(p',N')\text{ iff }p=p'.$$ 

Notice that for $M\in\K^{am}_\mu$, if $(p,N)\sim (p',N')\in\St(M)$, then
necessarily
$p=p'$.  To see this, take $M'\in\K^{am}_\mu$ some extension of $M$ and
$q\in\gaS(M')$ such that
$q$ extends both $p$ and $p'$ and does not $\mu$-split over $N$ and
$N'$.  Then $q\restriction M=p$ and $q\restriction M=p'$.  So $p$ and
$p'$ must be equal.
\end{remark}

However we do not know that $(*)$ holds in our context.  Shelah
has implicitly shown, with much work, that it does hold in categorical
AECs which satisfy the amalgamation property \cite{Sh 394}.  It is a
consequence of transitivity of non-splitting.

Property $(*)$ implies that relatively full towers
are full.  We use relatively full towers  since the construction of full 
towers by an increasing chain of
towers in this context has been seen
to be problematic.

The following proposition is immediate from the definition of relative
fullness.

\begin{proposition}\label{restriction of full}
Let $\alpha$ and $\delta$ be limit ordinals $<\mu^+$.  If $(\bar M,\bar
a,\bar N)\in\sq{+}\K^\theta_{\mu,\alpha\times\delta}$ is full relative to
$\langle\bar
M_{\beta,i}\mid(\beta,i)\in\alpha\times\delta\rangle$, then for every limit
ordinal
$\beta<\alpha$, we have that the restriction
$(\bar M,\bar a, \bar N)\restriction \beta\times\delta$ is full
relative to $\langle \bar
M_{\beta',i'}\mid(\beta',i')\in\beta\times\delta\rangle$.
\end{proposition}

The following theorem is proved in \cite{ShVi 635} for full towers
(Theorem 3.2.4 of their paper).  Our strengthening provides us with an
alternative characterization of limit models as the top of a relatively
full tower.

\begin{theorem}\label{full is limit}
Let $\alpha$ be an ordinal $<\mu^+$ such that $\alpha=\mu\cdot\alpha$.
Suppose $\delta<\mu^+$. If
$(\bar M,\bar a,\bar N)\in\sq{+}\K^\theta_{\mu,\alpha\times\delta}$ is
full relative to
$\langle\bar
M_{\beta,i}\mid(\beta,i)\in\alpha\times\delta\rangle$
and 
$\bar M$ is continuous, then $M:=\Union_{i<\alpha\cdot\delta}M_i$ is a
$(\mu,\cf(\alpha))$-limit model over $M_0$.
\end{theorem}

\begin{proof}
Let $M'\prec_{\K}\C$ be a
$(\mu,\alpha)$-limit over
$M_{0,0}$ witnessed by
$\langle M'_{i}\mid i<\alpha\rangle$.  By Weak Disjoint Amalgamation and
renaming elements, we can arrange that $\Union_{i<\alpha}M'_i\cap
\Union_{i<\alpha\cdot\delta}M_i
=M_{0,0}$ and that  for each $i<\alpha$ we can identify the universe of
$M'_i$ with $\mu(1+i)$.  Notice that since $\alpha=\mu\cdot\alpha$, we
have that
$i\in M'_{i+1}$ for every $i<\alpha$.
We will construct an
isomorphism from
$M$ into $M'$.

Now we define by induction on $i<\alpha$ a increasing and
continuous sequence of 
$\prec_{\K}$-mappings
$\langle h_i\mid i<\alpha\rangle$ such that
\begin{enumerate}
\item\label{hi cond} $h_i:M_{i,j}\rightarrow M'_{i+1}$ for some $j<\delta$
\item $h_0=id_{M_{0,0}}$ and
\item\label{put it all in} $i\in \rg(h_{i+1})$.
\end{enumerate}

For $i=0$ take $h_0=id_{M_{0,0}}$.  For $i$ a limit ordinal let
$\check h_i=\Union_{i'<i}h_{i'}$.  Since $\bar M$ is continuous, we know
that
$\Union_{\stackrel{i'<i}{j<\delta}}M_{i',j}$ is an amalgamation base.
Thus the induction hypothesis gives us that $h_i$ is a
$\prec_{\K}$-mapping from
$M_{i,0}=\Union_{\stackrel{i'<i}{j<\delta}}M_{i',j}$ into
$M'_i$ allowing us to satisfy condition (\ref{hi cond}) of the
construction.

Suppose that $h_i$
has been defined.  Let $j<\delta$ be such that $h_i:M_{i,j}\rightarrow
M'_{i+1}$.   There are two cases:  either $i\in\rg(h_i)$ or
$i\notin\rg(h_i)$. First suppose that $i\in\rg(h_i)$.  Since $M'_{i+2}$
is universal over
$M'_{i+1}$, it is also universal over $h_i(M_{i,j})$.  This allows us to
extend $h_i$ to $h_{i+1}:M_{i+1,0}\rightarrow M'_{i+2}$.

Now consider the case when $i\notin\rg(h_i)$.  We illustrate the
construction for this case:
\[\xymatrix@C10pt{
&&&i\ar@{}[d]|*+{\in}\ar@/^2pc/@{{|-}.>}[ddrrr]^{\check f_i}&&&\\
M'_0 \ar@2{->}[r]|{\dots}^{id}& M'_i \ar@2{->}[r]^{id}
& M'_{i+1}\ar@2{->}[r]^{id} & 
M'_{i+2}\ar@/^/@{.>}[drr]^{\check f_i}&&\\
&&&&&\check M &f_i(i)=f_a(h'(a_{i+1,j'}))\\
M_{0,0}\ar@{}[r]|*+{\prec_{\K}}\ar[uu]_{id}^{h_0} &
M_{i,0}\ar@{}[r]|*+{\prec_{\K}} & M_{i,j} \ar[uu]^{h_i}
\ar@{}[r]|*+{\prec_{\K}}&
M_{i+1,0}\ar@{}[r]|*+{\prec_{\K}} & 
M_{i+1,j'} \ar@/_/@{.>}[ur]_{\check f_a\circ h'}
\ar@{.>}[luu]^{h_{i+1}}_(0.4){\check f_i^{-1}\circ\check
f_a\circ h'} &&
\\
 &&&a_{i+1,j'}\ar@{}[u]|*+{\in}
\ar@/_2pc/@{{|-}.>}[uurrr]_{\check f_a\circ h'}&&&
}
\]

%

 Since $\langle M^\gamma_{i,j}\mid
\gamma<\theta\rangle$ witness that $M_{i,j}$ is a $(\mu,\theta)$-limit
model, by Fact \ref{non-split thm}, there exists $\gamma
<\theta$ such that
$\tp(i/M_{i,j})$ does not $\mu$-split over
$M^\gamma_{i,j}$.
By our choice of $\bar M'$ disjoint from $\bar M$ outside of $M_0$, we
know that $i\notin M_{i,j}$.  Thus $\tp(i/M_{i,j})$ is non-algebraic and
by relative fullness of
$(\bar M,\bar a,\bar N)$, there exists
$j'<\delta$ such that
$$(\tp(i/M_{i,j}),M^\gamma_{i,j})\sim(\tp(a_{i+1,j'}/M_{i+1,j'}),
N_{i+1,j'})\restriction M_{i,j}.$$ 
In particular we have that
$$(*)\quad\tp(a_{i+1,j'}/M_{i,j})=\tp(i/M_{i,j}).$$

We can extend $h_i$ to an automorphism $h'$ of $\C$.
An application of $h'$ to $(*)$ gives us

$$(**)\quad\tp(h'(a_{i+1,j'})/h_i(M_{i,j}))=
\tp(i/h_i(M_{i,j})).$$

 By $(**)$, there exist $M^*\in\K^{am}_\mu$ a $\K$-substructure of
$\C$ containing $M_{i,j}$ and $\prec_{\K}$-mappings
$f_a:h'(M_{i+1,j'+1})\rightarrow M^*$ and $f_i:M'_{i+2}\rightarrow M^*$
such that
$f_a(h'(a_{i+1,j'}))=f_i(i)$ and $f_a\restriction
h_i(M_{i,j})=f_i\restriction h_i(M_{i,j})=id_{h_i(M_{i,j})}$.
Since  $M'_{i+2}$ is universal over $M'_{i+1}$, it is also universal over
$h_i(M_{i,j})$.  So we may assume that $M^*=M'_{i+2}$.
Since
$\C$ is a $(\mu,\mu^+)$-limit model, we can extend $f_a$ and
$f_i$ to automorphisms of $\C$, say $\check f_a$ and $\check f_i$.
Let
$h_{i+1}:M_{i+1,j'+1}\rightarrow M'_{i+2}$ be defined as $\check
f_i^{-1}\circ \check f_a\circ h'$.  Notice that $h_{i+1}(a_{i+1,j'})=i$.
This completes the construction.

Let $h:=\Union_{i<\alpha}h_i$.  Clearly
$h:M\rightarrow M'$.  To see that $h$ is an
isomorphism, notice that condition (\ref{put it all in}) of the
construction forces $h$ to be surjective.

\end{proof}

\begin{remark}
Theorem \ref{full is limit} can be improved by replacing the assumption
of continuity of
$(\bar M,\bar a,\bar N)$ with niceness.  The same proof works with a
minor adjustment at the limit stage.  We lift the requirement that
$\langle h_i\mid i<\alpha\rangle$ is continuous and use the fact that
$M'_{i+1}$ is universal over $M'_i$ to carry out the construction at
limits.
\end{remark}

\bigskip

\section{Existence of Continuous
$<_{\mu,\alpha}^c$-extensions}\label{s:reduced and full}

 Our proof of the uniqueness of limit models will
involve a $<^c_{\mu,\alpha}$-increasing chain of continuous towers such
that the index sets of the towers grow throughout the chain.  
The purpose of this section and of Section \ref{s:dense towers} is to
develop the machinery that will allow us to construct such a chain of
continuous towers while refining the index sets along the way.
While 
we will only use the fact that every continuous tower has a
continuous extension,
we prove the 
stronger statement 
to fuel the induction of Theorem \ref{exist cont ext}.

The claim that every continuous tower has a continuous extension still
alludes a full solution.  Hypothesis 1 is sufficient to
derive the extension property.  It is an open problem if this hypothesis
can be removed.

\begin{description}
\item[Hypothesis 1]  Every continuous tower of length $\alpha$ inside $\C$
has an amalgamable $<^c_{\mu,\alpha}$-extension inside $\C$.
\end{description}

\begin{theorem}[Existence of Continuous Extensions]
\label{exist cont ext}
Let $(\bar M,\bar a,\bar N)$ be a nice tower of length $\alpha$
in $\C$.  Under Hypothesis 1, there exists a continuous, amalgamable
tower $(\bar M^*,\bar a,\bar N)$ inside $\C$ such that
$(\bar M,\bar a,\bar N)<^c_{\mu,\alpha}(\bar M^*,\bar a,\bar N)$.  

Furthermore, 
if $(\bar M',\bar a,\bar N)\in\sq{+}\K^*_{\mu,\beta}$ is a
continuous partial extension of $(\bar M,\bar a,\bar N)$, then there
exist a $\prec_{\K}$-mapping $f$ and a continuous tower $(\bar M^*,\bar
a,\bar N)$ extending
$(\bar M,\bar a,\bar N)$ so that $f(M'_i)\preceq_{\K}M^*_i$ for all
$i<\beta$.

\end{theorem}

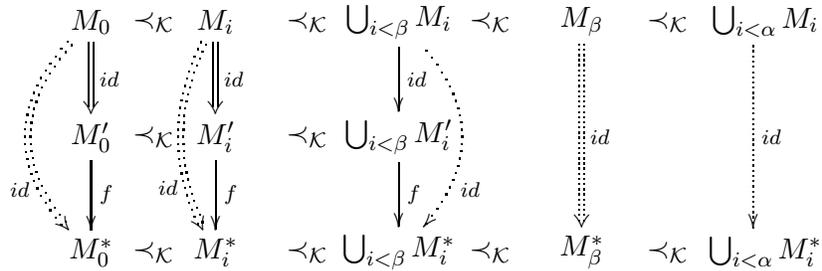
\[{\xymatrix 
{
M_0
\ar@2{->}[d]^{id}\ar@{}[r]|*+{\prec_{\K}}
\ar@/_2pc/@2{.>}[dd]_(.7){id}&
M_i\ar@2{->}[d]^{id}\ar@{}[r]|*+{\prec_{\K}}
\ar@/_1.1pc/@2{.>}[dd]_(.7){id} & \;\;\;\Union_{i<\beta}M_i\;\;\;
\ar[d]^{id}\ar@{}[r]|*+{\prec_{\K}} \ar@/^2pc/@{.>}[dd]^(.7){id}&
M_{\beta}\ar@2{.>}[dd]^{id}
\ar@{}[r]|*+{\prec_{\K}}&
\;\;\;\Union_{i<\alpha}M_i\ar@{.>}[dd]^{id}\\ 
M'_0\ar@{}[r]|*+{\prec_{\K}} \ar[d]^{f}& M'_i
\ar@{}[r]|*+{\prec_{\K}}\ar[d]^{f}& 
\Union_{i<\beta}M'_i \ar[d]^{f}& 
&  \\ M^*_0
\ar@{}[r]|*+{\prec_{\K}}& M^*_i\ar@{}[r]|*+{\prec_{\K}} &
\;\;\;\Union_{i<\beta}M^*_i\;\;\;
\ar@{}[r]|*+{\prec_{\K}}&
M^*_{\beta}\ar@{}[r]|*+{\prec_{\K}}
&\;\;\;\Union_{i<\alpha}M^*_i
}
}\]

The proof of Theorem \ref{exist cont ext} is by induction on $\alpha$. 
Notice that for $\alpha\leq\omega$, there is little to do since all
towers of length $\leq\omega$ are vacuously continuous.  If $\alpha$ is
the successor of a successor, then the induction hypothesis and the
extension property for non-$\mu$-splitting types (Theorem \ref{ext
property for non-splitting}) produce a continuous extension.   We take
care of the case that
$\alpha$ is a limit ordinal by taking direct limits of partial continuous
extensions. The difficult case is when $\alpha$ is the successor of a
limit ordinal.  This case employs Hypothesis 1.
  We will build an increasing
chain of continuous towers throwing in a particular element at each stage
so that in the end we will have added enough ($\mu$-many, predetermined)
elements to have a universal extension over $\Union_{i<\delta}M_i$. 
The following proposition allows us to add in the new elements in this
stage of the inductive proof of Theorem \ref{exist cont ext} (when
$\alpha=\delta+1$ and
$\delta$ is a limit ordinal).

\begin{proposition}\label{b in prop}
Suppose that Theorem \ref{exist cont ext} holds
for all amalgamable towers of length $\delta$ for some
limit ordinal $\delta<\mu^+$. Let
$(\bar M,\bar a,\bar N)$ be an amalgamable tower of length
$\delta$ inside $\C$.  For every
$b\in\C$, there exists a continuous, amalgamable tower
$(\bar M^*,\bar a,\bar N)\in\sq{+}\K^*_{\mu,\delta}$ inside $\C$
such that
$b\in\Union_{i<\delta}M^*_i$ 
and $(\bar M,\bar a,\bar N)<^c_{\mu,\delta}(\bar
M^*,\bar a,\bar N)$.

Furthermore,
if $(\bar M',\bar a,\bar N)\in\sq{+}\K^*_{\mu,\beta}$ is a
continuous partial extension of $(\bar M,\bar a,\bar N)$, we can
choose $(\bar M^*,\bar a,\bar N)$ such that  there exist a
$\prec_{\K}$-mapping
$f$ with $f(M'_i)\preceq_{\K}M^*_i$ for all $i<\beta$.

\end{proposition}

\begin{proof}
We begin by defining by induction on $\zeta<\delta$ 
a $<^c_{\mu,\delta}$-increasing
and continuous sequence of 
towers, $\langle (\bar M,\bar a,\bar
N)^\zeta\in\sq{+}\K^*_{\mu,\delta}\mid\zeta\leq\delta\rangle$
such that
\begin{enumerate}
\item $(\bar M,\bar a,\bar N)\leq^c_{\mu,\delta}(\bar M,\bar a,\bar
N)^0$,
\item $(\bar M,\bar a,\bar N)^\zeta$ is continuous and
\item if we are given $(\bar M',\bar a,\bar N)\in\sq{+}\K^*_{\mu,\beta}$
a continuous partial extension of $(\bar M,\bar a,\bar N)$, then there is
a $\prec_{\K}$-mapping $f$ with $f(M'_i)\preceq_{\K}M^0_i$ for all
$i<\beta$.
\end{enumerate}
This produces a $\delta$-by-$(\delta+1)$-array of models which we will 
diagonalize. 

Why is this construction possible?  
Since $(\bar M,\bar a,\bar N)$ is amalgamable, by the hypothesis of the
proposition, $(\bar M,\bar a,\bar N)\restriction\delta$ has a  continuous
extension $(\bar M^0,\bar a,\bar N)\in\sq{+}\K^{*}_{\mu,\delta}$. 
Furthermore, if we are given $(\bar M',\bar a,\bar
N)\in\sq{+}\K^*_{\mu,\beta}$ as above, then by condition $(2)$ of Theorem
\ref{exist cont ext}, we may find $f$ such that $f(M'_i)\preceq_{\K}M^0_i$
for all
$i<\beta$.
 At successor stages we can find continuous extensions
by the hypothesis of the proposition and the fact that continuous towers
are nice.  When
$\zeta$ is a limit ordinal, we take unions.  The unions will be
continuous, since the union of an increasing chain of continuous towers is
continuous. 

Since $\Union_{i<\delta}M_i$ is an amalgamation base, we can find an
isomorphic copy of this chain of towers inside $\C$.  WLOG, for
$\zeta<\delta$, $M_\zeta^\delta\prec_{\K}\C$.

Consider the diagonal sequence
$\langle M^\zeta_\zeta\mid\zeta<\delta\rangle$.  Notice that this is a
$\prec_{\K}$-increasing sequence of amalgamation bases. 
For $\zeta<\delta$, we have $M_{\zeta+1}^{\zeta+1}$ is universal over
$M_\zeta^\zeta$. Why?  From the
definition of $<^c$, 
$M^{\zeta+1}_\zeta$ is universal over $M^\zeta_\zeta$.  Since
$M^{\zeta+1}_\zeta\prec_{\K}M^{\zeta+1}_{\zeta+1}$, we have that
$M^{\zeta+1}_{\zeta+1}$ is also universal over $M^\zeta_\zeta$ (see
Proposition \ref{monotonicity of universality}).

By construction, each $\bar M^\zeta$ is continuous.  Thus
the sequence $\langle M^\zeta_\zeta\mid\zeta<\delta\rangle$ is continuous.
Then $\langle M^\zeta_\zeta\mid\zeta<\delta\rangle$ witnesses that
$\Union_{\zeta<\delta}M^\zeta_\zeta$ is a $(\mu,\delta)$-limit model.  Let
$M^b$ be a limit model inside $\C$ that is universal over 
$\Union_{\zeta<\delta}M^\zeta_\zeta$ and contains $b$.

Because $\Union_{\zeta<\delta}M^\zeta_\zeta$ is a limit
model, we can apply Fact
\ref{non-split thm} to\\
$\tp(b/\Union_{\zeta<\delta}M^\zeta_\zeta,M^\delta_\delta)$. 
Let
$\xi<\delta$ be such that 
$$(*)\quad\tp(b/\Union_{\zeta<\delta}M^\zeta_\zeta,M^b)
\text{ does not }\mu\text{-split over }M^\xi_\xi.$$

Notice that $(\langle M^i_i\mid i<\xi\rangle,\bar a,\bar
N)\restriction\xi$  is a $<^c_{\mu,\xi}$-extension of
$(\bar M,\bar a,\bar N)\restriction\xi$. 
 
We will find a $<^c_{\mu,\delta}$-extension of $(\bar M,\bar a,\bar
N)$ by defining an $\prec_{\K}$-increasing chain of models $\langle
N^*_i\mid i<\alpha\rangle$ and an increasing chain of
$\prec_{\K}$-mappings
$\langle h_i\mid i<\alpha\rangle$
 with the intention that the pre-image of $N^*_i$ under an extension of
$\Union_{i<\alpha}h_i$ will form a sequence $\bar M^*$ such that $(\bar
M,\bar a,\bar N)<^c_{\mu,\delta}(\bar M^*,\bar a,\bar N)$, $b\in
M^*_{\xi+1}$ and $M^*_i=M^i_i$ for all $i<\xi$.  
 We choose by induction on $i<\delta$
a $\prec_{\K}$-increasing and continuous chain of limit models $\langle
N^*_i\in\K_\mu\mid i<\delta\rangle$ and an
increasing and continuous sequence of
$\prec_{\K}$-mappings $\langle h_i\mid i<\delta\rangle$
satisfying

\begin{enumerate}
\item $N^*_{i+1}$ is a limit model and is universal over $N^*_i$
\item $h_i:M^i_i\rightarrow N^*_i$
\item $h_i(M^i_i)\prec_{\K}M^{i+1}_i$
\item $\tp(h_{i+1}(a_i/N^*_i)$ does not $\mu$-split over $h_i(N_i)$
\item $M^b\prec_{\K}N^*_{\xi+1}$ and 
\item for $i\leq\xi$, $N^*_i=M^i_i$ with $h_i=id_{M^i_i}$.
\end{enumerate}

We depict the construction below.  The inverse image of the
sequence of $N^*$'s will form the required continuous
$<^c_{\mu,\delta}$-extension of
$(\bar M,\bar a,\bar N)$.  

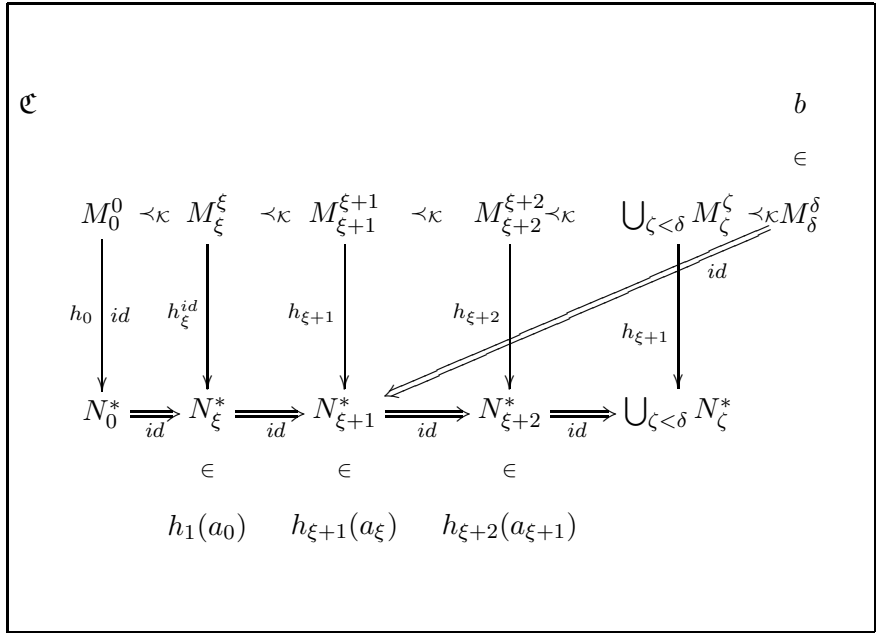
\[{\xymatrix@C10pt{
\\
\;\;\;\C&&&&&&b\ar@{}[d]|{\in}&\\
&M^0_0\ar[dd]_{h_{0}}^{id}\ar@{}[r]|{\prec_{\K}}
&M^{\xi}_{\xi} \ar[dd]_{h_{\xi}^{id}}\ar@{}[r]|{\prec_{\K}}&
M^{\xi+1}_{\xi+1}\ar[dd]_{h_{\xi+1}}\ar@{}[r]|{\prec_{\K}}&
M^{\xi+2}_{\xi+2}\ar[dd]_{h_{\xi+2}}
\ar@{}[r]|(0.3){\prec_{\K}}
 &\Union_{\zeta<\delta}M^\zeta_\zeta
\ar[dd]_(.6){h_{\xi+1}}\ar@{}[r]|(0.7){\prec_{\K}}&
M^\delta_\delta\ar@2{->}[ddlll]^(0.2){id}&&\\ &&&&&&\\
& N^*_0\ar@2{->}[r]_{id}
& N^*_{\xi} \ar@2{->}[r]_{id}&
N^*_{\xi+1}\ar@2{->}[r]_{id}&
N^*_{\xi+2}\ar@2{->}[r]_(0.4){id}&
\Union_{\zeta<\delta}N^*_{\zeta} &&\\
&&h_{1}(a_0\ar@{}[u]|{\in})&h_{\xi+1}(a_{\xi})\ar@{}[u]|{\in}&
h_{\xi+2}(a_{\xi+1})\ar@{}[u]|{\in}&
\\&&&&&&&
\save "1,1"."7,8"*[F]\frm{}
\restore}}
\] 

The requirements determine the definition of $N^*_i$ for $i\leq\xi$.  We
proceed with the rest of the construction by induction on $i$.  If
$i$ is a limit ordinal $\geq\xi$, let $N^*_i=\Union_{j<i}N^*_j$ and
$h_i=\Union_{j<i}h_j$.

Suppose that we have defined $h_i$ and $N^*_i$ satisfying the conditions
of the construction.  We now describe how to define $N^*_{i+1}$.
First, we extend $h_i$ to $\bar h_i\in\Aut(\C)$.  We can assume that
$\bar h_i(a_i)\in M^{i+2}_{i+1}$.  This is possible since $M^{i+2}_{i+1}$
is universal over $h_i(M^i_i)$ by construction.

 Since 
$\tp(a_i/M_i^i)$ does not $\mu$-split over $N_i$, by invariance we have
that
$\tp(\bar h_i(a_i)/h_i(M^i_i))$ does not $\mu$-split over $h_i(N_i)$.
We now adjust the proof of the existence property for non-splitting
extensions.

\begin{claim}\label{find g claim}
We can  find $g\in\Aut(\C)$ such that
$\tp(g(\bar h_i(a_i))/N^*_i)$ does not $\mu$-split over $h_i(N_i)$ and
$g(\bar h_i(M^{i+1}_{i+1}))\prec_{\K}M^{i+2}_{i+1}$.  
\end{claim}

\begin{proof}[Proof of Claim \ref{find g claim}]
First we find a $\prec_{\K}$-mapping $f$ such that $f:N^*_i\rightarrow
h_i(M^i_i)$ such that $f\restriction h_i(N_i)=id_{h_i(N_i)}$ which is
possible since $h_i(M^i_i)$ is universal over $h_i(N_i)$.  Notice that
$\tp(f^{-1}(\bar h_i(a_i))/N^*_i)$ does not $\mu$-split over $h_i(N_i)$
and
$$(+)\quad\tp(f^{-1}(\bar h_i(a_i))/h_i(M^i_i))=\tp(\bar
h_i(a_i)/h_i(M^i_i))$$ by a non-splitting argument as in the proof of
Theorem \ref{unique ext}.  

Let $N^+$ be
a limit  model of cardinality $\mu$ containing $f^{-1}(\bar h_i(a_i))$
with
$f^{-1}(\bar h_i(M^{i+1}_{i+1}))\prec_{\K}N^+$.
Now using the equality of types $(+)$ and the fact that $M^{i+2}_{i+1}$ is
universal over $h_i(M^i_i)$ with $\bar h_i(a_i)\in M^{i+2}_{i+1}$, we can
find a
$\prec_{\K}$-mapping
$f^+:N^+\rightarrow M^{i+2}_{i+1}$ such that
$f^+\restriction h_i(M^i_i)=id_{h_i(M^i_i)}$ and $f^+(f^{-1}(\bar
h_i(a_i)))=\bar h_i(a_i)$. 
Now set $g:=f^+\circ f^{-1}:\bar h(M^{i+1}_{i+1})\rightarrow
M^{i+2}_{i+1}$.

\end{proof}

Fix such a $g$ as in the claim and set $h_{i+1}:=g\circ \bar
h_i\restriction M^{i+1}_{i+1}$.  Let $N^*_{i+1}$ be a $\prec_{\K}$
extension of $N^*_i$, $M^b$ and $h_{i+1}(M^{i+1}_{i+1})$ of cardinality
$\mu$ inside $\C$.   Choose $N^*_{i+1}$ to additionally be a
limit model and universal over $N^*_i$.

 This completes the construction.

We now argue that the construction of these sequences is enough to find a
$<^c_{\mu,\delta}$-extension,
$(\bar M^*,\bar a,\bar N)$,
of
$(\bar M,\bar a,\bar N)$ such that
$b\in M^*_\zeta$ for some $\zeta<\delta$.  

Let $h_{\delta}:=\Union_{i<\delta}h_i$.
We will be defining for $i<\delta$,
$M^*_i$ to be pre-image of $N^*_i$ under some extension of $h_{\delta}$.  
The following claim allows us
to choose the pre-image so that $M^*_\zeta$ contains $b$ for some
$\zeta<\delta$.

\begin{claim}\label{exist the right aut claim}
There exists $h\in\Aut(\C)$ extending
$\Union_{i<\delta}h_i$ such that $h(b)=b$.
\end{claim} 
\begin{proof}[Proof of Claim \ref{exist the right aut claim}]
Let $h_\delta:=\Union_{i<\delta}h_i$.
Consider the increasing and continuous sequence
$\langle h_\delta(M^i_i)\mid i<\delta\rangle$.  By invariance,
$h_{\delta}(M^{i+1}_{i+1})$ is universal over
$h_{\delta}(M^i_i)$ and each $h_\delta(M^i_i)$ is a limit model.

Furthermore, from our  choice of $\xi$, we know   that $\tp(b/M^\delta_i)$
does not
$\mu$-split over $M^\xi_\xi$.  Since
$h_i(M^i_i)\prec_{\K}M^{i+1}_i\prec_{\K}\Union_{j<\delta}M^\delta_j$,
monotonicity of non-splitting allows us to conclude that 
$$\tp(b/h_\delta(M^i_i))\text{ does not }\mu\text{-split over
}M^\xi_\xi.$$

This allows us to apply 
Fact
\ref{non-split goes up},
to $\tp(b/\Union_{i<\delta}h_{\delta}(M^i_i))$ yielding

$$(**)\quad\tp(b/\Union_{i<\delta}h_{\delta}(M^i_i))\text{
does not }\mu\text{-split over }M^\xi_\xi.$$

Notice that $\Union_{i<\delta}M^i_i$ is a limit model
witnessed by $\langle M^j_j\mid j<\delta\rangle$.  So we can
apply Proposition 
\ref{mu,mu^+-limits are wmh} and extend
$\Union_{i<\delta}h_i$ to an automorphism $h^*$ of $\C$.  We will first
show that
$$(***)\quad\tp(b/h^*(\Union_{i<\delta}M^i_i),\C)
=\tp(h^*(b)/h^*(\Union_{i<\delta}M^i_i),\C).$$

By invariance and our choice of $\xi$ in $(*)$, 
$$\tp(h^*(b)/h^*(\Union_{i<\delta}M^i_i),\C)\text{ 
does not }\mu\text{-split over }M^\xi_\xi.$$

We will use non-splitting to derive $(***)$. 
To make the application of non-splitting more transparent, 
 let
$N^1:=\Union_{i<\delta}M^i_i$,
$N^2:=h^*(\Union_{i<\delta}M^i_i)$ and 
$p:=\tp(b/N^2)$.
By $(**)$, we have that $p\restriction N^2=h^{*}(p\restriction
N^1)$. In other words,
$$\tp(b/h^*(\Union_{i<\delta}M^i_i),\C)
=\tp(h^*(b)/h^*(\Union_{i<\delta}M^i_i),\C),$$ as desired.

From $(***)$ and Corollary \ref{type aut}, we can find an automorphism
$f$ of
$\C$ such that $f(h^*(b))=b$ and $f\restriction
h^*(\Union_{i<\delta}M^i_i)=
id_{h^*(\Union_{i<\delta}M^i_i)}$.  Notice that $h:=f\circ
h^*$ satisfies the conditions of the claim.

\end{proof}

Now that we have a automorphism $h$ fixing $b$ and
$\Union_{i<\delta}M_i$, we can define for each
$i<\delta$,
$M^*_i:=h^{-1}(N^*_i)$.  

\begin{claim}\label{its the right extension claim}
$(\bar M^*,\bar a,\bar
N)$ is a
$<^c_{\mu,\delta}$-extension of 
$(\bar M,\bar a,\bar N)$ such that
$b\in M^*_{\xi+1}$.
\end{claim}

\begin{proof}[Proof of Claim \ref{its the right extension claim}]
By construction $b\in M^\delta_\delta\subseteq N^*_{\xi+1}$.  Since
$h(b)=b$, this implies $b\in M^*_{\xi+1}$.  To
verify that we have a 
$\leq^c_{\mu,\delta}$-extension we need to show for $i<\delta$:
\begin{itemize}
\item[i.]  $M^*_i$ is universal over $M_i$ 
\item[ii.] $a_i\in M^*_{i+1}\backslash M_i$ for $i+1<\delta$ and
\item[iii.] $\tp(a_i/M^*_i)$ does not $\mu$-split over $N_i$ whenever
$i,i+1\leq\delta$.
\end{itemize} 
Item i. follows from the fact that $M^i_i$ is universal over $M_i$ and
$M^i_i\prec_{\K}M^*_i$.  Item iii. follows from invariance and
our construction of the $N^*_i$'s.  Finally, recalling that a
non-splitting extension of a non-algebraic type is also
non-algebraic (Remark
\ref{non-alg non-split}) we see that Item iii implies $a_i\notin M^*_i$.
By our choice of $h_{i+1}(a_i)\in
M^{i+2}_{i+1}\prec_{\K}N^*_{i+1}$, we have that $a_i\in M^*_{i+1}$.  Thus
Item ii is satisfied as well.

\end{proof}

\end{proof}

\medskip
Before beginning the proof of Theorem \ref{exist cont ext}, recall that
we will be building a directed system of partial extensions to take care
of the induction step when $\alpha$ is a limit ordinal.  Let us establish
a few facts about directed systems here.
Using the axioms of AEC and Shelah's Presentation Theorem, one can show
that Axiom \ref{union axiom} of the definition of AEC has an alternative
formulation (see
\cite{Sh 88} or Chapter 13 of
\cite{Gr2}): 

\begin{definition}\index{directed set}A partially ordered set $(I,\leq)$
is 
\emph{directed} iff for every $a,b\in I$, there exists
$c\in I$ such that $a\leq c$ and $b\leq c$.
\end{definition}

\begin{fact}[P.M. Cohn 1965]\label{direct limits}\index{direct
limit!existence of in AECs} Let $(I,\leq)$ be a directed set.  If $\langle
M_t\mid t\in I\rangle$ and $\{h_{t,r}\mid t\leq r\in I\}$ are such that
\begin{enumerate}
\item for $t\in I$, $M_t\in \K$
\item for $t\leq r\in I$, $h_{t,r}:M_t\rightarrow M_r$ is a
$\prec_{\K}$-embedding and
\item for $t_1\leq t_2\leq t_3\in I$, 
$h_{t_1,t_3}=h_{t_2,t_3}\circ h_{t_1,t_2}$ and $h_{t,t}=id_{M_t}$,
\end{enumerate}
then, whenever $s=\lim_{t\in I} t$, there exist
$M_s\in \K$ and $\prec_{\K}$-mappings
$\{h_{t,s}\mid t\in I\}$ such that
$$h_{t,s}:M_t\rightarrow M_s, M_s=\Union_{t<s}h_{t,s}(M_t)\text{ and}$$ 
$$\text{for }t_1\leq t_2\leq s, h_{t_1,s}=h_{t_2,s}\circ h_{t_1,t_2}
\text{ and }h_{s,s}=id_{M_s}.$$ 
\end{fact}
 
\begin{definition}
\begin{enumerate}

\item\index{directed system} $(\langle M_t\mid t\in I\rangle,\{h_{t,s}\mid
t\leq s\in I\})$ from Fact \ref{direct limits} is called a \emph{directed
system}.
\item\index{direct limit!definition}
We say that $M_s$ together with $\langle h_{t,s}\mid t\leq s\rangle$
satisfying the conclusion of Fact \ref{direct limits} is \emph{a
direct limit of $(\langle M_t\mid t<s\rangle,\{h_{t,r}\mid t\leq r<s\})$}.
\end{enumerate}
\end{definition}

Later we will generalize these systems by producing directed systems of
towers instead of models. 
%

Now we use Proposition \ref{b in prop} to prove Theorem \ref{exist cont
ext}.

\begin{proof}[Proof of Theorem \ref{exist cont ext}]
We prove that every amalgamable tower has a continuous extension by
induction on $\alpha$.  

\noindent $\pmb{\alpha=0:}$  By Theorem \ref{limits are ab} and Corollary
\ref{exist limit}, we can find a
$(\mu,\omega)$-limit over $M_0$.  Fix such a model and call it $M'_0$.

\noindent $\pmb{\alpha=\delta+1}$ {\bf and $\pmb{\delta}$ is a limit
ordinal:} 
The strategy is to start out with a continuous extension of $(\bar M,\bar
a,\bar N)\restriction\delta$ (which we call $(\bar M^{**},\bar
a\restriction\delta,\bar N\restriction\delta)$.)  If we are lucky, the top
of 
$(\bar M^{**},\bar
a\restriction\delta,\bar N\restriction\delta)$ will be universal over
$M_\delta$.  Since this cannot be guaranteed, we will repeatedly add new
elements into extensions of $(\bar M^{**},\bar a\restriction\delta,\bar
N\restriction\delta)$ until the top of one of these extensions is
universal over $M_\delta$.

By the induction hypothesis, we can find 
$(\bar M^{**},\bar a\restriction\delta,\bar
N\restriction\delta)\in\sq{+}\K^*_{\mu,\delta}$
 such that
\begin{itemize}
\item $(\bar M^{**},\bar a\restriction\delta,\bar N\restriction \delta)$
is a $<^c_{\mu,\delta}$-extension of $(\bar M,\bar a,\bar
N)\restriction\delta$ and
\item and if $(\bar M',\bar a\restriction\beta,\bar N\restriction\beta)$
is a continuous $<^c_{\mu,\beta}$-extension of $(\bar M,\bar a,\bar N)$,
then we can choose $\bar M^{**}$ such that there exists a
$\prec_{\K}$-mapping
$f$ with 
$f(M'_i)\prec_{\K}M^{**}_i$ for all $i<\beta$.
\end{itemize}
Notice that since $(\bar M^{**},\bar a\restriction\delta,\bar
N\restriction\delta)$ is continuous, we can apply
the induction hypothesis $\delta$-many times to find an
$<^c_{\mu,\delta}$-increasing chain of continuous towers of length
$\delta$.  In addition to being continuous, the top of this chain will be
an amalgamable extension of
$(\bar M,\bar a,\bar N)\restriction\delta$. Why?  The top of this tower
will be a $(\mu,\delta)$-limit model witnessed by the diagonal.  Thus
WLOG we may assume that $(\bar M^{**},\bar a\restriction\delta,\bar
N\restriction\delta)$ is amalgamable and continuous.

 We construct a continuous
$<^c_{\mu,\delta}$-extension of $(\bar M,\bar a,\bar N)$ by the
induction hypothesis and repeated applications of Proposition
\ref{b in prop}.  

Let $M'_\delta$ be a limit model and universal over $M_\delta$ inside
$\C$.   Enumerate
$M'_{\delta}$ as $\{b_\zeta\mid\zeta<\delta\mu\}$.  We will add these
elements into extensions of $(\bar M^{**},\bar a\restriction\delta,\bar
N\restriction\delta)$ by defining by induction on $\zeta\leq\delta\mu$ a
$<^c_{\mu,\delta}$-increasing and continuous chain of towers
$(\bar M,\bar a\restriction\delta,\bar
N\restriction\delta)^\zeta\in\sq{+}{\K^*_{\mu,\delta}}
$ such that
\begin{enumerate}
\item $(\bar M,\bar a\restriction\delta,\bar N\restriction\delta)^\zeta$
is a $<^c_{\mu,\delta}$-extension of $(\bar M^{**},\bar
a\restriction\delta,\bar N\restriction\delta)$
\item $(\bar M,\bar a\restriction\delta,\bar N\restriction\delta)^\zeta$
is continuous and 
\item $b_\zeta\in\Union_{i<\delta}M^{\zeta+1}_i\prec_{\K}\C$.
\end{enumerate}

The following diagram depicts the construction:

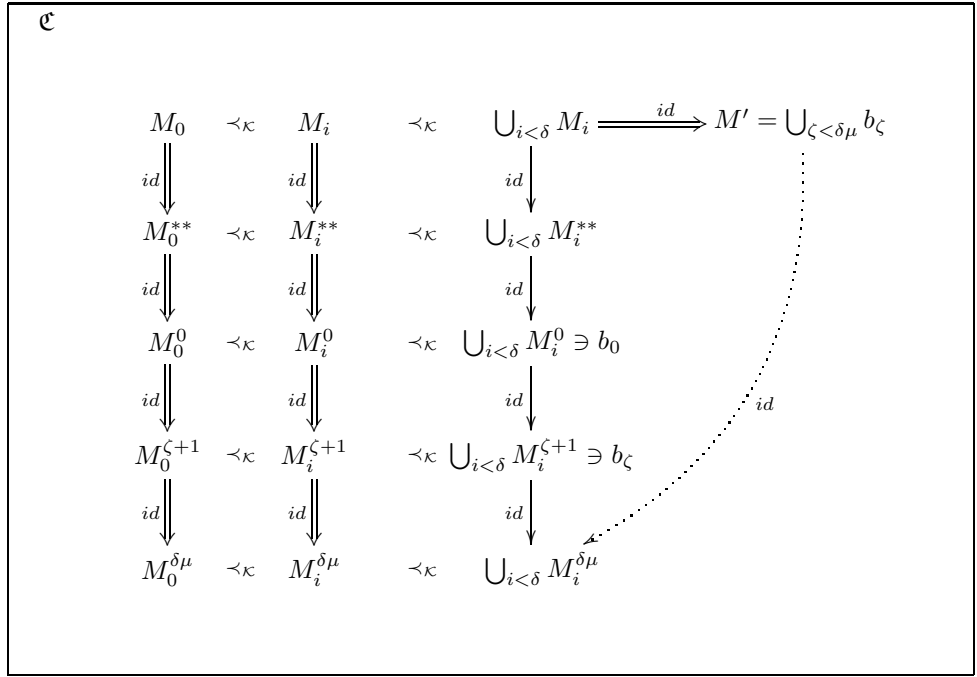
\[\small{\xymatrix{
\;\;\;\C&&&&&\\
&M_0\ar@2{->}[d]_{id}\ar@{}[r]|{\prec_{\K}}
&M_{i} \ar@2{->}[d]_{id}\ar@{}[r]|{\prec_{\K}}&
\;\;\;\Union_{i<\delta}M_{i}\ar[d]_{id}\ar@2{->}[r]^{id}&
M'=\Union_{\zeta<\delta\mu}b_\zeta\ar@/^3pc/@{.>}[ddddl]^{id}&\\
&M^{**}_0\ar@2{->}[d]_{id}\ar@{}[r]|{\prec_{\K}}
&M^{**}_{i} \ar@2{->}[d]_{id}\ar@{}[r]|{\prec_{\K}}&
\;\;\;\Union_{i<\delta}M^{**}_{i}\ar[d]_{id}&&
\\
&M^0_0\ar@2{->}[d]_{id}\ar@{}[r]|{\prec_{\K}}
&M^0_{i} \ar@2{->}[d]_{id}\ar@{}[r]|{\prec_{\K}}&
\;\;\;\Union_{i<\delta}M^0_{i}\ar[d]_{id}\ni b_0&&
\\
&M^{\zeta+1}_0\ar@2{->}[d]_{id}\ar@{}[r]|{\prec_{\K}}
&M^{\zeta+1}_{i} \ar@2{->}[d]_{id}\ar@{}[r]|{\prec_{\K}}&
\;\;\;\Union_{i<\delta}M^{\zeta+1}_{i}\ar[d]_{id}\ni b_\zeta&&
\\
&M^{\delta\mu}_0\ar@{}[r]|{\prec_{\K}}
&M^{\delta\mu}_{i} \ar@{}[r]|{\prec_{\K}}&
\;\;\;\Union_{i<\delta}M^{\delta\mu}_{i}&&
\\
&&&&&\\
\save "1,1"."7,6"*[F]\frm{}
\restore}}\]

The construction is possible by the induction hypothesis and Proposition
\ref{b in prop}:  

\noindent $\zeta=0$:  Since $\Union_{i<\delta}M^{**}_i$ is an amalgamation
base, we can apply Proposition \ref{b in prop} and find a
$<^c_{\mu,\delta}$-extension
$(\bar M^0,\bar a\restriction\delta,\bar N\restriction\delta)$ in $\C$
such that $b_0\in \Union_{i<\delta}M^0_i$.

\noindent $\zeta+1$:  Suppose that $(\bar M,\bar a\restriction\delta,\bar
N\restriction\delta)^\zeta$ has been defined.  It is a continuous tower
of length $\delta$.  If $\Union_{i<\delta}M^\zeta_i$ is an amalgamation
base, by the induction hypothesis we can apply Proposition
\ref{b in prop} to find a $<^c_{\mu,\delta}$-extension of $(\bar M,\bar
a\restriction\delta,\bar N\restriction\delta)^\zeta$, say $(\bar M,\bar
a\restriction\delta,\bar N\restriction\delta)^{\zeta+1}$ inside $\C$ such
that
$b_\zeta\in\Union_{i<\delta}M^{\zeta+1}_i$.  

Suppose on the other hand, that
$\Union_{i<\delta}M^\zeta_i$ is not an amalgamation base.  This may
occur when $\zeta$ is a limit ordinal of a different cofinality than
the cofinality of $\delta$.  By  Hypothesis 1, there is an
amalgamable extension of
$(\bar M,\bar a,\bar N)^\zeta$ inside $\C$.  Apply Proposition \ref{b in
prop} to the amalgamable extension and $b_\zeta$.  The proposition will
produce a 
$<^c_{\mu,\delta}$-extension of $(\bar M,\bar
a\restriction\delta,\bar N\restriction\delta)^\zeta$, say $(\bar M,\bar
a\restriction\delta,\bar N\restriction\delta)^{\zeta+1}$ inside $\C$ such
that
$b_\zeta\in\Union_{i<\delta}M^{\zeta+1}_i$.

\noindent \emph{$\zeta$ a limit ordinal}:  If $\zeta$ is a limit
ordinal we can set $(\bar M,\bar a\restriction\delta,\bar
N\restriction\delta)^\zeta:=\Union_{\xi<\zeta}(\bar M,\bar
a\restriction\delta,\bar N\restriction\delta)^\xi$.  It is a continuous
tower since all the $(\bar M,\bar a\restriction\delta,\bar
N\restriction\delta)^\xi$'s are continuous.  This completes the
construction.

 Now consider the tower $(\bar M^*,\bar a,\bar N)\in
\sq{+}{\K^*_{\mu,\delta+1}}$ defined by
$M^*_i:=M^{\delta\mu}_i$ for all $i<\delta$ and
$M^*_\delta:=\Union_{i<\delta}M^{\delta\mu}_i$.  Since $M^*_\delta$
contains
$M'_{\delta}$, it is universal over $M_{\delta}$.  Thus $(\bar M^*,\bar
a,\bar N)$ is a $<^c_{\mu,\delta+1}$-extension of $(\bar M,\bar a,\bar
N)$.  Since $(\bar M^{\delta\mu},\bar a\restriction\delta,\bar
N\restriction\delta)$ is continuous, we have that $(\bar M^*,\bar a,\bar
N)$ is also continuous.  Notice that $(\bar M^*,\bar a,\bar N)$ is
amalgamable as well.  By construction for every $i<\delta$, $M^*_\delta$
is a limit model.  For the case $i=\delta$, we see that
$M^*_\delta$ is a
$(\mu,\delta)$-limit model witnessed by the diagonal 
$\langle M^{i\mu}_{i}\mid i<\delta\rangle$.

\noindent$\pmb{\alpha=\delta+1}$ {\bf and $\pmb{\delta}$ is a successor
ordinal:}  
By the induction hypothesis we can find a continuous, amalgamable
extension
$(\bar M^{**},\bar a\restriction\delta,\bar N\restriction\delta)$ of
$(\bar M,\bar a,\bar N)\restriction\delta$ and if we are given $(\bar
M',\bar a\restriction\beta,\bar N\restriction\beta)$ as in part $(2)$ of
the statement of the theorem, we may assume that there is a
$\prec_{\K}$-mapping
$f^*$ such that $f^*(M'_i)\prec_{\K}M^{**}_i$ for all $i<\beta$.

Since $M^{**}_{\delta-1}$ and $M_\delta$ are both
$\K$-substructures of
$\C$, we can apply the Downward-L\"{o}wenheim Axiom for AECs to find
$M^{**}_{\delta}$ (a first approximation to $M^{*}_{\delta}$) a model
of cardinality
$\mu$ extending both $M^{**}_{\delta-1}$ and $M_\delta$.  WLOG by
Proposition
\ref{mu,mu+ limit is univ} and Lemma \ref{univ ext containing a} we may
assume that
$M^{**}_{\delta}$ is a limit model of cardinality $\mu$ and
$M^{**}_{\delta}$ is universal over both $M^{**}_{\delta-1}$ and
$M_\delta$. By Theorem \ref{ext property for non-splitting},
we can find a $\prec_{\K}$-mapping
$h:M^{**}_{\delta}\rightarrow\C$ such that
$h\restriction M_\delta=id_{M_\delta}$ and
$\tp(a_{\delta}/h(M^{**}_{\delta}))$ does not $\mu$-split over
$N_{\delta}$. Set $M^*_i:=h(M^{**}_{i})$ for all $i\leq\delta$.
Notice that by invariance $(\bar M^*,\bar a,\bar N)\restriction\delta$
is a $<^c_{\mu,\delta}$-extension of $(\bar M,\bar a,\bar N)$.  To
conclude that $(\bar M^*,\bar a,\bar
N)$ is the required
$<^c_{\mu,\alpha}$-extension of $(\bar M,\bar a,\bar N)$ with
$f=h\circ f^*$ if appropriate, it remains to check that   
\begin{subclaim}\label{a not in subclaim cont proof}
$a_{\delta}\notin M^*_\delta$
\end{subclaim}
\begin{proof}[Proof of Subclaim \ref{a not in subclaim cont proof}]
Suppose that $a_{\delta}\in M^*_\delta$.
Since $M_{\delta}$ is universal over $N_{\delta}$, there
exists a $\prec_{\K}$-mapping, $g:M^*_\delta\rightarrow M_{\delta}$
such that
$g\restriction N_{\delta}=id_{N_{\delta}}$.
Since $\tp(a_{\delta}/M^*_\delta)$ does not $\mu$-split
over $N_{\delta}$, we have that
$$(*)\quad\tp(a_{\delta}/g(M^*_\delta))=\tp(g(a_{\delta})/g(M^*_\delta)).$$
Notice that because $g(a_{\delta})\in g(M^*_\delta)$, $(*)$ implies
that $a_{\delta}=g(a_{\delta})$.  Thus
$a_{\delta}\in g(M^*_\delta)\prec_{\K}M_{\delta}$. This contradicts the
definition of towers: $a_{\delta}\notin M_{\delta}$.

\end{proof}

\noindent$\pmb{\alpha}$ {\bf is a limit ordinal }$\pmb{>\omega:}$
We will construct a 
directed system of partial extensions of $(\bar M,\bar a,\bar N)$,
$\langle (\bar M,\bar a,\bar N)^\zeta\mid\zeta<\alpha\rangle$
and $\langle f_{\xi,\zeta}\mid\xi\leq\zeta<\alpha\rangle$
satisfying the following conditions:

\begin{enumerate}
\item $(\bar M,\bar a,\bar N)\restriction\zeta<^c_{\mu,\zeta}(\bar
M,\bar a,\bar N)^\zeta$

\item\label{nice cond} $(\bar M,\bar a,\bar N)^\zeta$ is
continuous
\item $(\bar M,\bar a,\bar N)^\zeta$ lies in $\C$ 
\item $f_{\xi,\zeta}\restriction M_{i}^\xi:M_{i}^\xi\rightarrow
M_{i}^\zeta$ for $i<\xi\leq\zeta$ 
\item\label{univ over cond of cont ext constr} for all $\xi<\zeta$,
$M^\zeta_\xi
$ is universal
over $f_{\xi,\zeta}(\Union_{i<\xi}M^\xi_i)$ and 
\item $f_{\xi,\zeta}\restriction M_\xi=id_{M_\xi}$ for all
$\xi<\zeta<\alpha$.
\end{enumerate}

The construction is possible by the induction hypothesis and Proposition
\ref{b in prop}.  We provide the details here.  

\noindent $\zeta=0$: Set $\bar M^0$ equal to the empty sequence and
$f_{0,0}$ equal to the empty mapping.  

\noindent $\zeta=\xi+1$:  Suppose that $(\bar M,\bar a,\bar N)^\xi$ and
$\langle f_{\gamma,\gamma'}\mid \gamma\leq\gamma'\leq\xi\rangle $ have
been defined accordingly.  Then by the induction hypothesis applied to
$(\bar M,\bar a,\bar N)\restriction\zeta$ and the partial extension
$(\bar M,\bar a,\bar N)^\xi$, we can find a $\prec_{\K}$-mapping $f$ and a
continuous extension of $(\bar M,\bar a,\bar N)\restriction\zeta$.
By applying the induction hypothesis again to this continuous extension,
we can find $(\bar M,\bar a,\bar N)^\zeta\in\sq{+}\K^*_{\mu,\zeta}$
inside $\C$ such that for all $i<\xi$, $f(M^\xi_i)\prec_{\K}M^\zeta_i$,
$f\restriction M_i=id_{M_i}$  and $M^\zeta$ is universal over
$f(\Union_{i<\xi}M^\xi_i)$.  Notice that by setting
$f_{\gamma,\xi+1}=f\circ f_{\gamma,\xi}$ and
$f_{\zeta,\zeta}=id_{\Union_{\xi<\zeta}M^\zeta_\xi}$ we have completed the
successor stage of the construction.

\noindent $\zeta$ a limit ordinal:  
By the induction hypothesis we have constructed a directed system $\langle
\Union_{i<\gamma}M^\gamma_{i}\mid \gamma<\zeta\rangle$ with $\langle
f_{\gamma,\xi}\mid \gamma\leq\xi<\zeta\rangle$.  By Fact \ref{direct
limits} we can find a direct limit to this system, $M^{**}_\zeta\in\K$ and
$\prec_{\K}$-mappings
$\langle f^{**}_{\gamma,\zeta}\mid\gamma\leq\zeta\rangle$.
First notice that
\begin{subclaim}\label{fs incr in cont thm}
$\langle f^{**}_{\gamma,\zeta}\restriction
M_\gamma\mid\gamma<\zeta\rangle$ is increasing.
\end{subclaim}

\begin{proof}
Let $\gamma<\xi<\zeta$ be given.  By construction
$$f_{\gamma,\xi}\restriction M_\gamma=id_{M_\gamma}.$$
An application of $f^{**}_{\xi,\zeta}$ yields
$$f^{**}_{\xi,\zeta}\circ f_{\gamma,\xi}\restriction
M_\gamma=f^{**}_{\gamma,\zeta}\restriction M_\gamma.$$ 
Since $f^{**}_{\gamma,\zeta}$ and $f^{**}_{\xi,\zeta}$ come from a direct
limit of the system which includes the mapping $f_{\gamma,\xi}$,
we have
$$
f^{**}_{\gamma,\zeta}\restriction M_\gamma=
f^{**}_{\xi,\zeta}\circ f_{\gamma,\zeta}\restriction
M_\gamma.$$
Combining the equalities yields
$$  f^{**}_{\gamma,\zeta}\restriction
M_\gamma=f^{**}_{\xi,\zeta}\restriction M_\gamma.$$ 
This completes the proof of Subclaim \ref{fs incr in cont thm}.
\end{proof}

By the subclaim, we have that
$f:=\Union_{\gamma<\zeta}f^{**}_{\gamma,\zeta}\restriction M_\gamma$ is a
$\prec_{\K}$-mapping from
$\Union_{\gamma<\zeta}M_\gamma$ onto
$\Union_{\gamma<\zeta}f^{**}_{\gamma,\zeta}(M_\gamma)$.  Since $\C$ is a
$(\mu,\mu^+)$-limit model 
and since $\Union_{\gamma<\zeta}M_\gamma$ is an amalgamation base (as
$(\bar M,\bar a,\bar N)$ is nice) we can assume that 
$f$ is a partial automorphism of $\C$ and extend it
to
$F\in Aut(\C)$ by Proposition \ref{mu,mu^+-limits are
wmh}.  

Now
consider the direct limit defined by
$M^\zeta_\zeta:=F^{-1}(M^{**}_\zeta)$ with
$\langle f^*_{\xi,\zeta}:=F^{-1}\circ f^{**}_{\xi,\zeta}\mid \xi<
\zeta\rangle$ and $f^*_{\zeta,\zeta}=id_{M^*_\zeta}$. 
Let $M^\zeta_i:=f_{\xi,\zeta}(M^\xi_i)$ for all $i<\xi$.  This is
well-defined since $f_{\xi,\zeta}$ is part of the direct limit of a
directed system.
 Notice that
$f^*_{\xi,\zeta}\restriction M_\xi=F^{-1}\circ
f^{**}_{\xi,\zeta}\restriction M_\xi=id_{M_\xi}$ for
$\xi<\zeta$.

\begin{subclaim}\label{why a tower subclaim}
 $(\bar M,\bar a,\bar
N)\restriction\zeta<^c_{\mu,\zeta}(\bar M,\bar a,\bar N)^\zeta$.
\end{subclaim}
\begin{proof}[Proof of Subclaim \ref{why a tower subclaim}]
We need to verify that for all
$\xi<\zeta$, 
\begin{enumerate}
\item $M^\zeta_\xi\prec_{\K}M^\zeta_{\xi+1}$, 
\item $a_\xi\in
M^\zeta_{\xi+1}\backslash M^\zeta_{\xi}$ and
\item $\tp(a_\xi/M^\zeta_\xi)$ does not $\mu$-split over $N_\xi$.
\end{enumerate}
To see that $\bar M^\zeta$ is increasing, by the induction hypothesis,
$$f_{\xi,\xi+1}(\Union_{i<\xi}M^\xi_i)\prec_{\K}M^{\xi+1}_\xi.$$
 Applying $f_{\xi+1,\zeta}$ to both sides of this equation gives us for
every $j<\xi$, 
$$M^\zeta_j\prec_{\K}f_{\xi,\zeta}(\Union_{i<\xi}M^\xi_i)=f_{\xi+1,\zeta}(f_{\xi,\xi+1}
(\Union_{i<\xi}M^\xi_i))\prec_{\K}f_{\xi+1,\zeta}(M^{\xi+1}_\xi)
=M^\zeta_\xi.$$

By the induction hypothesis for all
$\xi<\zeta$, $a_\xi\notin M^{\xi+2}_\xi$ and $\tp(a_\xi/M^{\xi+2}_\xi)$
does not $\mu$-split over $N_\xi$.  Since $f_{\xi+2,\zeta}\restriction
M_{\xi+1}=id_{M_{\xi+1}}$, invariance gives us
$f_{\xi+2,\zeta}(a_\xi)=a_\xi\notin
f_{\xi+2,\xi}(M^{\xi+2}_\xi)=M^\zeta_\xi$ and
$\tp(a_\xi/M^\zeta_\xi)$ does not $\mu$-split over $N_\xi$.

\end{proof}

 Notice that $(\bar
M,\bar a,\bar N)^\zeta$ is continuous since it is formed from the direct
limit of a continuous system.  To see that $(\bar M,\bar a,\bar N)^\zeta$
is amalgamable, notice that condition (\ref{univ over cond of cont ext
constr}) of the construction guarantees that
$\Union_{\xi<\zeta}M^\zeta_\xi$ is a
$(\mu,\zeta)$-limit witnessed by $\langle
f_{\xi,\zeta}(\Union_{i<\xi}M_i^\xi)\mid \xi<\zeta\rangle$.
This completes the construction.

Why is the construction sufficient to produce $(\bar M',\bar a,\bar N)$ as
required?
We have constructed a directed
system
$\langle \Union_{i<\gamma}M^\gamma_{i}\mid
\gamma\leq\xi<\alpha\rangle$ with
$\langle f_{\gamma,\xi}\mid \gamma\leq\xi<\alpha\rangle$.  By Fact
\ref{direct limits} and Subclaim \ref{fs incr in cont thm} we can find a
direct limit to this system,
$M^{*}_\alpha$ and
$\prec_{\K}$-mappings
$\langle f_{\gamma,\alpha}\mid\gamma\leq\alpha\rangle$ such that
$f_{\gamma,\alpha}\restriction M_i=id_{M_i}$ for all
$i<\alpha$.  
If $(\bar M,\bar a,\bar N)$ is amalgamable, then $M^{*}_\alpha$ can be
chosen to lie in $\C$.
Define for all $\zeta<\alpha$,
$M^*_\zeta:=f_{\zeta+1,\alpha}(M^{\zeta+1}_\zeta)$.  Notice that as in
Subclaim \ref{why a tower subclaim}, $(\bar M,\bar a,\bar
N)<^c_{\mu,\alpha}(\bar M^*,\bar a,\bar N)$.  And, as in the limit stage
of the construction, we see that $(\bar M^*,\bar a,\bar N)$ is continuous
and amalgamable.

The second part of the statement of the theorem
 is obtained by modifying our construction by setting $(\bar
M,\bar a,\bar N)^\beta=(\bar M',\bar a,\bar N)$ and proceeding with the
construction from $\beta+1$.

\end{proof}

\bigskip

\section{Refined Orderings on Towers}\label{s:dense towers}

In this section we further develop the machinery of towers which will be
used to construct a relatively full tower in Section
\ref{s:unique limits}.

\begin{definition}For ordinals $\alpha,\alpha',\delta,\delta'<\mu^+$ with
$\alpha\leq\alpha'$ and $\delta\leq\delta'$.  We say that $(\bar M',\bar
a',\bar N')\in\sq{+}\K^*_{\mu,\alpha'\times\delta'}$ is a $<^c$-extension
of
$(\bar M,\bar a,\bar N)\in\sq{+}\K^*_{\mu,\alpha\times\delta}$ iff 
\begin{itemize}
\item for every
$\beta<\alpha$ and every $i<\delta$, $M'_{\beta,i}$ is universal over
$M_{\beta,i}$
\item for every $\beta<\alpha$ and $i+1<\delta$,
$a_{\beta,i}=a'_{\beta,i}$ and $N_{\beta,i}=N'_{\beta,i}$.
\end{itemize}
\end{definition}

The following theorem is used to construct relatively full towers by
adding realizations of strong types between $M_{\beta,i}$ and
$M_{\beta+1,0}$ in an $<^c$-extension of the tower $(\bar M,\bar a,\bar
N)\in\sq{+}\K^*_{\mu,\alpha\times\delta}$.

\begin{theorem}\label{adding the as}
Under Hypothesis 1, given $\alpha$ an ordinal $<\mu^+$ and  a nice tower,
$(\bar M,\bar a,\bar N)\in\sq{+}\K^*_{\mu,\alpha\times\mu\alpha}$, we can
find an amalgamable, continuous extension
$(\bar M',\bar a',\bar N')\in\sq{+}\K^*_{\mu,\alpha+1\times\mu(\alpha+1)}$
of $(\bar M,\bar a,\bar N)$ such that for  a fixed
enumeration,
$\{(p,N)^{\zeta}_l\mid
l<\mu\}$,  of
$\Union_{i<\mu\alpha}\St(M_{\zeta,i})$ for each $\zeta<\alpha$, we have
that
$$(*)\quad(p,N)_l^{\zeta}\sim(\tp(a_{\zeta+1,l+1}/
M'_{\zeta+1,l+1}),N_{\zeta+1,l+1})
\restriction \dom(p_l^{\zeta}).$$ 
\end{theorem}

\begin{proof}
We begin by constructing $(\bar M',\bar a,\bar N)$, a continuous,
amalgamable
$<^c_{\mu,\alpha\times\mu\alpha}$-extension  of
$(\bar M,\bar a,\bar N)$, such that for $\zeta+1<\alpha$, $M'_{\zeta+1,0}$
is a $(\mu,\mu)$-limit over $\Union_{i<\mu\alpha}M'_{\zeta,i}$.  The
construction of $(\bar M',\bar a,\bar N)$ is done by defining a directed
system of amalgamable, continuous partial extensions of $(\bar M,\bar
a,\bar N)$ using Theorem \ref{exist cont ext}. 
Specifically, Theorem \ref{exist cont ext} allows us to define by
induction on
$\zeta$, a directed system
$\langle (\bar M,\bar a,\bar N)^\zeta\mid1\leq\zeta\leq\alpha\rangle$
and $\langle f_{\xi,\zeta}\mid1\leq\xi\leq\zeta\leq\alpha\rangle$
satisfying the following conditions:

\begin{enumerate}
\item $(\bar M,\bar a,\bar
N)\restriction(\zeta\times\mu\alpha)<^c_{\mu,\zeta\times\mu\alpha}(\bar
M,\bar a,\bar N)^\zeta$

\item\label{nice cond} $(\bar M,\bar a,\bar N)^\zeta$ is
continuous and amalgamable
\item $(\bar M,\bar a,\bar N)^\zeta$ lies in $\C$ for $\zeta<\alpha$
\item $M^{\zeta+1}_{\zeta+1,0}$ is a $(\mu,\mu)$-limit over
$\Union_{i<\mu\alpha}M^\zeta_{\zeta,i}$
\item\label{univ over cond of new element const} for all $\xi<\zeta$,
$M^\zeta_\xi
$ is universal
over $f_{\xi,\zeta}(\Union_{i<\xi}M^\xi_i)$ 
\item $f_{\xi,\zeta}\restriction M_{i}^\xi:M_{i}^\xi\rightarrow
M_{i}^\zeta$ for $i<\xi\leq\zeta$ and 
\item $f_{\xi,\zeta}\restriction M_\xi=id_{M_\xi}$ for all
$\xi<\zeta<\alpha$.
\end{enumerate}

The details of the direct limit construction are similar to the direct
limit construction in the limit case of Theorem \ref{exist cont ext}.

The construction is sufficient:
 Let $(\bar M',\bar a,\bar N):=(\bar
M,\bar a,\bar N)^\alpha$.  For each
$\zeta+1<\alpha$, fix a sequence $\langle M^*_{\zeta,i}\mid i< \mu\rangle
$ witnessing that $M'_{\zeta+1,0}$ is a $(\mu,\mu)$-limit over
$\Union_{i<\mu\alpha}M'_{\zeta,i}$.  Define
$M'_{\zeta,\mu\alpha+i}:=M^*_{\zeta,i}$ for each $i<\mu$ and
$\zeta+1<\alpha$.

\bigskip
\noindent$\small{\xymatrix@C=15pt@M=1pt{
M_{0,0}\ar@2{->}[d]^{id}\ar@{}[r]|{\prec_{\K}} & M_{0,i}
\ar@2{->}[d]^{id}\ar@{}[r]|{\prec_{\K}}&
\;\;\Union_{i<\mu\alpha}M_{0,i}
\ar@2{->}[d]^{id}\ar@{}[rrrr]|{\prec_{\K}}
&&&&M_{1,0}\ar@2{->}[d]^{id}\\
M'_{0,0}\ar@{}[r]|{\prec_{\K}} & M'_{0,i}
\ar@{}[r]|{\prec_{\K}}&
\;\;\Union_{i<\mu\alpha}M'_{0,i}\ar@{=}[r]&
M^*_{0,0}\ar@2{.>}[r]^(.5){id}&
M^*_{0,j}\ar@2{.>}[r]^(.4){id}
&M^*_{0,j+1}\ar@2{.>}[r]^(.3){id}
&\Union_{j<\mu}M^*_{0,j}=M'_{1,0}\\
}}$

\bigskip

For each $\zeta+1<\alpha$ and each $l<\mu$, by the Theorem \ref{ext
property for non-splitting}, we can find $q\in
\gaS(M'_{\zeta+1,\mu\alpha+l})$  extending $p^{\zeta}_l$ such that $q$
does not
$\mu$-split over $N^\zeta_l$.  Since $M'_{\zeta+1,\mu\alpha+l+1}$ is
universal over $M'_{\zeta+1,\mu\alpha+l}$, there is $a\in
M'_{\zeta+1,\mu\alpha+l+1}$ realizing $q$.  Set
$a_{\zeta+1,\mu\alpha+l}=a$ and $N_{\zeta+1,\mu\alpha+l}=N^{\zeta}_l$.
This gives us a definition of $(\bar M',\bar a,\bar
N)\in\sq{+}\K^*_{\mu,\alpha\times\mu(\alpha+1)}$.  To extend this tower to
a tower with index set $(\alpha+1)\times\mu(\alpha+1)$, we use the fact
that
$(\bar M',\bar a,\bar N)$ is amalgamable to fix
$M^*$ a $(\mu,\mu(\alpha+1))$-limit model over
$\Union_{i<\mu\alpha,\zeta<\alpha}M'_{\alpha,i}$.  Let $\langle
M'_{\alpha,i}\mid i<\mu(\alpha+1)\rangle$ witness this.  WLOG  we
may assume that $M'_{\alpha,i+1}$ is a $(\mu,\omega)$-limit over
$M'_{\alpha,i}$ for each $i<\mu(\alpha+1)$. For each
$i<\mu(\alpha+1)$, fix
$a_{\alpha,i}\in M'_{\alpha,i+1}\backslash M'_{\alpha,i}$.  By Fact
\ref{non-split thm} and our choice of $M'_{\alpha,i}$ as a limit model,
there is a
$N\prec_{\K}M'_{\alpha,i}$ such that $M'_{\alpha,i}$ is universal over
$N$ and $\tp(a_{\alpha,i}/M'_{\alpha,i})$ does not $\mu$-split over
$N$.  Set $N_{\alpha, i}=N$.  Notice that $(\bar M',\bar a,\bar
N)\in\sq{+}\K^*_{\mu,(\alpha+1)\times\mu(\alpha+1)}$ is as required.
\end{proof}

\bigskip

\section{Uniqueness of Limit Models}\label{s:unique limits}

Recall the running assumptions:
\begin{enumerate}
\item $\K$ is an abstract elementary class,
\item $\K$ has no maximal models, 
\item $\K$ is categorical in some $\lambda>LS(\K)$,
\item GCH and $\Phi_{\mu^+}(S^{\mu^+}_{\cf(\mu)})$ holds
for every cardinal $\mu<\lambda$.
\end{enumerate}

Under these assumptions and Hypothesis 1, we can prove the uniqueness of
limit models using the results from Sections 
\ref{s:full} and \ref{s:dense towers}.

\begin{theorem}[Uniqueness of Limit Models]\label{uniqueness of limit
models}\index{uniqueness of limit models}\index{uniqueness of limit
models!conjecture!solution} Let $\mu$ be a cardinal
$\theta_1,\theta_2$ limit ordinals such that
$\theta_1,\theta_2<\mu^+\leq\lambda$.
Under Hypothesis 1, if $M_1$ and $M_2$ are $(\mu,\theta_1)$ and
$(\mu,\theta_2)$ limit models over $M$,
respectively, then there exists an isomorphism $f:M_1\cong M_2$ such that
$f\restriction M=id_M$.
\end{theorem}

\begin{proof}Let $M\in\K^{am}_\mu$ be given.
By Fact \ref{unique limits}, it is enough to show that there exists a
$\theta_2$ such that for every $\theta_1$ a limit ordinal $<\mu^+$, we
have that a $(\mu,\theta_1)$-limit model over $M$ is isomorphic to
a $(\mu,\theta_2)$-limit model over $M$.
Take $\theta_2$ such that $\theta_2=\mu\theta_2$.  Fix $\theta_1$ a limit
ordinal
$<\mu^+$. By Fact \ref{sigma and cf(sigma) limits}, we may assume that
$\theta_1$ is regular. Using Fact \ref{unique limits} again, it is enough
to construct a model
$M^*$ which is simultaneously a $(\mu,\theta_1)$-limit model over $M$ and
a
$(\mu,\theta_2)$-limit model over $M$.

 The idea is to build a (scattered) array of models such that
at some point in the array, we will find a model which is a
$(\mu,\theta_1)$-limit model witnessed by its height in the array and is
a $(\mu,\theta_2)$-limit model witnessed by its horizontal position in
the array, relative fullness and continuity. 
  We will
define a chain of length
$\mu^+$ of continuous towers while increasing the index
set of the towers in order to realize strong types as we proceed with the
goal of producing many relatively full rows.

Define by induction on $0<\alpha<\mu^+$ the
$<^c$-increasing 
sequence of  towers, $\langle(\bar M,\bar a,\bar
N)^{\alpha}\in\sq{+}\K^*_{\mu,\alpha\times\mu\alpha}\mid
\alpha<\mu^+\rangle$, such that

\begin{enumerate}
\item $M\prec_{\K}M^{\alpha}_{0,0}$,
\item $(\bar M,\bar a,\bar N)^{\alpha}$ is 
continuous and amalgamable,
\item $(\bar M,\bar a,\bar
N)^{\alpha}:=\Union_{\beta<\alpha}(\bar M,\bar
a,\bar N)^{\beta}$ for $\alpha$ a limit ordinal and
\item\label{putting in reps} In successor stages in new intervals of
length
$\mu$, put in representatives of all $\St$-types from the previous stages.
More formally, 
if $(p,N)\in\St(M^\alpha_{\beta,i})$ for $i<\mu\alpha$ and
$\beta<\alpha$, there exists
$j\in[\mu\alpha,\mu(\alpha+1)]$ such that
$$(p,N)\sim(\tp(a_{\beta+1,j}/M_{\beta+1,j}^{\alpha+1}),N_j)\restriction
M^\alpha_{\beta,i}.$$

\end{enumerate}
This construction is possible:\\
$\alpha=1$:  We can choose $\bar M^*=\langle M^{*}_i\mid i<\mu\rangle$
to be an  $\prec_{\K}$ increasing continuous sequence of limit models of
cardinality $\mu$ with $M^{*}_0=M$ and $M^*_{i+1}$ universal over
$M^*_i$.  For each
$i<\mu$, fix
$a^{1}_{0,i}\in M^{*}_{i+1}\backslash M^{*}_{i}$.  Now
consider
$\tp(a^{1}_{0,i}/M^{*}_{i})$.  
Since $M^{*}_{i}$ is a
limit model, we can apply Fact \ref{non-split thm} to fix
$N^{1}_{0,i}\in\K^{am}_\mu$ such that 
$\tp(a^{1}_{0,i}/M^{*}_{i})$ does not $\mu$-split over 
$N^{1}_{0,i}$ and $M^{*}_{i}$ is universal over
$N^{1}_{0,i}$.  Let $\bar a^1:=\langle a^1_{0,i}\mid i<\mu\rangle$
and $\bar N^1=\langle N^1_{0,i}\mid i<\mu\rangle$.
\\
$\alpha$ a limit ordinal:  Take $(\bar M,\bar a,\bar
N)^{\alpha}:=\Union_{\beta<\alpha}(\bar M,\bar
a,\bar N)^{\beta}$.  Clearly $(\bar M,\bar a,\bar N)^\alpha$ is
continuous.  To see that $(\bar M,\bar a,\bar N)^\alpha$ is also
amalgamable, we notice that
$\Union_{\beta,i\in\alpha\times\mu\alpha}M^\alpha_{(\beta,i)}$
is a $(\mu,\alpha)$-limit model witnessed by 
$\langle \Union_{i<\mu\beta}M^\beta_{\beta,i}\mid\beta<\alpha\rangle$.
\\
$\alpha=\beta+1$:  Suppose that $(\bar M,\bar a,\bar
N)^{\beta}$ has been defined.  
 By Fact
\ref{St small}, for every $\gamma<\beta$, we can enumerate
$\Union_{k<\mu\beta}\St(M^\beta_{\gamma,k})$
as $\{(p,N)^{\gamma}_l\mid l<\mu\}$.
By
Theorem \ref{adding the as}, we can find a continuous, amalgamable
extension
$(\bar M,\bar a,\bar
N)^{\beta+1}\in\sq{+}\K^*_{\mu,\beta+1\times\mu(\beta+1)}$
of $(\bar M,\bar a,\bar N)^\beta$
such that for every $l<\mu$ and $\gamma<\beta$,
$$(p,N)^{\gamma}_l\sim(\tp(a_{\gamma+1,\mu\beta+l}/
M^{\beta+1}_{\gamma+1,\mu\beta+l}),N_{\gamma+1,\mu\beta+l})
\restriction \dom(p^\gamma_l).$$ 
  This completes the
construction.

We now want to identify all the rows of the array which are relatively
full.

\begin{claim}\label{identify rel full claim}
For $\delta$ a limit ordinal
$<\mu^+$, we have that $(\bar M,\bar a,\bar N)^{\delta}$
is full relative to 
$\langle \bar M^\delta_{\beta,i}\mid
(\beta,i)\in\delta\times\mu\delta\rangle$ where
$\bar M^\delta_{\beta,i}:=\langle M^\gamma_{\beta,i}\mid
\gamma<\delta\rangle$.
\end{claim}
\begin{proof}
Let $(p,N)\in\St(M^{\delta}_{\beta,i})$ be given such that
$N=M^\gamma_{\beta,i}$ for some $\gamma<\delta$, $\beta<\delta$ and
$i<\mu\delta$.  Since our construction is increasing and continuous, there
exists
$\delta'<\delta$ such that
$(\beta,i)\in\delta'\times\mu\delta'$ and $\gamma<\delta'$.  
Notice then that $M^{\delta'}_{\beta,i}$ is universal over $N$. 
Furthermore, $p\restriction M^{\delta'}_{\beta,i}$ does not $\mu$-split
over $N$.  Thus $(p,N)\restriction
M^{\delta'}_{\beta,i}\in\St(M^{\delta'}_{\beta,i})$. 
By condition (\ref{putting in reps}) of the construction,
there exists $j<\mu(\delta'+1)$, such that
$$(p,N)\restriction
M^{\delta'}_{\beta,i}\sim(\tp(a_{\beta+1,j}/M^{\beta+1}_{\beta+1,j}),
N_{\beta+1,j})\restriction
M^{\delta'}_{\beta,i}.$$

Since $M^{\beta+1}_{\beta+1,j}\prec_{\K}M^{\delta}_{\beta+1,j}$ and
$\tp(a_{\beta+1,j}/M^\delta_{\beta+1,j})$ does not $\mu$-split over
$N_{\beta+1,j}$, we can replace $M^{\beta+1}_{\beta+1,j}$ with
$M^{\delta}_{\beta+1,j}$:
$$(p,N)\restriction
M^{\delta'}_{\beta,i}\sim(\tp(a_{\beta+1,j}/M^{\delta}_{\beta+1,j}),
N_{\beta+1,j})\restriction
M^{\delta'}_{\beta,i}.$$
Let $M'$ be a universal extension of $M^{\delta}_{\beta+1,j}$.
By definition of $\sim$, there exists $q\in\gaS(M')$
such that
$q$ extends $p\restriction
M^{\delta'}_{\beta,i}=\tp(a_{\beta+1,j}/M^{\delta'}_{\beta,i})$ and $q$
does not $\mu$-split over $N$ and $N_{\beta+1,j}$.  By the uniqueness of
non-splitting extensions (Theorem
\ref{unique ext}), since $p$ does not $\mu$ split over $N$, we have that
$q\restriction M^\delta_{\beta,i}=p$.  Also, since
$\tp(a_{\beta+1,j}/M^{\delta}_{\beta+1,j})$ does not $\mu$-split over
$N_{\beta+1,j}$, Theorem \ref{unique ext} gives us
$q\restriction
M^\delta_{\beta+1,j}=\tp(a_{\beta+1,j}/M^{\delta}_{\beta+1,j})$. By
definition of
$\sim$ and Lemma \ref{enough to consider M'}, $q$ also witnesses that
$$(\tp(a_{\beta+1,j}/M^{\delta}_{\beta+1,j}),
N_{\beta+1,j})\restriction M^\delta_{\beta,i}\sim(p,N).$$  Since $(p,N)$
was chosen arbitrarily, we have verified that $(\bar M,\bar a,\bar
N)^{\delta}$ satisfies the definition of relative fullness.
\end{proof}

Take $\langle\delta_\zeta<\mu^+\mid\zeta\leq\theta_1\rangle$ to be an
increasing and continuous sequence of limit ordinals $>\theta_2$.  
We will consider the restrictions (in the sense of Notation \ref{rest
notation special}) of
$(\bar M,\bar a,\bar N)^{\delta_\zeta}$ to
$\theta_2\times\mu\delta_\zeta$:

\begin{notation}\label{rest notation special}
For $\theta$ and $\delta$ ordinals $<\mu^+$ and a sequence $\bar M$
indexed by a superset of $\theta\times\mu\delta$, we will abbreviate
$\langle M_{\beta,i}\mid
\beta<\theta\text{ and }i<\mu\delta\}$ by $\bar M\restriction^{
\theta\times\mu\delta}$.
\end{notation}

Define
$$M^*:=\Union_{\zeta<\theta_1}
\Union_{i\in\theta_2\times\mu\delta_\zeta}M^{\delta_{\zeta}}_i=
\Union_{i\in\theta_2\times\mu\delta_{\theta_1}}M^{\delta_{\theta_1}}_i.$$


We will now verify that $M^*$ is a $(\mu,\theta_1)$-limit over $M$ 
and a $(\mu,\theta_2)$-limit over $M$.

Notice that $\langle
\Union_{i\in\theta_2\times\mu\delta_\zeta}M^{\delta_\zeta}_i
\mid\zeta<\theta_1\rangle$ witnesses that $M^*$ is a $(\mu,\theta_1)$
limit.
Since $M\prec_{\K}M^{\delta_0}_{0,0}$, $M^*$ is a $(\mu,\theta_1)$-limit
over $M$.

By Claim \ref{identify rel full claim} and the fact that the restriction
of a relatively full tower is relatively full (Proposition
\ref{restriction of full}), we have that
$$(\bar M,\bar a,\bar N)^{\delta_{\theta_1}}\restriction
^{\theta_2\times\mu\delta_\zeta}\text{ is
full relative to }
\langle \bar
M^{\delta_{\theta_1}}_{\beta,i}\mid
(\beta,i)\in\theta_2\times\mu\delta_{\theta_1}\rangle,$$
where $\bar M^{\delta_{\theta_1}}_{\beta,i}:=\langle
M^\gamma_{\beta,i}\mid
\gamma<\delta_{\theta_1}\rangle$.
  Furthermore, we see that $(\bar M,\bar a,\bar
N)^{\delta_{\theta_1}}\restriction^{\theta_2\times\mu\delta_{\theta_1}}$
is continuous.  Since $\theta_2=\mu\cdot\theta_2$,
we can apply Theorem
\ref{full is limit} to conclude that $M^*$ is a
$(\mu,\theta_2)$-limit model over $M$.
\end{proof}

\begin{remark}\label{relatively full
exist} The above proof implicitly shows the decomposition of a relatively full
tower into a resolution of $\theta'$ many towers for every limit
$\theta'<\mu^+$.
\end{remark}

\bigskip

\part{Conclusion}

We provide a partial proof of Hypothesis 1.  We also
discuss reduced towers, which appear in the \cite{ShVi 635} and may be
useful as a tool to prove the amalgamation property for categorical AECs
with no maximal models.  We will continue to make Assumption \ref{assm
intro}.

\bigskip

\section{$<^c_{\mu,\alpha}$-Extension Property for
Nice Towers}\label{s:<b extension property}

%
%
%

In \cite{ShVi 635}, Shelah and Villaveces claim that every tower in
$\sq{+}{\K}^*_{\mu,\alpha}$ has a proper $<^c_{\mu,\alpha}$ extension. 
This proof does not converge.  
Here we prove a weaker extension
property.  Namely, we show that every \emph{nice} tower in
$\sq{+}{\K^*_{\mu,\alpha}}$  has a proper
$<^c_{\mu,\alpha}$-extension (Corollary \ref{c-ext of nice}).  
This is a proof of an approximation to the statement of Hypothesis 1 which
states that every continuous tower has an amalgamable extension inside
$\C$.

\begin{theorem}\label{partial ext}
Let $\mu$ be a cardinal  and
$\alpha,\gamma$ ordinals such that
$\gamma<\alpha<\mu^+\leq\lambda$.
If
$(\bar M,\bar a,\bar N)\in\sq{+}\K^{*}_{\mu,\alpha}$ is nice
and $(\bar M{''},\bar a,\bar N)\restriction\gamma$ is an amalgamable
partial extension of
$(\bar M,\bar a,\bar N)$, then there exists an amalgamable
$(\bar M^*,\bar a,\bar N)\in\sq{+}\K^{*}_{\mu,\alpha}$
and a $\prec_{\K}$-mapping $f$
 such that
\begin{enumerate}
\item $(\bar M,\bar a,\bar N)<^c_{\mu,\alpha}(\bar M^*,\bar a,\bar N)$
\item $f(M''_i)= M^*_i$ for all
$i<\gamma$ and
\item $f\restriction M_i=id_{M_i}$ for all $i<\gamma$.
\end{enumerate}

Furthermore if $\Union_{i<\alpha}M_i\prec_{\K}\C$ and
$\bar b\in\sq{\leq\mu}\C$ is such that $\bar
b\cap\Union_{i<\alpha}M_i=\emptyset$, then we can find $(\bar M^*,\bar
a,\bar N)$ as above with $\bar b\cap \Union_{i<\alpha}M^*_i=\emptyset$.

\end{theorem}

\begin{remark}
If $(\bar M,\bar a,\bar N)$ is amalgamable and
$\Union_{i<\alpha}M_i\prec_{\K}\C$, then we can find an extension $(\bar
M',\bar a,\bar N)$ such that
$\Union_{i<\alpha}M'_i\prec_{\K}\C$.

\end{remark}

Theorem \ref{partial ext} is stronger than the
$<^c_{\mu,\alpha}$-extension property since it allows us to avoid
$\mu$-many elements $(\bar b)$.  This is possible due to 
Weak Disjoint Amalgamation, Fact \ref{wda}. 

%
%
%

\begin{proof}[Proof of Theorem \ref{partial ext}]
Let an amalgamable
$(\bar M,\bar a,\bar N)\in\sq{+}\K^{*}_{\mu,\alpha}$ be given.  

As in the proofs of Theorem \ref{exist cont ext} and \ref{adding the as},
we will define by induction on $i<\alpha$ 
a direct system of models 
$\langle M'_i\mid i<\alpha\rangle$ and 
$\prec_{\K}$-mappings, $\langle f_{j,i}\mid j<i<\alpha\rangle$ such that
for $i\leq\alpha$:
\begin{enumerate}
\item\label{cond c-ext} $(\langle f_{j,i}(M'_j)\mid j\leq i\rangle,\bar
a\restriction i+1,\bar N\restriction i+1)$ is a $<^c_{\mu,i+1}$-extension
of
$(\bar M,\bar a,\bar N)\restriction (i+1)$,
\item $(\langle M'_j\mid j<i\rangle,\langle f_{j,i}\mid j\leq i\rangle)$
forms a directed system,
\item
$M'_{i}$ is universal over $M_i$,
\item\label{M' univ over previous image} $M'_{i+1}$ is universal over
$f_{i,i+1}(M'_i)$,
\item\label{f's are id on M}$f_{j,i}\restriction M_j=id_{M_j}$,
\end{enumerate}

Notice that the $M'_i$'s will not necessarily form an
extension of the tower $(\bar M,\bar a, \bar N)$.  Rather, for
each $i<\alpha$, we find some image of $\langle M'_j\mid j<i\rangle$
which will extend the initial segment of length $i$ of
$(\bar M,\bar a,\bar N)$ (see condition (\ref{cond c-ext}) of the
construction). 

The construction is possible:

\emph{$i=0$}:  Since $M_0$ is an amalgamation base, we can find
$M''_0\in\K^*_\mu$ (a first approximation of the desired $M'_0$) such
that
$M''_0$ is universal over $M_0$.   
By Theorem \ref{ext property for non-splitting},
we may assume that $\tp(a_0/M''_0)$ does not 
$\mu$-split over $N_0$ and $M''_0\prec_{\K}\C$.
Since $a_0\notin M_0$ and
$\tp(a_0/M_0)$ does not $\mu$-split over $N_0$,
we know that $a_0\notin M''_0$.  But, we might have that for
some $l>0$, $a_l\in M''_0$ or $\bar b\cap M''_0\neq\emptyset$.  We use
Weak Disjoint Amalgamation to avoid $\{a_l\mid 0<l<\alpha\}$ and $\bar
b$.  By the Downward L\"{o}wenheim-Skolem Axiom for AECs (Axiom
\ref{DLS}) we can choose 
$M^2\in\K_\mu$ such that $M_0''$, $M_1\prec_{\K}M^2\prec_{\K}\C$.

By Corollary \ref{wda improv1} (applied to $M_1$, $M_\alpha$,
$M^2$ and $\langle a_l\mid 0<l<\alpha\rangle\cup\bar b)$,
we can find a $\prec_{\K}$-mapping $h$ such that
\begin{itemize}
\item $h:M^2\rightarrow \C$
\item $h\restriction M_1=id_{M_1}$
\item $h(M^2)\cap(\{a_l\mid 0<l<\alpha\}\cup\bar b)=\emptyset$
\end{itemize} 

Define $M'_0:=h(M''_0)$.  Notice that $a_0\notin M'_0$
because $a_0\notin M''_0$ and $h(a_0)=a_0$.  Clearly 
$M'_0\cap(\{a_l\mid 0\leq l<\alpha\}\cup\bar b)=\emptyset$, since
$M''_0\prec_{\K}M^2$ and $h(M^2)\cap\{a_l\mid 0<l<\alpha\}=\emptyset$.
We need only verify that
$\tp(a_0/M'_0)$ does not $\mu$-split over $N_0$.
By invariance,
$\tp(a_0/M''_0)$ does not $\mu$-split over $N_0$ implies
that 
$\tp(h(a_0)/h(M''_0))$ does not $\mu$-split over $N_0$.
But recall $h(a_0)=a_0$ and $h(M''_0)=M'_0$.
Thus $\tp(a_0/M'_0)$ does not $\mu$-split over $N_0$.

Set 
$f_{0,0}:=id_{M'_0}$.

Below is a diagram of the successor stage of the
construction.

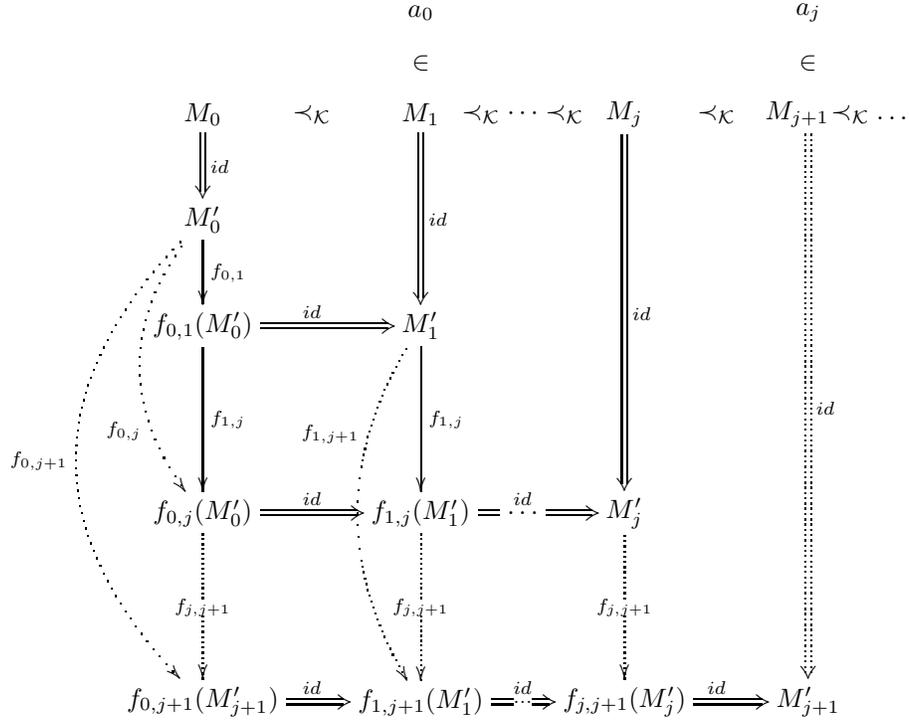
\[\small{{\xymatrix 
{&&&&\\
 & a_0\ar@{}[d]|*+{\in}&&a_j\ar@{}[d]|*+{\in}&\\
M_0\ar@2{->}[d]^{id}
\ar@{}[r]|*+{\prec_{\K}}
&M_1\ar@2{->}[dd]^{id}\ar@{}[r]|*+{\prec_{\K}\dots\prec_{\K}}&
M_j\ar@{}[r]|*+{\prec_{\K}}\ar@2{->}[dddd]^{id}&
M_{j+1}\;\;\;\ar@{}[r]|*+{\prec_{\K}\dots}
\ar@2{.>}[dddddd]^{id}&{}&\\
M'_{0} \ar[d]^{f_{0,1}}\ar@/_2pc/@{.>}[ddd]_(.7){f_{0,j}}
\ar@/_4pc/@{.>}[ddddd]_{f_{0,j+1}}&&&&&\\
f_{0,1}(M'_{0})\ar[dd]^{f_{1,j}}\ar@2{->}[r]^{id}
&M'_1\ar[dd]^{f_{1,j}}\ar@/_2pc/@{.>}[dddd]_(.3){f_{1,j+1}} &&&&\\
 & &&&&\\
f_{0,j}(M'_{0})\ar@{.>}[dd]|{f_{j,j+1}}\ar@2{->}[r]^{id}
&f_{1,j}(M'_1)\ar@{.>}[dd]|{f_{j,j+1}}
\ar@2{->}[r]|*+{\dots}^{id}
&M'_{j}\ar@{.>}[dd]|{f_{j,j+1}}&&&\\
 & &&&&\\
f_{0,j+1}(M'_{j+1})\ar@2{->}[r]^{id}
&f_{1,j+1}(M'_1)\ar@2{->}[r]|{\dots}^{id}
&f_{j,j+1}(M'_j)\ar@2{->}[r]^{id}
&M'_{j+1}&&
}}}\]

\emph{$i=j+1$}:
Suppose that we have completed the construction for all $k\leq j$.
Since $M'_j$ and $M_{j+1}$ are both $\K$-substructures of $\C$, we
can apply the Downward-L\"{o}wenheim Axiom for AECs to find
$M'''_{j+1}$ (a first approximation to $M'_{j+1}$) a model of cardinality
$\mu$ extending both $M'_j$ and $M_{j+1}$.  WLOG by Proposition
\ref{mu,mu+ limit is univ} and Lemma \ref{univ ext containing a} we may
assume that
$M'''_{j+1}$ is a limit model of cardinality $\mu$ and
$M'''_{j+1}$ is universal over $M_{j+1}$ and $M'_j$.
By Theorem \ref{ext property for non-splitting},
we can find a $\prec_{\K}$-mapping
$f:M'''_{j+1}\rightarrow\C$ such that
$f\restriction M_{j+1}=id_{M_{j+1}}$ and
$\tp(a_{j+1}/f(M'''_{j+1}))$ does not $\mu$-split over $N_{j+1}$.
Set $M''_{j+1}:=f(M'''_{j+1})$.  
\begin{subclaim}\label{a not in subclaim}
$a_{j+1}\notin M''_{j+1}$
\end{subclaim}
\begin{proof}[Proof of Subclaim \ref{a not in subclaim}]
Suppose that $a_{j+1}\in M''_{j+1}$.
Since $M_{j+1}$ is universal over $N_{j+1}$, there
exists a $\prec_{\K}$-mapping, $g:M''_{j+1}\rightarrow M_{j+1}$
such that
$g\restriction N_{j+1}=id_{N_{j+1}}$.
Since $\tp(a_{j+1}/M''_{j+1})$ does not $\mu$-split
over $N_{j+1}$, we have that
$$\tp(a_{j+1}/g(M''_{j+1})=\tp(g(a_{j+1})/g(M''_{j+1})).$$
Notice that because $g(a_{j+1})\in g(M''_{j+1})$, we have
that $a_{j+1}=g(a_{j+1})$.  Thus
$a_{j+1}\in g(M''_{j+1})\prec_{\K}M_{j+1}$. This contradicts the
definition of towers: $a_{j+1}\notin M_{j+1}$.

\end{proof}

$M''_{j+1}$ may serve
us well if it does not contain any $a_l$ for $j+1\leq l<\alpha$ or any
part of $\bar b$, but this is not guaranteed.  So we need to make an
adjustment. Let $M^2$ be a model of cardinality $\mu$ such that
$M_{j+2}, M''_{j+1}\prec_{\K}M^2\prec_{\K}\C$.
Notice that $\C$ is universal over $M_{j+2}$.
Thus we can apply Corollary \ref{wda improv1} to
$M_{j+2}$, $M_{\alpha}$, $M^2$ and $\langle a_l\mid 
j+2\leq l<\alpha\rangle\cup\bar b$.  This yields a 
$\prec_{\K}$-mapping $h$ such that
\begin{itemize}
\item $h:M^2\rightarrow \C$
\item $h\restriction M_{j+2}=id_{M_{j+2}}$ and
\item $h(M^2)\cap(\{a_l\mid j+2\leq l<\alpha\}\cup\bar b)=\emptyset$.
\end{itemize}

Set $M'_{j+1}:=h(M''_{j+1})$.  Notice that by invariance,
$\tp(a_{j+1}/M''_{j+1})$ does not $\mu$-split over $N_{j+1}$ implies
that $\tp(h(a_{j+1})/h(M''_{j+1}))$ does not $\mu$-split over
$h(N_{j+1})$.  Recalling that 
$h\restriction M_{j+2}=id_{M_{j+2}}$ we have
that
$\tp(a_{j+1}/M''_{j+1})$ does not $\mu$-split over $N_{j+1}$.
We need to verify that $a_{j+1}\notin M'_{j+1}$.
This holds because $a_{j+1}\notin M''_{j+1}$ and
$h(a_{j+1})=a_{j+1}$.

Set $f_{j+1,j+1}=id_{M_{j+1}}$  and
$f_{j,j+1}:=h\circ f\restriction M'_j$.
To guarantee that we have a directed system, for $k<j$, define
$f_{k,j+1}:=f_{j,j+1}\circ f_{k,j}$.

\emph{$i$ is a limit ordinal}:
Suppose that 
$(\langle M'_j\mid j<i\rangle, \langle f_{k,j}\mid k\leq j<i\rangle)$
 have been defined.
Since it is a directed system, we can take
direct limits. 
\begin{subclaim}\label{dir limit construction}
We can choose a direct limit $(M^*_i,\langle f^*_{j,i}\mid
j\leq i\rangle)$ of $(\langle M'_j\mid j<i\rangle, \langle f_{k,j}\mid
k\leq j<i\rangle)$  such that

\begin{enumerate}
\item $M^*_i\prec_{\K}\C$
\item $f^*_{j,i}\restriction M_j=id_{M_j}$ for every $j<i$.
\end{enumerate}
\end{subclaim}

\begin{proof}[Proof of Subclaim \ref{dir limit construction}]
This follows from Subclaim \ref{fs incr in cont thm} and the assumption
that $(\bar M,\bar a,\bar N)$ is nice.
\end{proof}

By Condition (\ref{M' univ over previous image}) of the construction,
notice that $M^*_i$ is a $(\mu,i)$-limit model 
witnessed by $\langle f^*_{j,i}(M'_j)\mid j<i\rangle$.
Hence $M^*_i$ is an amalgamation base.
Since $M^*_i$ and $M_i$ both live inside of $\C$,
we can find $M'''_i\in\K^*_\mu$ which is universal over $M_i$ 
and universal over $M^*_i$.  

By Theorem \ref{ext property for non-splitting}
we can find a $\prec_{\K}$-mapping
$f:M'''_{i}\rightarrow\C$ such that
$f\restriction M_i=id_{M_i}$ and
$\tp(a_i/f(M'''_i))$ does not $\mu$-split
over $N_i$.  Set $M''_i:=f(M'''_i)$.
By a similar argument to Subclaim \ref{a not in subclaim},
we can see that $a_i\notin M''_i$.

$M''_i$ may contain some $a_l$ when $i\leq l<\alpha$ or part of $\bar b$.
We need to make an adjustment using Weak Disjoint Amalgamation.
Let $M^2$ be a model of cardinality $\mu$ such that
$M''_i,M_{i+1}\prec_{\K}M^2\prec_{\K}\C$.
By
Corollary \ref{wda improv1} applied to
$M_i$, $M_\alpha$, $M^2$ and $\langle a_l\mid i<l<\alpha\rangle\cup\bar b$
we can find $h:M''_i\rightarrow \C$
such that $h\restriction M_{i+1}=id_{M_{i+1}}$ and
$h(M^2)\cap(\{a_l\mid i< l<\alpha\}\cup\bar b)=\emptyset$.

Set $M'_i:=h(M''_i)$.  We need to verify that
$a_i\notin M'_i$ and
$\tp(a_i/M'_i)$ does not $\mu$-split over $N_i$.
Since $a_i\notin M''_i$ and $h(a_i)=a_i$,
we have that $a_i\notin h(M''_i)=M'_i$.
By invariance of non-splitting, 
$\tp(a_i/M''_i)$ not $\mu$-splitting over $N_i$
implies that
$\tp(h(a_i)/h(M''_i))$ does not $\mu$-split over $h(N_i)$.
Recalling our definition of $h$ and $M'_i$, this yields
$\tp(a_i/M'_i)$ does not $\mu$-split over $N_i$.

As in the proof of Theorem \ref{exist cont ext}, we see that $(\langle
(f_{j,i}(M_j)\mid j\leq i\rangle, \bar a\restriction i,\bar
N\restriction i)$ is a $<^c_{\mu,i}$-extension of $(\bar M,\bar a,\bar
N)\restriction i$.

Set $f_{i,i}:=id_{M_i,i}$,
and
for $j<i$, $f_{j,i}:=h\circ f\circ f^*_{j,i}$.
  This completes the construction.

The construction is enough: 
We have constructed a directed
system
$\langle \Union_{i<\gamma}M'_{i}\mid
i<\alpha\rangle$ with
$\langle f_{i,j}\mid i\leq j<\alpha\rangle$.  By Fact
\ref{direct limits} and Subclaim \ref{fs incr in cont thm} we can find a
direct limit to this system,
$M^*_\alpha$ and
$\prec_{\K}$-mappings
$\langle f_{i,\alpha}\mid i\leq\alpha\rangle$ such that
$f_{i,\alpha}\restriction M_i=id_{M_i}$ for all
$i<\alpha$ and $M^*_\alpha$ avoids $\bar b$.
If $(\bar M,\bar a,\bar N)$ is amalgamable, then $M^*_\alpha$ can be
chosen to lie in $\C$.
Define for all $j<\alpha$,
$M^*_j:=f_{j+1,\alpha}(M'_j)$.  Notice that as in
Subclaim \ref{why a tower subclaim}, $(\bar M,\bar a,\bar
N)<^c_{\mu,\alpha}(\bar M^*,\bar a,\bar N)$.  And, as in the limit stage
of the construction, we see that $(\bar M^*,\bar a,\bar N)$ is continuous
and amalgamable.

\end{proof}

\begin{remark}
Notice that in Theorem \ref{partial ext} if the partial extension $(\bar
M',\bar a,\bar N)$ is continuous, then we can choose $\bar M''$ such that
it is continuous below $\gamma$, that is for every $i<\gamma$ with $i$ a
limit ordinal,
$M''_i=\Union_{j<i}M''_j$.
\end{remark}

\begin{corollary}[The $<^c_{\mu,\alpha}$-extension property for
nice towers]\label{c-ext of
nice}\index{$<^c_{\mu,\alpha}$-extension property|(} If\\
$(\bar M,\bar a,\bar N)\in\sq{+}\K^{*}_{\mu,\alpha}$ is nice, then
there exists an amalgamable
$(\bar M',\bar a,\bar N)\in\sq{+}\K^{*}_{\mu,\alpha}$ such that
$(\bar M,\bar a,\bar N)<^c_{\mu,\alpha}(\bar M',\bar a,\bar N)$.

\end{corollary}
\begin{proof}
Take $\gamma=0$ in Theorem \ref{partial ext}
\end{proof}

\begin{remark}
Notice that Hypothesis 3 implies that every tower is amalgamable.  Thus
Hypothesis 3 together with Corollary \ref{c-ext of nice} imply the 
$<^c_{\mu,\alpha}$-extension property for all towers.  
\end{remark}

\bigskip

\section{Reduced Towers}\label{s:reduced new section}

Shelah and Villaveces introduce the notion of reduced towers in order to
show the density of continuous towers.  While
there are difficulties with Shelah and Villaveces' approach, we discuss
reduced towers because they have characteristics similar to strongly
minimal types in first-order model theory.  
Additionally, they generalize reduced triples used in
\cite{Sh 576} to develop a notion of non-forking.

\begin{definition}\label{reduced defn}\index{reduced towers}
A tower $(\bar M,\bar a,\bar N)\in\sq{+}\K^*_{\mu,\alpha}$ is said to 
be \emph{reduced} provided that for every $(\bar M',\bar a,\bar
N)\in\sq{+}\K^*_{\mu,\alpha}$ with
$(\bar M,\bar a,\bar N)\leq^c_{\mu,\alpha}(\bar M',\bar a,\bar
N)$ we have that for every
$i<\alpha$,
$$(*)_i\quad M'_i\cap\Union_{j<\alpha}M_j = M_i.$$
\end{definition}

If we take a $<^c$-increasing
chain of reduced towers, the union
will be reduced.  The following proposition appears in \cite{ShVi 635}
(Theorem 3.1.14 of \cite{ShVi 635}) for reduced towers.  We provide the
proof for completeness.

\begin{theorem}\label{union of reduced is reduced}
If
$\langle (\bar M,\bar a,\bar N)^\gamma\in\sq{+}\K^*_{\mu,\alpha}\mid
\gamma<\beta\rangle$ is $<^c_{\mu,\alpha}$-increasing and continuous
sequence of
reduced towers, then the union of this sequence of towers is a
reduced tower.
\end{theorem}
\begin{proof}
Denote by $(\bar M,\bar a,\bar N)^\beta$ the union of
the sequence of towers.  That is $\bar a^\beta=\bar a^0$, $\bar
N^\beta=\bar N^0$ and $\bar M^\beta=\langle M^\beta_i\mid i<\alpha\rangle$
where $M^\beta_i=\Union_{\gamma<\beta}M^\gamma_i$.

  Suppose that $(\bar M,\bar a,\bar N)^\beta$ is not reduced. 
Let
$(\bar M',\bar a,\bar N)\in\sq{+}\K^*_{\mu,\alpha}$ witness this.
Then there exists an $i<\alpha$ and an element $b$ such that
$b\in (M'_i\cap\Union_{j<\alpha}M^\beta_j)\backslash M^\beta_i$.
There exists $\gamma<\beta$ such that
$b\in \Union_{j<\alpha}M^\gamma_j\backslash M^\gamma_i$.  
Notice that $(\bar M',\bar a,\bar N)$ witnesses that
$(\bar M,\bar a,\bar N)^\gamma$ is not reduced.
\end{proof}

The following 
appears in \cite{ShVi 635} (Theorem 3.1.13).

\begin{fact}[Density of reduced towers]\label{density of
reduced}\index{reduced towers!density of} There exists a reduced
$<^c_{\mu,\alpha}$-\\extension of every nice tower in
$\sq{+}\K^*_{\mu,\alpha}$.

\end{fact}
\begin{proof}

Suppose for the
sake of contradiction that no
$<^c_{\mu,\alpha}$-extension of the tower
$(\bar M,\bar a,\bar N)$  is 
reduced.  This allows us to construct a
$\leq^c_{\mu,\alpha}$-increasing and continuous sequence of towers 
$\langle(\bar M,\bar a,\bar
N)^\zeta\in\sq{+}\K^*_{\mu,\alpha}\mid
\zeta<\mu^+\rangle$ such that 
 $(\bar M,\bar a,\bar
N)^{\zeta+1}$ witnesses that
$(\bar M,\bar a,\bar N)^\zeta$ is not reduced

The construction:
Since $(\bar M,\bar a,\bar N)$ is nice, we can apply Corollary
\ref{c-ext of nice} to find $(\bar M,\bar a,\bar N)^0$ a
$<^c_{\mu,\alpha}$ extension of $(\bar M,\bar a,\bar N)$.  

Suppose that $(\bar M,\bar a,\bar N)^\zeta$ has been defined.  Since it
is a $<^c_{\mu,\alpha}$-extension of $(\bar M,\bar a,\bar N)$, we know
it is not reduced. Let
$(\bar M,\bar a,\bar N)^{\zeta+1}\in\sq{+}\K^*_{\mu,\alpha}$ be a
$\leq^c_{\mu,\alpha}$-extension of
$(\bar M,\bar a,\bar N)^\zeta$, witnessing this.

For $\zeta$ a limit ordinal, let $(\bar M,\bar a,\bar
N)^\zeta=\Union_{\gamma<\zeta}(\bar M,\bar a,\bar N)^\gamma$. This
completes the construction.

For each $b\in\Union_{\zeta<\mu^+,i<\alpha}M^\zeta_i$ define 
$$i(b):=\min\big\{i<\alpha\mid
b\in\Union_{\zeta<\mu^+}\Union_{j<i}
M^\zeta_j\big\}\text{ and}$$
$$\zeta(b):=\min\big\{\zeta<\mu^+\mid
b\in M^\zeta_{i(b)}\big\}.$$

$\zeta(\cdot)$ can be viewed as a function from $\mu^+$ to $\mu^+$.  Thus
there exists a club $E=\{\delta<\mu^+\mid \forall b\in
\Union_{i<\alpha}M^\delta_i,\; \zeta(b)<\delta\}$.  Actually, all we
need is for $E$ to be non-empty.

Fix $\delta\in E$.  By construction
$(\bar M,\bar
a,\bar N)^{\delta+1}$
witnesses the fact that $(\bar M,\bar a,\bar
N)^\delta$ is not reduced.  
So we may fix $i<\alpha$ and $b\in
M^{\delta+1}_i\cap\Union_{j<\alpha}M^{\delta}_j$ such that $b\notin
M^\delta_i$.  Since $b\in M^{\delta+1}_i$, we have that $i(b)\leq i$.
Since $\delta\in E$, we know that there exists $\zeta<\delta$ such that
$b\in M^{\zeta}_{i(b)}$.  Because $\zeta<\delta$ and $i(b)<i$, we have
that
$b\in M^\delta_i$ as well.  This contradicts our choice of
$i$ and $b$ witnessing the failure of $(\bar M,\bar a,\bar
N)^\delta$ to be reduced.
\end{proof}

A variation of the following theorem was claimed in \cite{ShVi 635} for
reduced towers.  Unfortunately, their proof does not converge.  
Under Hypothesis 3, we
resolve their problems here.  

\begin{theorem}[Reduced towers are
continuous]\label{reduced are cont}\index{reduced towers!are
continuous}Under Hypothesis 3, if $(\bar M,\bar a,\bar
N)\in\sq{+}\K^*_{\mu,\alpha}$ is reduced, then it is continuous.

\end{theorem}

The keys to resolving problems of
\cite{ShVi 635} are  the extra conditions in the main
construction and the following lemma which is a consequence of Theorem
\ref{partial ext} and the definition of reduced tower. 

\begin{lemma}\label{monotonicity of C-red}
Suppose that $(\bar M,\bar a,\bar N)\in\sq{+}\K^*_{\mu,\alpha}$ is
reduced and nice, then for every $\beta<\alpha$, $(\bar M,\bar
a,\bar N)\restriction\beta$ is reduced.

\end{lemma}

Notice that without the full $<^c_{\mu,\alpha}$-extension property, it is
conceivable to have a discontinuous reduced tower with non-reduced
restrictions.

\begin{proof}[Proof of Theorem \ref{reduced are cont}]

Suppose the claim fails for $\mu$
and
$\delta$ is the minimal limit ordinal for which it fails.  More
precisely, $\delta$ is the minimal element of 
$$S=\left\{\delta<\mu^+\begin{array}{l|l}
&\delta\text{
is a limit ordinal such that there exists}\\
&\text{an }\alpha<\mu^+\text{ and}\\
&\text {a nice, reduced tower }(\bar M,\bar a,\bar
N)\in\sq{+}\K^*_{\mu,\alpha}\\ 
&
\text{with }M_\delta\succneqq_{\K}\Union_{i<\delta}M_i
\end{array}\right\}$$

 Let $\alpha$ 
witness that $\delta\in S$. Hypothesis 3 implies that every tower is
amalgamable.  Thus we can apply Lemma
\ref{monotonicity of C-red}, to assume that $\alpha=\delta+1$. 
Fix
$(\bar M,\bar a,\bar N)\in\sq{+}\K^*_{\mu,\delta+1}$ 
witnessing that $\delta\in S$. 
 Let $b\in
M_\delta\backslash\Union_{i<\delta}M_i$ be given. 
By Fact \ref{density of
reduced}, Hypothesis 3
 and the minimality of
$\delta$, every nice tower of length $\delta$ has a continuous extension.
Combining this with the fact that $(\bar M,\bar a,\bar
N)\restriction\delta$ is amalgamable, we can apply Proposition \ref{b in
prop} to
$(\bar M,\bar a,\bar N)\restriction\delta$ and $b$ to find a
$<^c_{\mu,\delta}$-extension of $(\bar M,\bar a,\bar
N)\restriction\delta$, say $(\bar M',\bar a\restriction\delta,\bar
N\restriction\delta)\in\sq{+}\K^*_{\mu,\delta}$, in
$\C$ containing
$b$.  Let $M'_\delta\prec_{\K}\C$ be a limit model universal over
$M_\delta$ containing $\Union_{i<\delta}M'_i$.  Notice that $(\bar
M',\bar a,\bar N)\in\sq{+}{\K^*_{\mu,\delta+1}}$ is an extension of $(\bar
M,\bar a,\bar N)$ witnessing that $(\bar M,\bar a,\bar N)$ is not reduced.

\end{proof}

Positive solutions to the following questions would allow us to adjust
the previous proof to conclude that every nice tower has a continuous
extension without any extra hypothesis.

\begin{question}
Is it possible to remove Hypothesis 3 in the proof of
Theorem \ref{reduced are cont}?  Alternatively, can one show the density
of nice,
reduced towers?
\end{question}

The next step towards Shelah's Categoricity Conjecture is to show that
the uniqueness of limit models implies the amalgamation property in this
context.

\end{document}